\documentclass[11pt,reqno]{amsart}

 \makeatletter
\@namedef{subjclassname@2020}{%
  \textup{2020} Mathematics Subject Classification}
\makeatother

 \usepackage{curves} 
\usepackage{graphicx} 
\usepackage{amsmath, amsfonts, amssymb, amscd,graphics, graphpap, tikz, color, xcolor, cancel,enumerate,enumitem, hyperref, mathabx, mathrsfs, mathtools}
\usepackage[scr=rsfs,cal=boondox]{mathalfa}

\usetikzlibrary{arrows.meta,arrows}
\usetikzlibrary{decorations.pathreplacing,decorations.markings}
\tikzset{arrow data/.style 2 args={%
      decoration={%
         markings,
         mark=at position #1 with \arrow{#2}},
         postaction=decorate}
      }%
      
\usepackage[mathscr]{euscript}
  \usepackage[normalem]{ulem}
  \usepackage{soul}

\usepackage{curves} 
\usepackage{graphicx}

\makeatletter
\newcommand{\doublewidetilde}[1]{{%
  \mathpalette\double@widetilde{#1}%
}}

\textwidth= 
148mm
\numberwithin{equation}{section}

\theoremstyle{plain}
\newtheorem{theo}{Theorem}[section]
\newtheorem{lem}[theo]{Lemma}
\newtheorem{prop}[theo]{Proposition}

\newtheorem{cor}[theo]{Corollary}

\theoremstyle{definition}
\newtheorem{rem}[theo]{Remark}

\newtheorem{example}[theo]{Example}
\newtheorem{definition}[theo]{Definition}

\newenvironment{pf}{\noindent{\it Proof.\,}}{\hfill $\square$\par \medskip}

\theoremstyle{plain}

\theoremstyle{definition}

\newcommand{\rank}{\operatorname{rank}}
\newcommand{\beq}{\begin{equation}}
\newcommand{\eeq}{\end{equation}}
\renewcommand{\a}{\alpha}
\renewcommand{\b}{\beta}

\renewcommand{\d}{\delta}

\newcommand{\f}{\varphi}
\newcommand{\g}{\gamma}
\newcommand{\h}{\eta}

\renewcommand{\l}{\lambda}
\renewcommand{\o}{\omega}

\renewcommand{\r}{\rho}
\newcommand{\s}{\sigma}
\renewcommand{\t}{\tau}

\newcommand{\z}{\zeta}

\newcommand{\bC}{\mathbb{C}}
\newcommand{\bR}{\mathbb{R}}

\newcommand{\bF}{\mathbb{F}}

\newcommand{\bN}{\mathbb{N}}

\newcommand{\bW}{\mathbb{W}}

\newcommand{\bT}{\mathbb{T}}
\newcommand{\bV}{\mathbb{V}}
\newcommand{\bY}{\mathbb{Y}}

\newcommand{\gq}{\mathfrak{q}}

\newcommand{\gr}{\mathfrak{r}}

\newcommand{\gy}{\mathfrak{y}}

\newcommand{\gQ}{\mathfrak Q}

\newcommand{\cA}{\mathscr{A}}
\newcommand{\cB}{\mathcal B}
\newcommand{\cC}{\mathcal{C}}
\newcommand{\cD}{\mathscr{D}}
\newcommand{\cE}{\mathscr{E}}

\newcommand{\cG}{\mathscr{G}}

\newcommand{\cK}{\mathscr{K}}
\newcommand{\cL}{\mathscr{L}}
\newcommand{\cM}{\mathscr{M}}
\newcommand{\cN}{\mathscr{N}}

\newcommand{\cQ}{\mathscr{Q}}
\newcommand{\cR}{\mathscr{R}}
\newcommand{\cS}{\mathscr{S}}
\newcommand{\cT}{\mathscr{T}}
\newcommand{\cU}{\mathscr{U}}
\newcommand{\cV}{\mathscr{V}}

\newcommand{\p}{\partial}

\renewcommand{\square}{\kern1pt\vbox
{\hrule height 0.6pt\hbox{\vrule width 0.6pt\hskip 3pt
\vbox{\vskip 6pt}\hskip 3pt\vrule width 0.6pt}\hrule height0.6pt}\kern1pt}

\DeclareMathOperator{\Span}{Span}
\DeclareMathOperator{\ad}{ad}

\newcommand{\wt}{\widetilde}

\newcommand{\wh}{\widehat}
\newcommand{\wc}{\widecheck}

\newcommand{\bt}{\begin{theo}\ \ }
\newcommand{\et}{\end{theo}}
\newcommand{\bp}{\begin{prop}\ \ }
\newcommand{\ep}{\end{prop}}
\newcommand{\bc}{\begin{cor}\ \ }
\newcommand{\ec}{\end{cor}}
\newcommand{\bl}{\begin{lem}\ \ }
\newcommand{\el}{\end{lem}}
\newcommand{\bd}{\begin{definition}}
\newcommand{\ed}{\end{definition}}
\newcommand{\be}{\begin{equation}}
\newcommand{\ee}{\end{equation}}

\def\<#1,#2>{\langle\,#1,\,#2\,\rangle}
\newcommand{\arr}{\begin{array}{rlll}}
\newcommand{\ea}{\end{array}}
\newcommand{\bea}{\begin{eqnarray}}
\newcommand{\eea}{\end{eqnarray}}
\newcommand{\bean}{\begin{eqnarray*}}
\newcommand{\eean}{\end{eqnarray*}}

\renewcommand{\=}{:=}

\newcommand{\ve}{\varepsilon}

  \newcommand{\vertiii}[1]{{\left\vert\kern-0.25ex\left\vert\kern-0.25ex\left\vert #1 
    \right\vert\kern-0.25ex\right\vert\kern-0.25ex\right\vert}}

\newcommand{\Surr}{\cS\text{\it urr}^{(\cU, \bT,T)}}
\newcommand{\Attain}{\cM\text{\it -Att}}

\newcommand{\Reach}{\cR\text{\it each}}

\def\sideremark#1{\ifvmode\leavevmode\fi\vadjust{
\vbox to0pt{\hbox to 0pt{\hskip\hsize\hskip1em
\vbox{\hsize3cm\tiny\raggedright\pretolerance10000
\noindent #1\hfill}\hss}\vbox to8pt{\vfil}\vss}}}

\title[Distributions and controllability  problems (I)]{Distributions and controllability  problems (I)}
 \author[Cristina Giannotti,  Andrea Spiro and Marta Zoppello]{Cristina Giannotti \quad  Andrea Spiro \quad Marta Zoppello}
 \subjclass[2020]{93B03, 93B05, 34H05}
 \thanks{{\it Data availability statement.} 
No datasets were generated or analysed during the current study.}
 \keywords{Controllability of non-linear systems; Chow-Rashevski\u\i\ Theorem; Kalman criterion; Reachable sets}
\begin{document}

\begin{abstract} 
 We consider a   non-linear  real analytic control   system   of first order $\dot q^i = f^i(t, q, w)$, with controls $w = (w^\a)$ in  a connected open set  $\cK \subset \bR^m$  and  configurations $q = (q^i)$ in  $\cQ \= \bR^n$. The set of points in the extended   space-time  $\cM = \bR \times \cQ \times \cK$, which  can be reached  from  a triple  $x_o = (t_o , q_o, w_o) \in \cM$  through  a continuous graph completion    $\g(s) = \big(t(s), q(t(s)), w(t(s))\big)$ of  the graph of   a solution   $t \to (q(t), w(t))$, $t \in [t_o ,t_o +  T]$,   with   piecewise real analytic  controls,  is called  the  {\it $\cM$-attainable set  of $x_o$ in time $T$}.  We   prove that  if   $y_o$ is an $\cM$-attainable point of $x_o$, a large set of  other nearby  $\cM$-attainable points of $x_o$ can be  determined  starting   directly   from $y_o$  and applying  an  appropriate ordered composition of  flows  of  vector fields in a distinguished  distribution $\cD^{II} \subset T \cM$,   canonically  associated with   the control system.  We  then determine   sufficient conditions  for such  neighbouring points  to constitute  an   orbit of  the pseudogroup of local diffeomorphisms generated by  the vector fields in $\cD^{II}$. If   such   conditions are satisfied and if  the tangent spaces of these  orbits  have maximal rank  projections  onto $\cQ$,  the control system  is locally accessible and  has the   small time local controllability property  near the state   points of equilibrium.  These results   lead to  new proofs of     classical local 
 controllability  criterions   and  yield  new  methods to establish the    accessibility  and  the small time local   controllability  of  non-linear control systems. 
\end{abstract}
\maketitle

\setcounter{section}{0}
\section{Introduction}
Investigating the controllability of a non-linear control system is often a quite hard task. And  most of the known criterions for the  accessibility  or the  small time controllability of a non-linear  system (as e.g.  the linear Kalman test, the Chow theorem for driftless systems, the  Sussmann criterion, the Coron return method, etc.) are developed only for systems  that are  affine  in  the controls. 
In this paper we  tackle the general  problem of  the accessibility and  the small time local  controllability of non-linear  real analytic  systems,  not necessarily  with  affine  controls.  Our approach  starts with a discussion   of  the points in the extended  space $\cM = \bR \times \cQ \times \cK$, given by  the triples $(t,q,w)$ 
of times, states and controls, which can be reached from a triple $x_o = (t_o, q_o, w_o)$  through graph completions of  graphs of solutions determined by piecewise real analytic  controls.  Let us call such points {\it $\cM$-attainable from $x_o$}.  Our first main result consists in the proof that, given an $\cM$-attainable  point $y_o$   from a fixed $x_o \in \cM$, a  large set of  points, which are $\cM$-attainable  from $x_o$ and are   in proximity of  $y_o$,  can be determined  in a  direct way  applying to  $y_o$    an    appropriate ordered compositions of  flows of  certain  vector  fields, which we call   {\it surrogate vector  fields}. These  vector   fields   take values in   a  particular distribution $\cD^{II} \subset T \cM$, called {\it secondary distribution}, which is canonically associated with the control system  and is very easy to be   determined.  We  then   prove that  near any $\cM$-attainable  point $y_o$, there exists a set of  generators for the secondary distribution, which consists  only of surrogate  vector fields. Combining these  two fact,  we are able to show   that  a very large set of   $\cM$-attainable points in a neighbourhood       of    $y_o$  fill  a  (possibly, proper)  subset of the (local) orbit of $y_o$   under the action of  the pseudogroup of local diffeomorphisms generated  by   the vector fields in $\cD^{II}$.   Using  a modification of  Rashevski\u\i's proof of the  celebrated  Chow-Rashevski\u\i-Sussmann Theorem \cite{Ra, GSZ1}, we establish a couple of   sufficient conditions  for   the above described set of  $\cM$-attainable points from  $x_o$  to  coincide with  the local orbit of $y_o$ under the action of  such pseudogroup. If these conditions are satisfied,  this set of neighbouring $\cM$-attainable  points 
is an open   subset of  a maximal  integral leaf of the  involutive distribution $\cE^{II\operatorname{(Lie)}} \subset T\cM$  generated by the iterated Lie brackets of the vector   fields in $\cD^{II}$.  It follows that, under these sufficient conditions, whenever the maximal  integral leaves of   $\cE^{II\operatorname{(Lie)}} $ have  maximal rank projections onto   $\cQ$,   the control system is accessible and has the small time local controllability property near its stable points. In this way  we  get not only   new proofs of    classical local
 controllability  criterions, but new criterions which can be used 
to establish the  accessibility and the  small time local   controllability of certain non-linear control systems, for which, at the best of our knowledge, all so far known criterions are inconclusive. \par
 \smallskip
We now go into a more detailed  description of our results. Consider a real analytic  controlled system in normal form 
\beq \label{controlsystem-intro} 
\dot q^i(t) =  f^i\big(t, q(t), w(t)\big)\ , \qquad 1 \leq i \leq n\ ,\eeq
for curves $q(t)$ in the  space of states  $\cQ = \bR^n$ and  control curves $w(t) = (w^\a(t))$ in a region $\cK$ of   some $\bR^m$.   Let $\bT$, $\cD$ and $\cD^I \subset \cD$ be the vector field and the two constant rank distributions  on $\cM =  \bR \times \cQ \times \cK$ given by 
\beq \label{primary} 
 \bT= \frac{\p}{\p t} + f^i(t, q, w) \frac{\p}{\p q^i} \ ,
\ \   \cD = \left\langle \bT,  \frac{\p}{\p w^1},  \ldots, \frac{\p}{\p w^m}\right\rangle\ ,\ \   \cD^I = \left\langle \frac{\p}{\p w^1},  \ldots, \frac{\p}{\p w^m}\right\rangle
\eeq
and  denote by  $\cD^{II}$  the (possibly singular) distribution,  generated by the vector fields
\beq \label{gen11} \frac{\p}{\p w^\a} \qquad \text{and}\qquad \underset{\text{$k$-times}}{\underbrace{\left[\bT, \left[\bT, \ldots, \left[ \bT, \frac{\p}{\p w^\a}\right]\ldots\right]\right]}}\ ,\qquad  \ 1 \leq k < \infty\  ,\quad  1 \leq \a \leq m\ .\eeq
We call  $\cD$ and  $\cD^{II}$ the {\it primary} and the  {\it secondary distribution}, respectively,   of  \eqref{controlsystem-intro}. \par
Given   a  point $x_o = (t_o , q_o, w_o) \in \cM$ and $T> 0$,  we   denote by  $\operatorname{Attain}^{(T)}(x_o)$   the set of all final points of   graph completions $\g(s) = \big(t(s), q(t(s)), w(t(s))\big) $, corresponding to  graphs of  solutions  $t \to (q(t), w(t))$, $t \in [t_o ,t_o +  T]$,  with initial data $(q(t_o), w(t_o)) = (q_o, w_o)$  and    piecewise   real analytic  controls $w(t)$. 
Note  that, in case  $t_o = 0$,  the projection of  the set $\operatorname{Attain}^{(T)}(x_o)$  onto  $\cQ$  is nothing but  the set of   reachable points from $q_o$,  which are determined by piecewise real analytic   controls and in   time exactly equal to $T$. Due to this, via  projections onto $\cQ$,  information on  the  sets $\operatorname{Attain}^{(T)}(x_o) \subset \cM $    provide  important information on the  reachable sets of the control system. \par
\smallskip
Our first main  result is given by the following  theorem. Here, given a   local vector field $X$ on  a relatively compact  open set $\cU \subset \cM$, we denote by   $\Phi^X_s: \cU \to \cM$, $s \in (- \ve , \ve)$,    a one-parameter family of   local diffeomorphisms given   by  the flow of  $X$. \par
\noindent{\bf Theorem A.}   {\it Let $T> 0$, $x_o \in \cM$ and $y_o$ a point in  $\operatorname{Attain}^{(T)}(x_o)$, $T >0$, and  denote by $\cU \subset \cM$ a relatively compact neighbourhood  of $y_o$. Let also $\d T \in (0, T)$ small enough, so that the points with  times $t \in [T-\d T, T]$ of the  graph of a  piecewise $\cC^\o$ solution, that  starts from $x_o$ and ends in $y_o$,  constitute the trace of a real analytic curve.  If $\d T$ and  $\cU$  are sufficiently small, then: 
\begin{itemize}[leftmargin = 20pt]
\item[(1)] There are  intervals $(-\ve_j, \ve_j) \subset  \bR$ such that  $\operatorname{Attain}^{(T)}(x_o)$ contains all points $y \in \cU$  which   can be obtained from $y_o$ by the expression 
\beq
\label{surrogat}
  y= \Phi^{Y^{(\t_p)}_p}_{s_p}\hskip  -0.5 cm  \circ \ldots \circ \Phi^{Y^{(\t_1)}_1}_{s_1}(y_o)\eeq
for some  $ s_j \in (-\ve_j, \ve_j)$ and  vector fields having the form $Y^{(\t_j)}_j \= \Phi^{\f_*(\bT)}_{\t_j*}\left( \l_j^\a \frac{\p}{\p w^\a}\right)$, where    $\lambda_j^\a \in \bR$,  $\f$ is a real analytic diffeomorphism, mapping  the integral leaves of $\cD^I$ into themselves,    and the $\t_j$'s are real number     satisfying the inequalities: 
\beq \label{constraint1}  \d T > \t_1 > \t_2  >  \ldots >\t_p > 0\ . \eeq
\item[(2)] Any vector field $Y^{(\t_j)}_j \= \Phi^{\f_*(\bT)}_{\t_j*}\left( \l_j^\a \frac{\p}{\p w^\a}\right)$, $\l^\a_j \in \bR$, $\t_j \in (0, T)$,  takes values in $\cD^{II}|_{\cU}$; 
\item[(3)]  There is a set of generators for  $\cD^{II}|_{\cU}$ which is made of vector fields  of the form  (2). 
\end{itemize}
}
The  vector fields in   (2) are the {\it surrogate fields}, that we  mentioned at the beginning of this introduction.  We stress the fact that Theorem A implies that  for any sufficiently small neighbourhood $\cU$ of a   point $y_o \in \operatorname{Attain}^{(T)}(x_o)$,   the set  $ \operatorname{Attain}^{(T)}(x_o) \cap \cU$
contains a subset of points -- let us denote it by  $\wt{ \operatorname{Attain}^{(T)}(x_o)}\cap \cU$  --belonging to the  orbit $\operatorname{Orb}^{\cG}(y_o)$   of   $y_o$ of the local diffeomorphisms  in the  pseudo-group $\cG$ generated by the  flows of the  vector fields in $\cD^{II}$.  
Note also that  the proof of Theorem A gives evidences that very likely the set $\wt{ \operatorname{Attain}^{(T)}(x_o)}\cap \cU$
is essentially equal to $\operatorname{Attain}^{(T)}(x_o)\cap \cU$ for sufficiently small neighbourhoods.
Thus,  {\it  whenever this  subset $ \wt{\operatorname{Attain}^{(T)}(x_o) } \cap \cU$  coincides with   $\operatorname{Orb}^{\cG}(y_o)\cap \cU $}, 
several  useful information on  the local structure  of $\operatorname{Attain}^{(T)}(x_o)$    are immediately derivable from information on  the pseudogroup  orbits   $\operatorname{Orb}^{\cG}(y_o)\cap \cU $ (the latter being  much  easier  to be determined).  \par
\smallskip
Unfortunately,  in general, one has that  $ \wt{\operatorname{Attain}^{(T)}(x_o)} \cap \cU \subsetneq \operatorname{Orb}^{\cG}(y_o) \cap \cU$, regardless on the size of the neighbourhood $\cU$. The  main reason for not   having the equality between the two sets is    the fact that    the parameters  $\t_j$, which occur in the definition of the surrogate fields,   are    constrained by the inequalities  \eqref{constraint1} (indeed,  if there were not   such a constraint,   the equality could   be very easily established).   The   points $y_o \in  \operatorname{Attain}^{(T)}(x_o)$  for which  there is a    neighbourhood $\cU$  with the property that    $\wt{\operatorname{Attain}^{(T)}(x_o) } \cap \cU = \operatorname{Orb}^{\cG}(y_o) \cap \cU$   are named {\it good points}.
\par
\smallskip
We are now facing  two crucial  problems: (a) Find sufficient conditions for a  point $y_o$  to be a   good point; (b) Determine the  orbits $\operatorname{Orb}^{\cG}(y_o) \cap \cU$ under the pseudo-group $\cG$  generated by the  local flows of the  vector fields in $\cD^{II}$. \par
\smallskip
At the best of our knowledge,  the most general method to tackle a problem  as in  (b)  is  provided by   Sussmann's   results in  \cite{Su}.   According to them,  the orbits $\operatorname{Orb}^{\cG}(y_o)$  coincide with the maximal integral leaves of the smallest    distribution   $\cE^{II(\text{Suss)}} $ which contains $\cD^{II}$ and  is invariant under the flows of the vector fields in $\cD^{II}$.  Since we deal with real analytic vector fields. such a distribution $\cE^{II(\text{Suss)}} $ coincides with the  involutive distribution   $\cE^{II(\text{Lie})} \supset \cD^{II}$,  which is generated  by  all  linear combinations and  iterated Lie brackets of the local vector fields in $\cD^{II}$ (see \cite[Thm. 5.16 \& Cor. 5.17]{AS}).  These results essentially gives the way to  answer   the problem (b). For what concerns the problem (a), the same facts  together with an argument taken from the proof of the Chow-Rashevsk\u\i i Theorem (see \cite{GSZ1}) led us to the second main result of this paper. In order to  state it in a simpler way, it is convenient to preceed it by  the  following notion.  \par
\smallskip  Consider  a   distribution $\wc \cD$ on a manifold $\cM$ and the (generalised) involutive distribution  $\wc \cE^{(\text{Lie})}$ spanned at each point $y \in \cM$ by  the values at $y$ of all finite linear combinations   of the  vector fields in $\wc \cD$ and all possible their   iterated Lie brackets. Given an open subset  $\cU \subset \cM$,  we call  {\it decomposition of  $\cU$  into  $\wc \cD$-strata of  maximal  $\wc \cD$-depth    $\mu$}  any expression  of $\cU$ as  a  finite union  of  disjoint subsets  $\cU = \cU_0 \dot \cup \cU_1 \dot \cup  \ldots  \dot \cup \cU_p$   such that:  
\begin{itemize}[leftmargin = 20pt]
\item[(i)] for a fixed  $0 \leq j \leq p$,    all  spaces   $\wc \cE^{(\text{Lie})}|_{y} \subset T_y\cM$, $y \in \cU_j$,  have   the  same dimension  and the   integral leaves  of $\wc \cE^{(\text{Lie})}|_{\cU}$ through the points  $y \in \cU_j$ are     contained  in $\cU_j$; 
\item[(ii)] there are  integers $1 \leq \mu_j \leq \mu$,   $0 \leq j \leq p$,  (called {\it $\wc \cD$-depths} -- one per each $\cU_j$) such that  for any  $y \in \cU_j$ the  space   $\wc \cE^{(\text{Lie})}|_{y}$ is spanned  by  vectors of the form  $[Y_{i_1}, [Y_{i_2}, [\ldots [Y_{i_{r-1}}, Y_{i_r}]\ldots ]]]|_y$,  where each   $Y_{j_\ell}$ is in $\wc \cD|_{\cU}$ and  the integer   $r$ is less than or equal to $\mu_j (\leq \mu)$.   
\end{itemize}
As a direct consequence of  the Noetherianity of the  rings of  real analytic functions (\cite[Thm. 3.8]{Ma}), one can prove that  for any   $x_o \in \cM$ there is  a neighbourhood $\cU$  admitting a decomposition into  $\wc \cD$-strata of an appropriate maximal depth.  We  say that 
a real analytic control system \eqref{controlsystem-intro} is   {\it of   type $\mu$  on some  $\cU \subset \cM = \bR \times \cQ \times \cK$} if $\cU$ admits a decomposition into $\cD^{II}$-strata  of maximal $\cD^{II}$-depth  $\mu$.
\par
\smallskip
We can now state our  second main result, which gives a pair of sufficient conditions for the goodness of reachable points of  systems of  this  kind. \par
\noindent{\bf Theorem B.}   {\it  Let  $y_o$ be a point in $\operatorname{Attain}^{(T)}(x_o)$, $T >0$, and assume that the system is of  type $\mu= 1$ or $2$   on a neighbourhood $\cU$ of $y_o$. Then $y_o$ is a good point if: \\[-20pt]
\begin{itemize}[leftmargin = 20pt]
\item[($\a$)] Either  $\mu = 1$ or  
\item[($\b$)]  $\mu = 2$ and there exists  a  set of  generators for the distribution  $\cE^{II(\text{\rm Lie})}|_{\cU}$  satisfying the  conditions of  Theorem \ref{criterione} below (see \S \ref{sect8.3} for details) . 
\end{itemize}
}\par
\smallskip
We stress the fact  that the conditions  stated in ($\a$) and ($\b$)  are just  conditions -- easy to be checked -- on  the Lie brackets  between the generators \eqref{gen11}  for   $\cD^{II}$. Moreover,  our proof of Theorem B  is designed to 
allow  generalisations to the cases in which the control system  is locally of any type  $\mu\geq 3$.  Such  generalisations are left to  future work.  But it is   remarkable that,  as we mentioned above,     Theorems A and  B are good enough  to get   new proofs of the  classical Kalman criterion and  Kalman linear test together with new methods to establish small time local controllability of certain systems, for which at the best of our knowledge all  so-far known   criterions are inconclusive.   A  short overview of such applications   is given in the  concluding section  \S \ref{section8}.  The details  are postponed to   \cite{GSZ2}, which is the natural  continuation of this paper.    
\par
As a concluding remark,  we  would like to point out  that, despite of the fact that  our results concern only real analytic  control systems,  it is reasonable to expect  that,  by means of  approximation techniques,  several  parts of our methods and  results can be extended  to  a large class of control systems of  class $\cC^k$ for   large  $k$.     We hope to address   this issue in the near future. \par
\smallskip
 The paper is divided in three parts. Part I starts with  a  preliminary section on oriented curves  and distributions.   We then introduce the notion of  {\it completed graphs}  of  solutions   with  piecewise real analytic controls  and  we geometrically characterise them as curves tangent to a so-called  {\it rigged distribution} -- see \S \ref{thirdsection}. In \S \ref{section4} we  introduce the  {\it sets of $\cM$-attainable points}  and we briefly explain their relations with the reachable sets of  the considered   control system.  Part II consists of just the section      \S \ref{thefifthsection},  where after introducing    the notions of   {\it secondary distribution}  and  of  {\it $\bT$-surrogate vector  fields},  we  prove  a few   properties, of which Theorem A is a direct consequence. Part III begins with the section  \S \ref{sect8},  where  we prove the existence of sets of local generators for the secondary distribution, which are made just  of $\bT$-surrogate fields,  and it continues with   \S \ref{theeigthsection}, where we introduce the notion of {\it good points} and prove the two criterions, which correspond to the conditions ($\a$) and ($\b$)  of  Theorem B. In  \S 8  we provide the above mentioned   short survey of  applications of our main results, referring to  \cite{GSZ2} for details. We also  state a few  conjectures and open problems. The paper ends with a couple of  appendices, where the proofs of two technical lemmas are given. \par
\bigskip
\addcontentsline{toc}{section}{Part I}
\centerline{\large \it PART I}
\ \\[-40pt]
\section{Oriented curves tangent to  distributions} 
\subsection{Piecewise regular oriented curves} Given an  $N$-dimensional manifold $\cM$,  a  {\it curve in $\cM$}  is    the trace $\g([a, b])$   of a differentiable and regular parameterised curve $\g(t)$, i.e. of a map $\g:[a,b] \subset \bR  \longrightarrow \cM$ of class $\cC^1$  and with nowhere vanishing velocity $\dot \g(t) \neq 0$. 
 Two differentiable and regular parameterised curves  $\g(t)$ and $\wt \g(s)$ with same trace  are called  {\it consonant} if  one is obtained  from the other by a   change of  parameter $t = t(s)$ with 
   $\frac{dt}{ds} > 0$ at all points. Consonance  is an equivalence relation  and an {\it orientation} of a  curve  is a choice of one of the    two   equivalence classes of its parameterisations. 
 An {\it oriented curve} is a curve  with an orientation. We   indicate the orientation  by  one of its consonant   parameterisations $\g(t)$. \par
\smallskip
Let $\g_1(t)$, $\g_2(s)$, $t \in [a, b]$,  $s \in [c,d]$, be  two  (oriented)  curves, such that the final endpoint   $\g_1(b)$ of the first curve is equal to the 
initial  endpoint  $\g_2(c)$ of the second curve.  The  {\it (oriented) composition} $\g_1\ast \g_2$ is the  union of  the  two (oriented) curves determined by the two parameterisations. The curves  $\g_1$, $\g_2$ are called {\it regular arcs} of $\g_1 \ast \g_2$.  We define in a very similar way  the  {\it (oriented) composition} of a finite number of (oriented) curves $\g_1$, $\g_2$, $\ldots$, $\g_r$,  each of them sharing its final endpoint with the initial endpoint of the succeeding one. A connected subset of $\cM$, which is obtained as (oriented) composition of a finite collection of (oriented) curves,  is called {\it piecewise regular (oriented) curve}.\par
\smallskip
\subsection{Regular and singular distributions and their tangent curves}  
 A {\it regular distribution of rank $p$}  on  $\cM$  is a smooth family $\cD$  of subspaces $\cD_x \subset T_x \cM$ of the tangent spaces of $\cM$ of constant dimension $p$. Here,  with the term ``smooth family''   we   mean  that  for any  point $x_o \in \cM$ 
there is a neighbourhood $\cU$ and a  set of $\cC^\infty$ pointwise linearly independent vector fields $X_1, \ldots, X_p$ on $\cU$,  such that 
$$\cD_x = \langle X_1|_x, \ldots, X_p|_x \rangle\qquad \text{for any}\ x \in \cU\ .$$
Given a regular  distribution $\cD \subset T\cM$ and a (local) vector field $X$,  we are going to   use the  notation ``$X \in \cD$'' to indicate that $X_x$  is in  $\cD_x$  at any point $x$ where the field $X$ is defined. \par
\smallskip
Generalisations of the notion of ``regular distribution'',  in which the  condition $\dim \cD_x = \text{const.}$  is not assumed, are possible,  but  demand some  care.  In this paper we   adopt  the  following definition,  which is a  variant of   those  considered in   \cite{Na, SJ, Su}. 
\begin{definition}\label{quasi-regular-def}  A {\it  quasi-regular set of   $\cC^\infty$ (resp.  $\cC^\o$) vector fields of rank $p$\/} on $\cM$ is a set $V$ of local vector fields of the following kind. There exist  
  an open cover  $\{\cU_A\}_{A \in J}$  of $\cM$ and a family    $\{ (X^{(A)}_{1}, \ldots, X^{(A)}_{p_A})\}_{A \in J}$  of tuples of cardinalities $p_A \geq p$,   of $\cC^\infty$ (resp. $\cC^\o$) vector fields  -- one tuple for each   open set $\cU_A$ --  each of them containing a $p$-tuple, made of vector fields that are 
 pointwise linearly independent on some open and dense subset
of $\cU_A$,  and  satisfying the following  conditions: 
\begin{itemize}[leftmargin = 10pt]
\item  if $\cU_A \cap \cU_B \neq \emptyset$,  then  for any $x \in \cU_A \cap \cU_B$   one has $ X^{(A)}_{i}|_x = \cA^{(AB)}{}_{i}^j\big|_x X^{(B)}_{j}\big|_x$ for a $\cC^\infty$ (resp. $\cC^\o$) matrix valued map
$\cA^{(AB)} : \cU_{A} \cap \cU_B \to \bR_{p_A \times p_B}$
\item the vector fields $Y \in V$  are exactly the local vector fields,   for which  any  restriction  $Y|_{\cU \cap \cU_A}$ to the intersection  between the   domain $\cU$  and a set  in  $\{\cU_A\}_{A \in J}$, has   the form
\beq Y|_{\cU_A \cap \cU} = Y^{(A)i} X^{(A)}_{i}\eeq
for  some $\cC^\infty$ (resp. $\cC^\o$)  functions $Y^{(A)i}$.  
\end{itemize}
The tuples $(X^{(A)}_{1}, \ldots, X^{(A)}_{p_A})$  are   called {\it sets of local  generators for $V$}.\par 
 A 
  {\it  smooth (resp. real analytic) generalised distribution of rank $p$}   is a pair $(V, \cD^V)$ given by  a quasi-regular set of $\cC^\infty$ ($\cC^\o$) vector fields    $V$ of rank $p$  and  
the  associated family $\cD^V$ of  tangent  subspaces, defined   by 
  $\cD^V_x  = \{X_x\ ,\ X \in V \}$, $x \in \cM$.
 If $\dim \cD^V_x $  is constant,  $(V, \cD^V)$ is called  {\it regular}, otherwise it is called {\it singular}.
\end{definition}
  \par
\smallskip
 Note that if   $(V, \cD^V)$ is   regular, the set  $V$ coincides with the full set of local vector fields with  values in  $ \cD^V$ and  the pair $(V, \cD^V)$ is  fully determined by  $ \cD^V$.  On the contrary,  if  $(V, \cD^V)$ is   singular, the set $V$ is no longer determined by $\cD^V$, because  there might be several different   generalised distributions having the same    family of tangent subspaces $\cD^V$. For instance, 
 the quasi-regular sets of local vector fields on $\bR^2$ 
 \begin{align*}
 &V = \Big\{\ X = f^1(x) \frac{\p}{\p x^1}  + f^2(x) x^1 \frac{\p}{\p x^2}\ ,\ f^i \ \text{smooth}\ \Big\}\ ,\\
 & \wt V = \Big\{\ X = g^1(x) \frac{\p}{\p x^1}  + g^2(x) (x^1)^2 \frac{\p}{\p x^2}\ ,\ g^i \ \text{smooth}\ \Big\}\ 
 \end{align*}
 generate the same family of spaces $\cD^V = \cD^{\wt V}$,   but $(V, \cD^V) \neq (\wt V, \cD^{\wt V})$ . 
\par
A  generalised distribution $(V, \cD^V)$   is called {\it non-integrable} if there exists at least one  pair of vector fields $X$, $Y \in V$, whose  Lie bracket $[X, Y]$  is not  in $V$.  Otherwise  it is called {\it integrable} or {\it involutive}.
  \par
\smallskip
A  curve $\g$  is said to be {\it tangent to the generalised distribution $(V,\cD =  \cD^V)$}  if for one (hence,  for all) regular parameterisation $\g(t)$ of  the curve, the velocities $\dot \g(t)$  are such that  $\dot \g(t) \in \cD_{\g(t)}$  for any $t$.  A piecewise regular curve $\g = \g_1 \ast \g_2 \ast \ldots \ast \g_r$  is    {\it tangent to $(V, \cD = \cD^V)$} if each of its regular arcs has this property. We  call  any such $\g$  a {\it  $\cD$-path}.
 \par
 The equivalence classes in $\cM$  of the relation
 \beq x \simeq x'\qquad \Longleftrightarrow\qquad \text{there exists a $\cD$-path joining $x$ to $x'$}\eeq
 are the  {\it $\cD$-path connected components} of $\cM$. We recall that, by the usual   proof of Frobenius Theorem (see e.g. \cite{Wa}), whenever $\cD$ is regular and involutive, the $\cD$-path connected components coincide with the maximal integral leaves of  the distribution. \par
 \bigskip
\section{First order control  systems and  rigged distributions} \label{thirdsection}
\subsection{First order control systems and completed graphs of solutions}
Consider  a  first order system of control  equations on curves $q(t) = (q^i(t))$ in $\cQ = \bR^n$ of the form 
\beq \label{controlsystem} 
\dot q^i(t) =  f^i(t, q(t), w(t))\ , \qquad 1 \leq i \leq n\eeq
where the  control curves  $w(t) = (w^\a(t))$ take values in an  open and connected subset $\cK$  of $\bR^m$  and the $f^i$  are  smooth real functions $f^i: \bR \times \cQ \times \cK \to \bR$. 
 \par
\smallskip
A  {\it solution} of \eqref{controlsystem} is a  map   $t \mapsto (q(t), w(t)) \in \cQ \times \cK$, $t \in [a, b] \subset \bR$, in which  $q(t)$   is an absolutely continuous   map  with values  in $\cQ$ and   $w(t)$ is a   measurable  map  with values  in $\cK$,   satisfying \eqref{controlsystem} at almost every point.  
Since in this paper  we consider only  maps,   in which $w(t)$ is   piecewise $\cC^1$,   possibly not continuous but     with  only finite  jumps  at the points of discontinuity (and thus with also $q(t)$ 
 piecewise  $\cC^1$),  from now on we tacitly assume that a ``solution'' is a map with  this  additional assumption. Note that, for any solution $\big(t, q(t), w(t)\big)$, the corresponding  {\it  graph} 
 $$\{\big(t, q(t), w(t)\big)\ , t \in [a,b]\} \subset \bR \times \cQ \times \cK\ ,$$  is    union of a finite collection of {\it oriented}  curves $\g_1, \ldots, \g_r$ (with possibly one or two endpoints deleted)  having  the following property:  {\it  the standard projection 
   $$\pi^\bR: \bR \times \cQ \times  \cK \to \bR$$
maps bijectively each   $\g_i$  onto  a non-trivial interval $I$  of $ \bR$  and the orientation of such $\g_i$ is  the  usual one,  corresponding to  increasing times (more precisely, it is the orientation given by   the equivalence class of the natural parameterisation $(\pi^{\bR}|_{\g_i})^{-1}: I \subset \bR \to \g_i$)}.  Combining this fact with  Bressan and  Rampazzo's notion of graph completion \cite{BR},   we  are led  to the following
\begin{definition}  \label{completed} A {\it completed graph} of a solution of \eqref{controlsystem} is a piecewise regular  oriented curve  $\h_1 \ast \h_2 \ast  \ldots\ast \h_r$ in  $\bR \times \cQ \times  \cK$, 
whose  regular arcs $\h_i$ are mapped by   $\pi^\bR$  either onto a non-trivial interval   with the standard orientation or  onto a singleton  $\{t\}$,  and satisfy the following conditions 
\begin{itemize}[leftmargin = 15pt]
\item  if the $\bR$-projection of $\h_i$  is a non-trivial interval, then $\h_i$ is the graph of a smooth solution  $t \mapsto (q(t), w(t))$ of \eqref{controlsystem}; 
\item  the arcs $\h_1$ and $\h_r$  are both mapped onto a non-trivial interval; 
\item  for $2 \leq i \leq r -1$, whenever the $\pi^\bR$-projection of $\h_i$  is a  singleton,   also the image  of  the standard projection   $\pi^{\bR \times \cQ}: \bR \times \cQ \times \cK \to \bR \times \cQ$ is a singleton. 
\end{itemize}
The   regular arcs $\h_i$ which project onto singletons of $\bR$ constitute   the  {\it added part} of the completed graph. Two completed graphs of solutions are said to be {\it g-equivalent} if they differ only for their added parts.
\end{definition}  
One can  directly see  that the graph of {\it any} piecewise smooth   solution $(q(t), w(t))$ of  \eqref{controlsystem} is contained in at least one completed graph.   This and  any other g-equivalent completed graph    are the   {\it graph completions} of $(q(t), w(t))$.   Note  that the projection  onto $\bR \times \cQ$ of any graph completion of  a solution $t \mapsto (q(t), w(t))$   is    always the  graph of the  map $t \mapsto q(t)$ and   is always a piecewise regular oriented curve. \par
\medskip
 \centerline{
\begin{tikzpicture}
  \begin{scope}[white]
        \draw[fill=red!30, semitransparent] (0.5, 0) -- (0.5,2.5) -- (1,4) --  (1,1.5) --cycle;
        \draw[fill=red!30, semitransparent] (2, 0) -- (2,2.5) -- (2.5,4) --  (2.5,1.5) --cycle;
        \draw[fill=red!30, semitransparent] (3.5, 0) -- (3.5,2.5) -- (4,4) --  (4,1.5) --cycle;
\end{scope}
\draw[->, line width = 0.6, blue]  (0,0.5) to (5,0.5);
\draw[->, line width = 0.6, red]  (0.5,0) to (0.5,4);
\draw[->, line width = 0.2, red]  (0.4,0.2) to (1.2,2.6);
\node at  (5, 0.3) { \color{blue} \tiny$q$};
\node at  (0.3, 3.75) { \color{red} \tiny $w^1$};
\node at  (1.5, 2.5) { \color{red}  \tiny $ w^2$};
\draw [line width = 0.8, red, dashed](0.6,0.8)  to  [out=40, in=250] (0.9, 1.9)   ; 
\draw [line width = 0.8, red, dashed](0.9,2.2)  to  [out=87, in=270] (0.95, 3.8)   ; 
\draw[fill, red]  (0.9, 1.9) circle [radius = 0.04];
\draw[fill, red]  (0.9, 2.2) circle [radius = 0.04];
\draw [line width = 0.7, blue](0.6,0.8)  to  [out=-5, in=230] (2.2, 1.9)   ; 
\draw [line width = 0.7, blue](2.22,2.2)  to  [out=60, in=240] (4, 3.8)   ; 
\draw[fill, blue]  (2.2, 1.9) circle [radius = 0.04];
\draw[fill, blue]  (2.22, 2.2) circle [radius = 0.04];
 \begin{scope}[white]
        \draw[fill=red!30, semitransparent] (7.5, 0) -- (7.5,2.5) -- (8,4) --  (8,1.5) --cycle;
        \draw[fill=red!30, semitransparent] (9, 0) -- (9,2.5) -- (9.5,4) --  (9.5,1.5) --cycle;
        \draw[fill=red!30, semitransparent] (10.5, 0) -- (10.5,2.5) -- (11,4) --  (11,1.5) --cycle;
\end{scope}
\draw[->, line width = 0.6, blue]  (7,0.5) to (12,0.5);
\draw[->, line width = 0.6, red]  (7.5,0) to (7.5,4);
\draw[->, line width = 0.2, red]  (7.4,0.2) to (8.2,2.6);
\node at  (12, 0.3) { \color{blue} \tiny$q$};
\node at  (7.3, 3.75) { \color{red} \tiny $w^1$};
\node at  (8.5, 2.5) { \color{red}  \tiny $ w^2$};
\draw [line width = 0.8, red, dashed](7.6,0.8)  to  [out=40, in=250] (7.9, 1.9)   ; 
\draw [line width = 0.8, red, dashed](7.9,2.2)  to  [out=87, in=270] (7.95, 3.8)   ; 
\draw[fill, red]  (7.9, 1.9) circle [radius = 0.04];
\draw[fill, red]  (7.9, 2.2) circle [radius = 0.04];
\draw [line width = 0.7, blue](7.6,0.8)  to  [out=-5, in=230] (9.2, 1.9)   ; 
\draw [line width = 0.7, blue](9.22,2.2)  to  [out=60, in=240] (11, 3.8)   ; 
\draw[fill, blue]  (9.2, 1.9) circle [radius = 0.04];
\draw[fill, blue]  (9.22, 2.2) circle [radius = 0.04];
\draw [line width = 1, red, densely dotted](7.9,1.9)  to   [out=80, in=270] (7.9, 2.2)   ; 
\draw [line width = 1, blue, densely dotted](9.2,1.9)  to [out=80, in=270] (9.22, 2.2)   ; 
 \end{tikzpicture}
 }
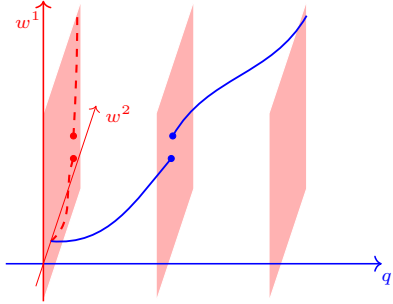
 \centerline{\tiny \bf \hskip - 0.3 cm Fig.1  $\cQ \times \cK$-projection of the graph of a solution \hskip 0.5cm  Fig.2 $\cQ \times \cK$-projection of a completed  graph }
\medskip
\subsection{A  characterisation of  the completed graphs of the solutions} \label{Sect5.3}
The purpose of this section is to  give  a purely differential-geometric characterisation of  the   graph completions  of   the solutions of a control system as   \eqref{controlsystem}. We start with the following \par
\begin{definition} \label{rigged-distribution} A {\it  rigged distribution  of  rank $m+1$} on  a  Riemannian manifold $(\cM, g^\cM)$ is a triple $\left(\cD,   \cD^I,    \cT \= \bT{\!\!\mod \!\cD^I}\right)$ given by:
\begin{itemize}
\item[(a)] a regular distribution $\cD$  of  rank $m+1$; 
\item[(b)] an involutive regular sub-distribution $\cD^I \subset \cD$ of  rank  $m$; 
\item[(c)] a  nowhere vanishing smooth section $\cT$  of  the quotient bundle 
$$\pi: \cD/\cD^I \longrightarrow \cM\ .$$
  \end{itemize}
 \end{definition}
  With the phrase   ``nowhere vanishing smooth section of $\cD/\cD^I$''  we    mean  that $\cT$ is a map 
 $x \mapsto \cT_x$, $x \in \cM$,  from $\cM$  into the union  $\bigcup_{y \in \cM}\cD_y/\cD^I_y$,   taking value in the quotient  $\cD_x/\cD^I_x$ for each $x \in \cM$ and   which  is  locally  of  the form
 \beq \label{theV}  \cT_x = \bT_x{\!\!\!\!\mod\!\cD^I}\eeq
 for some smooth  local vector field $\bT$ with values in  $\cD \setminus \cD^I$ (so that $\cT_x$ is not the trivial class of $\cD_x/\cD_x^I$ for any $x$).  Since $\bT$ is  determined by  the map  $x \mapsto \cT_x$ up to the addition of a  local vector field in $\cD^I$, we may also   denote  the section $\cT$ as   $\cT = \bT{\!\!\mod \!\cD^I}$. \par
 \smallskip
 The system  \eqref{controlsystem} is  naturally associated with the rigged distribution on the Riemannian manifold $\left(\cM = \bR \times \cQ \times \cK, g^E\right)$, with the standard Euclidean metric $g^E$ of   $\bR^{1 + n + m} (\supset \cM)$,   which is defined as follows.  Let $\bT$ and $W_1, \ldots, W_m$  be the vector fields of  $\bR \times \cQ \times \cK$,  defined   at each  $x = (t, q, w)$  by 
\beq \label{the-1} \bT|_x  \= \frac{\p}{\p t} \bigg|_x + f^i(x)\frac{\p}{\p q^i}\bigg|_x \ ,\qquad W_\a|_x \=  \frac{\p}{\p w^\a}\bigg|_{x} \  , \ \ 1 \leq \a \leq m\ ,\eeq 
and  denote by   $\cD^I$ and $\cD$ the regular distributions  defined by 
\beq\label{theriggedexample} \cD^I|_x =  \left\langle W_1|_x , \ldots, W_m|_x \right \rangle  \ ,\quad \cD|_x =  \left\langle  \bT|_x, W_1|_x , \ldots, W_m|_x \right \rangle\ .\eeq
 The triple $(\cD, \cD^I ,\bT{\!\!\mod \!\cD^I})$ is the rigged distribution {\it  canonically associated with  \eqref{controlsystem}}.
 \par
 \begin{rem}  \label{theremark} Notice that    {\it the  section $\cT = \bT {\!\!\mod \!\cD^I} $  of a manifold with a rigged distribution is invariant under  pushing-forward by   local diffeomorphisms that  preserve the integral leaves of  $\cD^I$}.  This can be checked as follows.  Being $\cD^I$  an involutive distribution,   around each point  $x_o$,  there is a  neighbourhood $\cU \subset \cM$ that  can  be  identified with a cartesian product $\cU \simeq \cU_1 \times \cU_2$ for some open sets $\cU_1 \subset \bR^{\dim \cM - m}$, $\cU_2 \subset \bR^{m}$, in which the fibers of the standard projection $\cU \simeq \cU_1 \times \cU_2 \to \cU_1$  are the   integral leaves  of $\cD^I|_{\cU}$. A local diffeomorphism  $\f:  \cU  \to   \f(\cU)$  mapping each integral leaf of   $\cD^I$  into itself, i.e. of the form $ \f(y, w) = (y, \wt \f(y,w))$, $y \in \cU_1$, $w \in \cU_2$, transforms the vector  field $\bT = \bT^i(y,w)\frac{\p}{\p y^i} + \bT^\a(y, w) \frac{\p}{\p w^\a}$ in $\cD|_{\cU}$ into 
 $\f_*(\bT) =  \bT^i(y, w)\frac{\p}{\p y^i} + \big(\bT(\wt \f)|_{(y, w)}\big)^\a \frac{\p}{\p w^\a}$. From this it follows   that the two sections 
 $$\cT|_{\cU}\= \bT|_{\cU} {\!\!\mod \!\cD^I} \ ,\qquad \f_*(\cT|_{\cU})   \=  \f_*(\bT|_{\cU}) {\!\!\mod \!\cD^I} $$
  are equal. \end{rem}
 \begin{lem} \label{lemma23} Consider  the system \eqref{controlsystem} with  associated rigged distribution $(\cD, \cD^I,$ $\bT{\!\!\mod \!\cD^I})$ on $\cM = \bR \times \cQ \times \cK$. A  piecewise regular {\rm oriented} curve  $\h = \h_1 \ast \ldots \ast \h_r$ of  $\cM$ is a completed graph  of  a solution  if and only if  it satisfies the following three conditions: 
 \begin{itemize}[leftmargin = 18 pt]
 \item[(1)]  it is a  $\cD$-path; 
 \item[(2)] for   $ 2 \leq i \leq r-1$,   the velocities  $\dot \h_i(t)$ of one (hence of any) of the consonant parameterisations of $\h_i$
 are such that  either
  \begin{itemize}  
   \item[(a)]  $\dot \h_i(t) \in \cD^I_{\h_i(t)}$ for all $t$ or
 \item[(b)]    there is an  everywhere positive  real  function $\l(t) > 0$ such that  
  \beq \label{positivity} \dot \h_i(t) =  \l(t) \bT_{\h_i(t)}{\!\!\!\!\mod \!\cD^I}\ ,\qquad \text{for any}\ t \ ;\eeq
 \end{itemize}
 \item[(3)]  for each of the two  curves $\h_i$ with  $i = 1$ or $i = r$,    there  exists an everywhere positive  $\l(t) > 0$  such that  \eqref{positivity} holds. 
 \end{itemize}
 \end{lem} 
 \begin{pf} If $\h_1 \ast \ldots \ast \h_r$ is a completed graph, for each   regular arc $\h_i$,    either  $\h_i$ is a graph of a smooth solution or the velocities $\dot \h_i(t)$ of  one of its parameterisations have zero  components along the coordinate vector  fields $\frac{\p}{\p t}$,  $\frac{\p}{\p q^j}$ and are therefore  linear combinations only of the vectors $W_\a|_{\h_i(t)}$. Hence each arc $\h_i$ is a $\cD$-path (thus (1) holds) and  the velocities  $\dot \h_i(t)$   satisfy  either  (2.a) or   (2.b) with $\l(t) \equiv 1$ (at least when $\h_i(t)$ is precisely the parameterisation $t \to (t, q(t), w(t))$ given by a solution $(q(t), w(t))$ of the system).  Since the first and the last arcs are both required to be graphs of solutions, it follows that  (2) and (3) hold. Conversely, assume that  $\h$ satisfies (1) -- (3). Then,  for any regular arc $\h_i$ that satisfies \eqref{positivity} for some parameterisation $\h_i(t)$, it is   possible to replace the parameterisation by  a new  consonant one, which satisfies \eqref{positivity} with the function $\l(t) \equiv 1$.  This new  parameterisation gives a map $t \to (q(t), w(t))$ which is a smooth solution of \eqref{controlsystem}. On the other hand,  any arc $\h_i$ that  satisfies (a) has  projections on $\bR$ and on $\bR \times \cQ$ which  are both singletons.  This means that $\h$ is a composition  of regular arcs  satisfying the  conditions of Definition \ref{completed} and is therefore a completed graph.
 \end{pf}
 \par
 \smallskip
 By Lemma \ref{lemma23},  any   graph completion of a  solution is a $\cD$-path of  the canonically associated  rigged distribution, equipped with the  orientation given  by the  parameterisations of the regular arcs $t \longmapsto (t, q(t), w(t))$ determined by their time components.  These   orientations  are also those that  establish the existence of  strictly positive functions $\l(t) > 0$  for which  \eqref{positivity} holds. This property  motivates  the  following definition, which can be used to discuss  oriented $\cD$-path on arbitrary rigged distributions, not necessarily associated with control systems.
 \par 
 \medskip
 \begin{definition} \label{positiveD-paths} Let $(\cD, \cD^I,   \bT{\!\!\mod \!\cD^I})$  be a rigged distribution on a manifold $\cM$ and, for any  $x \in \cM$ and   $ v \in \cD_x $,   let us denote by $\l^{(x,v)}$   the unique real number such that 
 $$v = \l^{(x,v)} \bT_x{\!\!\mod \!\cD^I_x}\ .$$
 The vectors $v \in \cD_x$ for which   $\l^{(x, v)} > 0$ (resp. $\l^{(x,v)} = 0$) are  called  {\it positive} (resp. {\it null)}.   An oriented $\cD$-path  $\h_1 \ast \ldots \ast \h_r$  is called {\it nonnegative} if  each of its regular arcs $\h_i$  satisfy the following condition for one (hence for all) of its consonant regular parameterisations $\h_i(t)$: {\it either all   velocities $\dot \h_i(t)$ are   positive or all of them are null and for the first and the last arcs $\h_1$,  $\h_r$  only   the first possibility  is allowed}.
   \end{definition}
By Lemma \ref{lemma23},  {\it the completed graphs of the piecewise smooth solutions of   \eqref{controlsystem}, equipped with their standard orientations, coincide with  the nonnegative $\cD$-paths of the canonically associated rigged distribution}. This is a useful  fact, because it  allows to rephrase any    question on  the  final configurations of  the solutions of a first order control system   as a problem  on  non-negative $\cD$-paths on  manifolds with  rigged distributions. 
\par
\medskip
\begin{rem} \label{remark36}  Given  a  regular $\cD$-path $\h$ in $\cU$ with positive velocities at all points, by standard arguments and possibly restricting $\cU$,   one can  see that there is at least one  local vector field  $\bT \in \cD|_\cU$ such that: (1)  $\cT = \bT  {\!\!\mod \!\cD^I} $ and (2)  one has  $\dot \h(t) = \bT_{\h(t)}   {\!\!\mod \!\cD^I} $ for any $t$.  Actually, modifying   $\bT$  by adding a vector field in $\cD^{I}$,  the vector field  $\bT$ can be chosen so that  it holds  (2')   $\dot \h(t) = \bT_{\h(t)}$.  Consider now a different  curve $\h'$ still with positive velocities at all points.  We claim that {\it  if $\cU$ is sufficiently small, for   an appropriate parametrisation of $\h'(t)$ and any sufficiently  small interval  $I$ of the parameter $t$,  there is a  local diffeomorphism $\f: \cV \subset \cU \to \cU$  such that  (a)   $\dot \h'(t) = (\f_*\bT) _{\h(t)}  $ for any $t$ in $I$  and (b) $\f$  maps  each integral leaf of $\cD^I$ into  itself and thus  is also such that $ \f_*\bT = \bT  {\!\!\mod \!\cD^I} $ }(see  Remark \ref{theremark}).  The proof is the following.  Assume that $\cU$  is small enough to be identifiable  with a  cartesian product $\cU \simeq \cU_1 \times \cU_2$ as in Remark \ref{theremark},  i.e. with  the subsets $\{y\} \times \cU_2$ equal to the  integral leaves   of  $\cD^I$ in $\cU$.   Under this identification,  the parameterised curves  $\h(t)$, $\h'(t)$ have the forms  $\h(t) = (y(t), w(t))$,  $\h'(t)  = (y'(t), w'(t))$, their  tangent vectors have the forms $\dot \h(t) = \dot y^i(t) \frac{\p}{\p y^i} + \dot w^\a(t) \frac{\p}{\p w^\a}$,    $\dot \h'(t) = \dot y'^i(t) \frac{\p}{\p y^i} + \dot w'^\a(t) \frac{\p}{\p w^\a}$,  and  the  vector field $\bT$  has the form $\bT = A^i(y, w) \frac{\p}{\p y^i} +  B^\a(y, w) \frac{\p}{\p w^\a} $ for some smooth functions $A^i$, $B^\a$ such that 
$$\dot y^i(t) =  A^i(y(t), w(t)) \ ,\qquad  \dot w^a(t) =  B^a(y(t), w(t))\ .$$ 
From the assumption that  $\h'(t) = (y'(t), w'(t))$ has positive velocities,  there is re-parameterisation of  $\h'(t)$ such that also the following equality holds
$$\dot y'^i(t) =  A^i(y'(t), w(t))\ .$$ 
Since the $\cU_1$-parts $\dot y^i(t) \frac{\p}{\p y^i} $,  $\dot y'^i(t) \frac{\p}{\p y^i} $ of the  vectors $\dot \h(t)$, $\dot \h'(t)$  are nowhere vanishing, there exist  two indices  $  i_1,  i_2 $  such that the real functions $t \to y^{i_1}(t)$,  $t \to y^{i_2}(t)$ have nowhere vanishing derivatives  and are therefore invertible.  This means  that the coordinate maps $y^{i_1}, y^{i_2}: \cU \to \bR$, when restricted to the curves $\h$, $\h'$,   determines real functions $y^{i_1}(t)$,   $y^{i_2}(t)$  that   give new parameters  for $\h$,  $\h'$. 
There are  therefore    two smooth maps $\psi, \psi': \wt I \subset \bR \to \cU_2$  from a suitably small   interval of $\bR$ into $\cU_2$ such that   
$\h(t) = \big(y(t),  \psi(y^{i_1}(t))\big)$,  $\h'(t) =  \big(y'(t),   \psi'(y^{i'_1}(t))\big)$.
If we now consider the local map $\f: \cU_1 \times \cU_2 \to \cM$ defined by  (here, $y^{i_1}$, $y^{i_2}$ are considered as    maps of the form $y^{i_1}, y^{i_2}: \cU  = \cU_1 \times \cU_2 \to \bR$)
$$\f(y, w) = \big(y, w - \psi\circ y^{i_1} + \psi' \circ y^{i'_1}\big)$$
we see that $\f$ satisfies the property (b) and is such that  the $\cU_2$-components of the  curves $\h'(t)$ and $\f \circ \h(t)$  are equal for all $t$ in an appropriate small interval $I \subset \bR$. In particular, we have that the $\cU_2$-components of the vectors $\dot \h'(t)$ and $\frac{d }{dt} (\f \circ \h)\big|_t  = \f_*(\dot \h(t)) = (\f_*\bT)_{\h'(t)}$ coincide for all $t\in I$. Thus, since we know that  the $\cU_1$-components of the three vectors $\dot \h'(t)$, $\bT_{\h'(t)}$ and of  $ (\f^*\bT)_{\h'(t)}$ are all equal, also  (a) holds. 
\end{rem} 
\par
\medskip
\section{ $\cM$-attainable sets  and reachable sets}
\label{section4}
\subsection{Reachable sets, accessibility and small time  local  controllability}
Let us now introduce  the notions of  reachable set, accessibility and small-time controllability  for  a  control system as in \eqref{controlsystem}.  The following definitions  are essentially equivalent to the most commonly adopted,  except for   additional restrictions on the regularity   of the considered solutions (see e.g.  \cite{Bl, AS, BP, MS}).
\par
\begin{definition}  Given a configuration $q_o \in \cQ$ and $T \in (0, + \infty)$
 the {\it reachable set   in time  exactly equal to $T$ and by means of  piecewise $\cC^k$  solutions}  (where the index $k$ is possibly equal  to $\infty$ or $\o$) is the subset of $\cQ$ defined by 
 \begin{multline} \Reach^{\cC^k}_{ T}(q_o) {:=}  \bigg\{ q \in \cQ\,:\, q\ \text{is the final point of a  piecewise $\cC^k$ curve $q(t)$}\\
 \text{ that starts at $q_o$ and is the projection on} \ \cQ\\
  \text{of a piecewise $\cC^k$ solution $(q(t), w(t))$ to \eqref{controlsystem} with}\  t \in [0, T]\   \bigg\}\ .
 \end{multline}
 The {\it reachable set in time  $T$} (resp. {\it in any time}) {\it through piecewise $\cC^k$  solutions} is 
\beq \Reach^{\cC^k}_{\leq T}(q_o) {=} \{q_o\} \cup  \bigcup_{\wt T \in (0,  T]} \Reach^{\cC^k}_{ \wt T}(q_o)\ \ \left( \ \text{resp.}\ \  \Reach^{\cC^k}_{< \infty}(q_o) {=} \bigcup_{ T \in (0, \infty)} \Reach^{\cC^k}_{ \leq T}(q_o)\ \right)\ .\eeq
 The system  \eqref{controlsystem}  is said to have  the   {\it hyper-accessibility} (resp.  {\it accessibility})  {\it property in the $\cC^k$ sense} if  for any $q_o \in \cQ$ and $T> 0$, the set  $\Reach^{\cC^k}_{ T}(q_o)$ (resp.  $\Reach^{\cC^k}_{\leq T}(q_o)$ ) is open (resp. has non-empty interior).  It is said to have the  {\it  small-time local controllability property  at $q_o \in \cQ$ in the $\cC^k$ sense} if  there is a  $T> 0$ such that $\Reach^{\cC^k}_{\leq T}(q_o)$  contains a neighbourhood of $q_o$. 
 \end{definition}
 Note that, according to the above definitions,  a  system having  the hyper-accessibility property in $\cC^k$-sense, is also accessible in $\cC^k$-sense and,  for any $q_o$, for which there is   a $\cC^k$-solution $(q(t), w(t))$,   $t \in [0, T]$, $T> 0$, satisfying $q(0) = q_o = q(T)$,  there is also the small-time local controllability at $q_o$.
   \par
 \smallskip
 In what follows, a point $q_o$ with the above  property will be briefly indicated as     {\it  a point with   the homing property}. Among them, there are  the  {\it stable points (in the $\cC^k$-sense)}, namely the points $q_o$ for which there exists  a $\cC^k$ solution $(q(t), w(t))$ with $q(t) \equiv q_o$. 
\par
\subsection{$\cM$-attainable sets  in $\cM$ and  their relation with the  reachable sets} 
\label{Sect62}
Let $(\cM, g^\cM)$ be a Riemannian manifold equipped  with the rigged distribution $ (\cD, \cD^I,  \cT \= \bT{\!\!\mod \!\cD^I})$.  For any  $x_o \in \cM$,  the  {\it $\cM$-attainable set of  $x_o$}  is the set 
 \begin{multline} \Attain^{\cC^k}_{x_o} {=} \bigg\{ y \in \cM\,:\, y\ \text{is the final endpoint of a  non-negative  $\cD$-path} \\
 \text{ that starts at}\ x_o \ \text{and  with regular arcs  of class}\ \cC^k\ \bigg\}\ .
 \end{multline}
In case $\cM = \bR \times \cQ \times \cK$ and  $ (\cD, \cD^I,  \cT \= \bT{\!\!\mod \!\cD^I})$ is the canonical   rigged distribution associated with  \eqref{controlsystem}, for any   $x_o \in \cM$ of the form $x_o = (t_o = 0, q_o, w_o)$, the corresponding $\cM$-attainable set  $\Attain^{\cC^k}_{x_o} $ is the set given by the  final endpoints $(t_{\text{fin}}, q(t_{\text{fin}}), w(t_{\text{fin}}))$ of  the graph completions  of the piecewise $\cC^k$ solutions $(q(t), w(t))$, $t \in [0, t_{\text{fin}}]$,  of \eqref{controlsystem}  with   initial conditions $q(0) = q_o$ and $w(0) = w_o$. 
Hence,  if we denote by
 $$\pi^\bR: \cM \to \bR\ ,\qquad \pi^\cQ: \cM \to \cQ$$
    the standard projections of $\cM = \bR \times \cQ \times \cK$ onto  $\bR$ and $\cQ$, respectively,    each reachable set $\Reach^{\cC^k}_T(q_o)$  is equal to 
\beq
\Reach^{\cC^k}_T(q_o) =    \pi^\cQ \left(  \bigg(\bigcup_{\smallmatrix  w_o  \in \cK\\
 \endsmallmatrix}  \Attain^{\cC^k}_{(0, q_o, w_o)} \bigg) \cap (\pi^{\bR})^{-1}( T) \right) \ .
 \eeq
  This  and the above discussion immediately implies  the following 
 \begin{prop} If for any $q_o \in \cQ$ and $T > 0$, there exists a $w_o \in \cK$ such that  $\pi^\cQ\big( \Attain^{\cC^k}_{(0, q_o, w_o)} \cap (\pi^\bR)^{-1}(T)\big) $  is open, then the system   \eqref{controlsystem} is   hyper-accessible in the $\cC^k$ sense and it is small-time locally controllable at the points with the homing property.
 \end{prop}
\vskip 25pt
\addcontentsline{toc}{section}{Part II}
\centerline{\large \it PART II}
\ \\[-45pt]
\section{Secondary  distributions  and  $\bT$-surrogate fields} \label{thefifthsection}
\subsection{The secondary  distribution of  a  rigged distribution} \label{sect63}
Let  us now focus on the case of a  {\it real analytic} Riemannian $N$-dimensional manifold  $(\cM, g^\cM)$    equipped   with  a {\it real analytic} rigged distribution  $(\cD, \cD^I,  \cT = \bT{\!\!\mod \!\cD^I})$ of rank $m+1$, that is consisting of two  real analytic distributions $\cD$, $\cD^I$ and a real analytic section $\cT$ of $\cD/\cD^I$.  
Note that the  rigged distribution on $\cM = \bR \times \cQ \times \cK$, which is canonically associated with  a real analytic control system  \eqref{controlsystem} automatically satisfies these conditions.  \par
\smallskip
 We  denote by    $V^{II} = V^{II}_\cD$  the  smallest family of  real analytic local vector fields on $\cM$, which is closed under sums and multiplications by  real analytic functions  and  contains
\begin{itemize}[leftmargin = 15pt]
\item   all real analytic  local vector fields $Z \in \cD^{I}$; 
\item   all   real analytic  fields  of the form $Z' = [Y_1, [Y_2, \ldots, [ Y_r, Z]\ldots]]$, $r \geq 1$, for some local  real analytic  $Z  \in \cD^I$ and   real analytic  local vector fields $Y_i$ having  the form 
\beq\label{condsecond}  Y_i = \bT  + W_i\ ,\qquad W_i \in \cD^I\ .\eeq
\end{itemize}
We may alternatively define  $V^{II}$   as follows. For any  $y \in \cM$,  let $\cU$ be  a neighbourhood  on which there is  a set of real analytic  generators $\{X_1, \ldots, X_m\}$  for  $\cD^I|_{\cU}$. 
By possibly restricting $\cU$, we may   assume that there is also a local real analytic  vector field $\bT$ on $\cU$,   such that  the section  $\cT|_{\cU}: \cU \to \cD/\cD^I $  coincides with the  family of  equivalence classes  $\cT_x  = \bT_x{\!\!\!\!\mod\!\cD^I}$.  On such $\cU$,  any local vector  fields   $Z \in \cD^I|_{\cU}$ and $Y_i = \bT +W_i$, $W_i \in \cD^I|_{\cU}$, have the form 
$$ Z =  \mu^\b X_\b \ ,\qquad Y_i = \bT +\l^\g_i X_\g .$$
Therefore the  Lie bracket $[Y_i, Z]$ has  the  form  
\beq [Y_i, Z] = [\bT +   \l^\g_i X_\g,  \mu^\b X_\b] =   \mu^\b [\bT, X_\b]  + \text{ vector fields in} \ \cD^I \ ,\eeq
meaning that   it  is a pointwise linear combination of  the  $X_\a$ and of  the Lie brackets  $[\bT, X_\b]$. 
This  and similar computations concerning  the iterated Lie brackets $[Y_1, [Y_2, \ldots, [ Y_r, Z]\ldots]]$ show that we  may consider the following  equivalent  definition for  $V^{II}_\cD$: {\it it   is the set of all local real analytic vector fields $Y$ having the form 
\beq \begin{split}
&Y =  \sum_{\ell = 0}^\nu f^{\b}_{(\ell)} X_\b^{(\ell)}\ ,\qquad f^{\b}_{(\ell)}: \cU \to \bR\ \text{real analytic}\ , \ \nu \in \bN\\
&\qquad  \text{where}\ \label{genfor} X_\b^{(\ell)} \= \underset{\text{$\ell$-times}}{\underbrace{[\bT, [\bT, \ldots, [ \bT, }}X_\b]\ldots]]\ ,
\end{split}
\eeq
for   some  set of real analytic generators $\{X^{(0)}_\a = X_\a\}$ of $\cD^I$ and a real analytic vector field $\bT$  for which   $\bT_x{\!\!\!\!\mod\!\cD^I}  = \cT_x $  at all  points $x$ where  the field $\bT$ is defined}. \par
\smallskip 
We now  recall that,   by the Noetherian property of the rings of real analytic functions (\cite[Thm. 3.8]{Ma}), for any $y \in \cM$, there  is a neighbourhood $\cU$ on which there are real analytic  generators $X_\a$, $1 \leq \a \leq m$,   for $\cD^I|_\cU$,  a vector field $\bT$ such that $\cT|_{\cU}  = \bT{\!\!\mod \!\cD^I}$ and a  finite set of vector fields of the form \eqref{genfor},  say $X^{(\ell_1)}_{\b_1}\ ,\ \ldots, \  X^{(\ell_{q})}_{\b_{q}} $
($q$  might  depend   on $\cU$),  such that any other vector field  $X^{(\ell)}_\b$  is  linear combination  of them  with real analytic components.
%
 From   this,   it follows  that  
 $V^{II}$ is a quasi-regular set of real analytic vector fields of rank $p$, with $p$ equal to the maximal dimension of the   spaces $\cD^{V^{II}}_x  = \{X_x\ ,\ X \in V^{II} \}$, $x \in \cM$. Hence    $(V^{II}, \cD^{II} \= \cD^{V^{II}})$ is a  generalised distribution in the sense of Definition \ref{quasi-regular-def}.
 We are now ready to introduce the following crucial notion, whose  relevance  will be   explained in Remark \ref{remark56}. 
\begin{definition}
The generalised 
 distribution $(V^{II}, \cD^{II} \= \cD^{V^{II}})$ defined above  is   the  {\it secondary  distribution}  associated with  the real analytic rigged distribution   $(\cD, \cD^I,   \cT \= \bT{\!\!\mod \!\cD^I})$.
  The rank of $(V^{II}, \cD^{II} )$ is called the    {\it  secondary rank} of
 the rigged distribution. Given a set of local generators $\{Y_i\}_{i = 1, \ldots, p_\cU}$  on an open set $\cU$ for $V^{II}|_{\cU}$, each of them having  the form $Y_i = X_{\a_i}^{(\ell_i)}$  for a vector field $\bT$  with  $  \cT \= \bT{\!\!\mod \!\cD^I}$ and a set of local generators $\{X_\a\}$ for $\cD^{I}|_{\cU}$, 
 the  integer $\nu \= \max\{\ell_1, \ldots, \ell_{p_\cU}\}$  is called  {\it  $\bT$-height of the set of generators}. 
The smallest value for the $\bT$-heights of sets of generators on  $\cU$  is called  {\it secondary height of the rigged distribution on $\cU$}.
\end{definition}
\begin{example} \label{example72}
Some of the  most important  examples of secondary distributions are provided by the {\it  linear control systems}, i.e. by the systems of the form  
\beq \dot q^i = A^i_j q^j + B^i_\a w^\a\ ,\qquad A = (A^i_j) \in \bR_{n \times n}\ ,\ \ B = (B^i_	\a) \in \bR_{n \times m}\ . \eeq
For a system of this kind,  the corresponding rigged distribution is determined as in \eqref{theriggedexample} by the {\it globally defined} vector fields 
\beq \label{KAL1} \bT  \= \frac{\p}{\p t}  + (A^i_j q^j + B^i_\a w^\a)\frac{\p}{\p q^i}\ ,\qquad W_\a \=  \frac{\p}{\p w^\a} \in \cD^I\ .\eeq
In this case, the set of vector fields $V^{II}$ is given by all real analytic  local vector fields which are 
finite linear combinations -- with coefficients given by real analytic functions -- of the vector fields 
 \begin{align*}
 &   W^{(0)}_1\=  W_1 \ ,\ \ldots \ ,&& \ldots\ , && W^{(0)}_m \=  W_m \ ,\\
 & W^{(1)}_1 \=  [\bT, W_1] =-  B^i_1\frac{\p}{\p q^i}\ ,\ \ldots\ ,&& \ldots\ , &&W^{(1)}_m \=  [\bT, W_m] = - B^i_m\frac{\p}{\p q^i} \ ,\\
 & W^{(2)}_1 \= [\bT, [\bT, W_1]] =   A^i_j B^j_1\frac{\p}{\p q^i}\ ,\ \ldots\ , &&  \ldots\ , &&   W^{(2)}_m \=  [\bT, [\bT, W_m] ]=   A^i_j B^j_m\frac{\p}{\p q^i}\ ,\\
 & \text{etc.} && && 
 \end{align*}
These  vector fields generate  $V^{II} $, have the form \eqref{genfor} and  have constant coefficients. Hence   {\it  the generalised distribution $(V^{II}, \cD^{II})$  is  regular}. Its rank   is equal to  the first integer $n_{\ell_o}$ after which the monotone sequence of dimensions $n_0 \leq n_1\leq n_2 \leq \ldots$   defined by   
\begin{multline*} n_\ell = m +  \dim \left\langle B_1, \ldots, B_m, \ \ A{\cdot}B_1, \ldots, A{\cdot}B_m, \ \ A {\cdot} A{\cdot}B_1, \ldots, A {\cdot} A{\cdot}B_m,\right. \\
 A {\cdot} A {\cdot} A{\cdot}B_1, \ldots, A {\cdot} A {\cdot} A{\cdot}B_m\ ,
\left.\ldots,  \underset{\text{$(\ell-1)$-times}}{\underbrace{A{\cdot} \ldots {\cdot}A}}{\cdot}B_1,\ \  \ldots, \ \ \underset{\text{$(\ell-1)$-times}}{\underbrace{A{\cdot} \ldots {\cdot}A}}{\cdot}B_m\right \rangle\ ,\end{multline*}
 stabilises.   
\end{example}
\begin{example} \label{example72*} An elementary  example  in which the secondary distribution  is {\it singular}  is given by the control equation on  curves $q(t)$  in $\cQ = \bR$ with  controls $w(t)$ in $\cK = \bR$ 
\beq \dot q = w^2\ .\eeq
In this situation, $\cM = \bR^3$ and the canonically associated rigged distribution is determined as in \eqref{theriggedexample} by the {\it globally defined} vector fields 
\beq \label{KAL1*} \bT  \= \frac{\p}{\p t}  +  w^2\frac{\p}{\p q}\ ,\qquad W \=  \frac{\p}{\p w} \in \cD^I\ .\eeq
It is quite immediate to verify that  
the  vector fields in  $V^{II}$  are  the real analytic  local vector fields that are 
pointwise  linear combinations of the  vector fields
\beq W^{(0)} = \frac{\p}{\p w}\ ,\qquad   W^{(1)} = \big[\bT, \frac{\p}{\p w}\big] = [\bT, W^{(0)}] = - 2  w\frac{\p}{\p q}\ .\eeq
Note that in this case the secondary  distribution $(V^{II}, \cD^{II})$ is singular: Indeed the dimension of  $ \cD^{II}_{x} $ is $2$ at  all points $x = (t, q,w) $ with $w \neq 0$ and is $1$ otherwise.
\end{example}
\par
The next lemma  shows that around any point  of $\cM$ there exists a set  of generators for the secondary distribution with certain useful properties and  admitting a convenient  order. It is a merely technical result, but  it  implies   a significant   simplification   for the arguments of our first main result.  The  proof  of this lemma   is postponed to the Appendix \ref{appendixb}.\par
 \begin{lem}  \label{lemmone} Let  $(\cD, \cD^I, \cT \= \bT{\!\!\mod \!\cD^I})$  be a  real analytic rigged distribution. Let also $\bT$ be a local vector field such that  $\cT_x =  \bT_x{\!\!\mod \!\cD^I_x}$ for all $x$ in an open set $\cV \subset \cM$.  For any sufficiently smaller open set $\cU \subset \cV$, there are    integers $R_a$,   $0 \leq a \leq \nu$,  $\nu$   secondary height  on $\cU$, and   a set of real analytic generators for  $(V^{II}|_{\cU}, \cD^{II}|_{\cU})$, indexed in the following way \\
\resizebox{0.98\hsize}{!}{\vbox{
 \begin{align*} & W_{0(0)1}, \ldots, W_{0(0)R_{0}},  &&W_{0(1)1}, \ldots, W_{0(1)R_{1}},  &&   \ldots, && \ W_{0(\nu -1)1}, \ldots, W_{0(\nu-1)R_{\nu-1}}, &&W_{0(\nu)1}, \ldots, W_{0(\nu)R_{\nu}}, \\[10pt]
&  & &W_{1(1)1}, \ldots, W_{1(1)R_{1}},  && \ldots && W_{1(\nu-1)1}, \ldots, W_{1(\nu-1)R_{\nu -1}},&&W_{1(\nu)1}, \ldots, W_{1(\nu)R_{\nu}}, \\[10pt]
 & && &&    && \ddots && \vdots && \\[10pt]
   & && && && W_{\nu-1(\nu -1)1}, \ldots, W_{\nu-1(\nu-1)R_{\nu-1}}, && W_{\nu-1(\nu)1}, \ldots, W_{\nu-1(\nu)R_{\nu}}, \\[10pt]
   & && && && &&  W_{\nu(\nu)1}, \ldots, W_{\nu(\nu)R_{\nu}}\ , 
   \end{align*}}}
   \\
   such that  the next conditions hold:
\begin{itemize}
\item[(1)] The vector fields appearing in the first row  of the previous formula, i.e.  
$$W_{0(0)1}\ , \ \ldots\ , \ W_{0(0)R_0}\ , W_{0(1)1}\ , \ \ldots\ , \ W_{0(1)R_1}\ ,\ \ldots\ ,\ W_{0(\nu)1}\ , \ \ldots\ , \ W_{0(\nu)R_\nu}$$ 
are  generators for   $\cD^{I}|_{\cU}$;  in particular,  $R_0 + R_1 + \ldots + R_\nu =  m$; 
\item[(2)]  For  any   $1\leq a\leq \nu$, $1 \leq \ell \leq a$ and $1 \leq j \leq R_{a}$,  the vector field $W_{\ell(a) j}$  is obtained  from the  vector field $W_{0(a)j}$ by applying $\ell$-times the operator $X \mapsto [\bT, X]$, i.e. 
\begin{multline} \label{defW} [\bT, W_{0(a)j}] = W_{1(a)j}\ ,\ \ [\bT, [\bT, W_{0(a)j}]]= W_{2(a)j}\ ,\ \ \ldots \ \\[3 pt]
\ldots \ \underset{\ell\text{\rm -times}}{\underbrace{[\bT, [\bT, [\bT,  \ldots [\bT, }}W_{0(a)j}] \ldots ]]]= W_{\ell(a)j}\ , \ \ \ldots \ \  \\[-10 pt]
\ldots\ , \ \ \underset{a\text{\rm -times}}{\underbrace{[\bT, [\bT, [ \bT, \ldots [\bT, }}W_{0(a)j}] \ldots ]]]= W_{a(a)j}\ .
\end{multline}
 \item[(3)] For any   $0 \leq a\leq \nu$ and $1 \leq j \leq R_{a}$ the iterated Lie bracket
$$   \underset{(a+1)\text{\rm -times}}{\underbrace{[\bT, [\bT, [\bT,  \ldots [\bT, }}W_{0(a)j}] \ldots ]]] $$
is pointwise a linear combination of the  generators $W_{\ell'(b)j'}$ having  $b \leq a$.
\end{itemize}
 \end{lem}
 Given a  local vector field  $\bT$ and an open set $\cU \subset \cM$  satisfying the hypotheses of   Lemma \ref{lemmone}, 
a  set of (local) generators $(W_{j(a)\ell})$   as in the statement  is  called  {\it $\bT$-adapted}. For  any such a set  of generators,     we  assume that the indices $``\ell (a) j"$  are     lexicographically ordered, i.e. we assume that 
$ ``\ell (a) j" < ``\wt \ell (\wt a)\wt j"$  if  $ \ell <\wt  \ell$ or  $\ell =\wt  \ell$ and  $a< \wt a$ or $\ell =\wt  \ell$, $a = \wt a$ and    $j < \wt j$.
In this way, if   ${\bf m}$ is the secondary rank of the rigged distribution,   each  index $``\ell(a)j"$ is uniquely representable by   its  position  $1 \leq A \leq {\bf m}$   within the set   of all such triples.   By (1) of Lemma \ref{lemmone},  the first $m$ vector fields $W_1, \ldots, W_m$ (i.e.  the vector fields
  $W_{0(a) j}$, $0 \leq a \leq \nu$, $1 \leq j \leq R_a$) are   generators for   $\cD^I|_\cU$. In order to easily distinguish the subset  $(W_1, \ldots, W_m)$ from its complementary subset, we   sometimes  denote them as 
  $(W^{(I)}_\a)$ and $(W^{(II \setminus I)}_B)$,  respectively.
 \par
 \medskip
 \subsection{Stepped $\bT$-paths and   $\bT$-surrogate fields} \label{stepandsurr}
 As in the previous section,   $  ( \cD, \cD^I,  \cT\= \bT{\!\!\mod \!\cD^I})$  is a  real analytic rigged distribution on   $(\cM,g^\cM)$.   We  assume that  $x_o  \in \cM$  is a fixed point and 
 $$\h = \h_1 \ast \ldots \ast \h_{r-1} \ast \h_r$$
is  a piecewise $\cC^\o$ non-negative  $\cD$-path   originating    from $x_o$. 
We  denote by  $x_o'$ and   $y_o$   the final points of  $\h_1 \ast \ldots \ast \h_{r-1}$  and of 
$\h = \h_1 \ast \ldots \ast \h_{r-1}\ast \h_r$, respectively  (they  are  therefore  both   points of    $\Attain^{\cC^\o}_{x_o} $ and, in addition,  $y_o$ is a point of $\Attain^{\cC^\o}_{x'_o} $).   
  Since  the final arc $\g = \h_{r}$ can be    considered as   a finite composition of  arbitrarily small positive  sub-arcs,  replacing $\h_r$ by  one of such  compositions, there is   no loss of generality if   we   assume that   $x_o'$  is  contained in  a  prescribed  small neighbourhood $\cU$ of  $y_o$, on which there are   a set of $\bT$-adapted generators $(W_A) = (W^{(I)}_\a,W^{(II \setminus I)}_B)$ for $(V^{II}, \cD^{II})$,    with $W^{(I)}_\a \in \cD^I$, as defined in \S\ref{sect63}, 
 and 
 a local real analytic section   $\bT$ of $\cT$  with the property that  $ \g $  is an integral curve of $\bT$ (i.e.   of the form  $\g(t) = \Phi^\bT_t(x_o')$, $ t \in [0, T]$ -- the existence of  $\cU$ and $\bT$ with this property  is proved in Remark \ref{remark36}). In what follows, we begin our  analysis of   the  local structure of the set    $\Attain^{\cC^\o}_{x_o} $ in proximity of $y_o$. \par
\smallskip 
\subsubsection{Stepped $\bT$-paths and their $\bT$-surrogates} 
 \label{sect7.2.2}
 Our goal is to show that in a sufficiently small neighbourhood of $y_o$, the set  $\Attain^{\cC^\o}_{x_o} $ contains a special subset of 
points, all of them joinable to $y_o$ through certain (non-oriented) curves tangent to the secondary distribution. For this, we first  need to recall a classical result  on composition of flows.\par
\smallskip
 Let $Y$  be a vector field on a manifold $\cN$ 
 and $f: \cV \subset \cN \to \cN$ a local diffeomorphism.  Then the   flows of $Y$  and of its pushed-forward  field $f_*(Y)$  satisfy the relation 
 \beq f\left(\Phi^Y_t(x)\right) = \Phi^{f_*(Y)}_t\left(f(x)\right) \eeq
 for any $x$,  $t$ where both sides are  defined (for checking this very simple relation, it suffices to observe    that the left and the right hand sides are  $t$-parameterised curves with same initial point and same velocity for  each $t$).  Therefore,  given  two vector fields  $X,Y$, for any sufficiently small $\ve> 0$ and for $t, s \in [0, \ve]$ we have that
 \beq \Phi^X_s \circ \Phi^Y_t (x) = \Phi^{Y^{s}}_t\circ \Phi^X_s(x)\ ,\qquad \text{where }\ \ Y^{s}\= (\Phi^X_s)_*(Y)\eeq
and  the two piecewise regular curves $\wc \s_1 \ast \wc \s_2$ and $\wt \s_2 \ast \wt \s_1$, with regular arcs  given by  
\beq \label{327}
\begin{split}
&\wc \s_1(t) \= \Phi^Y_t( x) \ ,\ t \in [0, t_o]\ ,\qquad \wc \s_2(s) = \Phi^X_{s} \left( \Phi^Y_{t_o}(x)\right) \ ,\ s \in [0, s_o]\ ,\\
&\wt \s_1(s) \= \Phi^X_s( x) \ ,\ s \in [0, s_o]\ ,\qquad \wt \s_2(t) = \Phi^{Y^{s_o}}_t \left( \Phi^X_{s_o}(x)\right)  \ ,\ t \in [0, t_o]\ ,
\end{split}
\eeq
with $t_o, s_o \in [0, \ve]$, have the same initial and final points.
This simple fact is the main ingredient of the following fundamental lemma. \par
 \begin{lem} \label{lemma4.10}  Let  $\g(t) = \Phi^{\bT}_t$, $t \in [0, T]$,  be the regular positive arc from  $x_o'$ to  $y_o$ as above and denote by  
 $$\wc \g = \wc \g_1 \ast \ldots \ast \wc \g_{2k} \ast \wc \g_{2k+1}\ ,$$
     a  non-negative  $\cD$-path that starts from $x'_o$  and   is defined as a  composition of $2k+1$  real analytic  arcs satisfying the following conditions:
     \begin{itemize}[leftmargin = 20pt]
     \item[(1)] the odd arcs  $\wc \g_{2 \ell + 1}$, $0 \leq \ell \leq k$,   are 
positive oriented  integral curves  of the vector field $\bT$ (i.e.    of the form $\wc \g_{2 \ell + 1}(t) = \Phi^\bT_t(x_{2 \ell +1})$  with $t$ running in some interval  $ [0, \s_{2 \ell +1}]$); 
\item[(2)] 
  the even  arcs $\wc \g_{2\ell}(t)$, $1 \leq \ell \leq k$,  are integral curves  of   vector fields  $\l^\a_{2\ell} W^{(I)}_\a$, $\l^\a_{2\ell} \in \bR$  (i.e.   of the form $\wc \g_{2 \ell}(t) = \Phi^{\l^\a_{2\ell} W^{(I)}_\a}_t(x_{2 \ell})$  with  $t$ varying in some interval  $ [0, \s_{2 \ell}]$); 
  \item[(3)] the widths  $\s_{2 \ell +1}$ of the ranges of  the parameters of  the odd arcs, are such that 
   \beq \label{3.28} \sum_{\ell = 0}^k \s_{2 \ell +1} = T \qquad \text{so that}\qquad  \Phi^{\bT}_{ \sum_{\ell = 0}^k \s_{2 \ell +1} = T}(x_o') = \g(T) =  y_o\ .\eeq
  \end{itemize}
 Then the  final point  of  $\wc \g$   coincides with the final point of  the  piecewise regular   curve with initial point $x'_o$
\beq \label{ass} \wt \g =   \g \ast   \wt \g_2 \ast \wt \g_4 \ast  \ldots \ast \wt \g_{2k}  \ , \eeq 
in which the sub-curve      $\wt \g_2 \ast  \wt \g_4 \ast \ldots \ast \wt \g_{2k}$ is a piecewise regular curve originating from $y_o$ and whose regular arcs   $\wt \g_{2 \ell}$    are integral curves of the vector fields 
$\l^\a_{2 \ell}W_\a^{(I)\t_{2\ell}}  \=  (\Phi^\bT_{\t_{2\ell}})_*(\l^\a_{2\ell}W^{(I)}_\a) $ with $\t_{2 \ell}$ equal to  $\t_{2 \ell} \= \sum_{j = \ell}^k \s_{2 j + 1}$ (see Fig. 4).
 \end{lem}
 \begin{rem} \label{remark56} In case of a rigged distribution associated with a real analytic control system  \eqref{controlsystem}, the curve $\g$ and any curve $\wc \g$ satisfying the conditions (1) and (2) of Lemma \ref{lemma4.10} are graphs or graph completions in $\cM = \bR \times \cQ \times \cK$ of   two piecewise smooth solutions, both with  the same initial condition and  both  defined   on an identical    time interval  $[t_o, t_o + T]$.   Thus,  Lemma \ref{lemma4.10} implies that  the  points of a neighbourhood of $y_o$ in $\cM$, which can be reached  through a graph completion   $\wc \g$,    are at the same  time  the   points which can be reached moving first along the curve   $\g$  joining  $x_o'$ to $y_o$  and then, starting from $y_o$,  following the integral curves of  certain new vector fields. These are  the $\Phi^\bT$-pushed-forward  fields  of  the  vector fields  $\lambda ^\a W_\a^{(I)}$ that  determine the ``added parts'' of the graph completions $\wc \g$.   As we will shortly see, those $\bT$-pushed-forward  fields are vector fields in the secondary generalised distribution $(V^{II}, \cD^{II})$. This means that   near  $y_o$, there is  an important class of  $\cM$-attainable points of $x_o$,   
 which  are  joined to  $y_o$ through  a well defined class of  paths that are tangent to   the secondary distribution $(V^{II}, \cD^{II})$ and with trivial projections on the time axis.   Note that  {\it there is no distinguished ``time'' orientation on such  $\cD^{II}$-paths} (along each of them the time coordinate is constant). This is a  feature  that makes them much more useful   in  the analysis of  the $\cM$-attainable sets than the  graph completions, because the latter are not just arbitrary  $\cD$-paths,  but {\it  $\cD$-paths equipped with a distinguished  time orientation}.  This is  essentially the main  reason  of   interest   for the distribution    $(V^{II}, \cD^{II})$ and for   the   $\cD^{II}$-paths.  
 \end{rem}
 \begin{pf} Since  the regular arc  $ \wc \g_{2k}$ is an integral curve of  the vector field $\l^\a_{2k} W^{(I)}_\a$ and   $\wc \g_{2k+1}$ is an integral curve of $\bT$, by   \eqref{327}  the piecewise regular curve  $ \wc \g_{2k} \ast \wc \g_{2k+1}$ has the same endpoints of an appropriate piecewise regular curve $\wt \g_{2k+1} \ast \wt \g_{2k}$, where  $\wt \g_{2k+1}$ is an integral curve of $\bT$  and the second regular arc is an integral curve of 
  $\l^\a_{2 k}W_\a^{(I)\t_{2k}}  \=  (\Phi^\bT_{\t_{2k}})_*(\l^\a_{2k}W^{(I)}_\a) $ 
with $\t_{2k} = \s_{2k +1}$.    This implies that $\wc \g = \wc \g_1 \ast \ldots \ast \wc \g_{2k -1}\ast \wc \g_{2k} \ast \wc \g_{2k+1}$ has the same endpoints of the piecewise regular curve $ \wc \g_1 \ast \ldots \ast \wc \g_{2k-1} \ast \wt \g_{2k+1} \ast \wt \g_{2k}$.  Note that the parameter $t$ of the flow which determines  $\wc \g_{2k+1}$  is the same parameter of the flow that gives  $\wt \g_{2k+1}$ and it is therefore running in $[0, \s_{2 k+1}]$. \par
In the new curve, the  adjacent  arcs $ \wc \g_{2k -1}$ and $\wt \g_{2k+1}$ are both integral curves of $\bT$ and  their composition is a single regular arc,  parameterised by  $t\mapsto \Phi^\bT_t(x_{2 k-1})$  with   $t \in [0, \s_{2 k-1} + \s_{2 k +1 } ]$.  For simplicity of notation,  let us  denote such longer arc  just by  $\wc \g_{2k-1}$.  Now, by the same  argument of before,   we see  that the piecewise regular curve $\wc \g_{2k-2} \ast \wc \g_{2k-1}$ has the same endpoints of a new piecewise regular curve of the form  $\wt \g_{2k-1} \ast \wt \g_{2k-2}$, where   
$\wt \g_{2k-1}$ is an integral curve of $\bT$ and $\wt \g_{2k-2}$  is an integral curve of  $\l^\a_{2 k-2}W_\a^{(I)\t_{2k-2}}$ with $\t_{2k-2} \= \s_{2 k -1} + \s_{2 k +1}$. We may therefore replace the  composed curve  $\wc \g_{2k-2} \ast \wc \g_{2k-1}$ with  $\wt \g_{2k-1} \ast \wt \g_{2k-2}$
and obtain a new piecewise regular curve 
$ \wc \g_1 \ast \ldots\ast (\wc \g_{2k -3} \ast \wt \g_{2k -1} )\ast \wt \g_{2k-2} \ast \wt \g_{2k}$, still with the same endpoints. Once again, the piecewise regular curve  $\wc \g_{2k -3} \ast \wt \g_{2k -1} $ is given by just one integral curve of $\bT$ and it is given by the flow of $\bT$ with parameter $t$ running in $[0, \s_{2 k -3} + \s_{2 k -1} + \s_{2k +1}]$. As before, we denote this arc simply by $\wc \g_{2 k -3}$. Iterating  $k$-times this argument, we conclude that the original piecewise regular curve  $\wc \g$ has the same endpoints of a new  piecewise regular curve of the form $\wt \g =   \wt \g_1 \ast   \wt \g_2 \ast  \wt \g_4 \ast \ldots \ast \wt \g_{2k}  $,   in which the first arc $  \wt \g_1 $ is an integral curve of $\bT$  with a parameterisation of the form $\Phi^{\bT}_{t}(x_o') $ with $t \in [0,  \sum_{\ell = 0}^k \s_{2 \ell  +1}]$. Since by assumption $ \sum_{\ell = 0}^k \s_{2 \ell  +1} = T$, such integral curve  coincides  with the smooth  curve $\g$.  Moreover, by construction,  all  arcs  $\wt \g_{2 \ell}$ are  integral curves of vector fields of the form $\l^\a_{2 \ell}W_\a^{(I)\t_{2\ell}} $ with   $\t_{2 \ell} = \sum_{j = \ell}^k \s_{2 j + 1}$,  and the composition of these arcs starts from $y_o$ and ends at the same final point of $\wc \g$,  as claimed. 
\end{pf}
 The property established in Lemma \ref{lemma4.10} motivates the following 
\begin{definition}   \label{def76} Let  $ \f =  \f_1 \ast \f_2 \ast \ldots \ast  \f_{r} $ be  a piecewise oriented regular curve in an open subset $\cU \subset \cM$, on which  there exists a  vector field $\bT$ such that  $\cT_x =  \bT_x{\!\!\mod \!\cD^I_x}$ for all $x \in \cU$,    and a set of $\bT$-adapted generators $(W_A) = (W^{(I)}_\a,W^{(II \setminus I)}_B)$ for $(V^{II}|_{\cU}, \cD^{II}|_{\cU})$.\par
\noindent (1)
 In case  $\f$ is a non-negative $\cD$-path,  consisting  of an odd number $r = 2 k +1$  of  regular arcs satisfying the following conditions: 
 \begin{itemize}[leftmargin = 20pt]
 \item  the   odd arcs $ \f_{2 \ell + 1}$  are non-trivial integral curves of  $\bT$, 
 \item    the even  arcs  $ \f_{2 \ell}$ are integral curves of vector fields   $\l^\a_{2\ell} W^{(I)}_\a$ for some $\l^\a_{2\ell} \in \bR$, 
 \end{itemize}  
  the curve   is called a  {\it stepped (non-negative) $\bT$-path}. The widths    $ \s_{2 \ell +1}$ of the ranges $ [0, \s_{2 \ell +1}]$ for   the parameters of the maps $\f_{2 \ell + 1}(t) = \Phi^\bT_t(x_{2 \ell +1})$ which parameterise  the odd arcs,  are  called  {\it $\bT$-lengths} of  $\f$. The positive
 real numbers $\t_{ \ell} \= \sum_{j = \ell}^k \s_{2 j + 1}$  are the {\it  $\bT$-depths of $\f$ (relative to the final point of $\f$)}.   \par
\noindent (2)
 In case  each regular arc  $\f_i$   is an  integral curve of a vector field of the form 
 \beq \label{defs} (\l^\a_{2\ell} W^{(I)}_\a)^{\t_i} \= \Phi^\bT_{\t_i *}(\l^\a_{2\ell} W^{(I)}_\a)\eeq
  for an ordered set of real numbers   $(\t_1, \ldots, \t_r)$  satisfying  
  $$ \t_1 > \t_2 > \ldots > \t_r > 0\ ,$$ 
 the  curve is called {\it  a surrogate for a  stepped $\bT$-path with  $\bT$-depths  $(\t_1, \ldots, \t_r)$} or,  shorter,   {\it a $\bT$-surrogate}.  In this case, the ordered set  $(\t_1, \ldots, \t_r)$
 is called the  {\it  tuple of $\bT$-depths} of the $\bT$-surrogate (see Fig.\ 3).
\end{definition}
%
 \centerline{
\begin{tikzpicture}
\draw[<->, line width = 0.7] (4,5) to (4,2.5) to (12.5,2.5);
\draw[->, line width = 0.7]  (4,2.5) to  (2,0.5);
\node at  (12, 2.3) { \tiny$t$};
\node at  (3.8, 4.8) {\tiny $w$};
\node at  (2.3, 0.5) {\tiny $q$};
\draw [line width = 0.8, blue, densely dotted](4.7, 2.6)  to  (5.2, 2.6)    ; 
\draw [line width = 0.7, blue] (5.2, 2.6)  to (5.2, 3.3)    ; 
\draw [line width = 0.8, blue, densely dotted](5.2,3.3)  to  (6.8, 3.3)   ; 
\draw [line width = 0.7, blue](6.8, 3.3)  to  (6.8, 4.3)   ; 
\draw [line width = 0.8, blue, densely dotted](6.8, 4.3)  to  (11.5, 4.3)   ; 
\begin{scope}[very thick,decoration={
    markings,
    mark={at position 0.3 with {\arrow{latex}}}} ] 
\draw [line width = 1, purple, densely dotted, postaction={decorate} ] (3.2,0.8)  to  [out=-5, in=200] (4, 1.4)   to [out=20, in=200] (9.7, 1.6) ; 
\end{scope}
\begin{scope}[very thick,decoration={
    markings,
    mark={at position 0.5 with {\arrow{latex}}}} ] 
\draw [line width = 1, purple, densely dotted, postaction={decorate} ] (3.2,0.8)  to  [out=-5, in=200] (4, 1.4)   to [out=20, in=200] (9.7, 1.6) ; 
\end{scope}
\begin{scope}[very thick,decoration={
    markings,
    mark={at position 0.7 with {\arrow{latex}}}} ] 
\draw [line width = 1, purple, densely dotted, postaction={decorate} ] (3.2,0.8)  to  [out=-5, in=200] (4, 1.4)   to [out=20, in=200] (9.7, 1.6) ; 
\end{scope}
\node at  (7.5, 1.1) {\color{purple} \tiny  \bf  "undeformed" curve
};
\draw [line width = 1, black]  (9.7, 1.6)   [out=80, in=220] to  (9.9, 2.4)  [out=100, in=280] to     [out=90, in=260]  (10.5, 3.4)   ; 
\draw[fill, black]  (9.7, 1.6) circle [radius = 0.05];
\draw[fill, black]  (9.9, 2.4) circle [radius = 0.05];
\draw[fill, black]  (10.5, 3.4) circle [radius = 0.05];
\node at  (9.9, 1.5) {\color{black} \tiny $y_o$};
\node at  (11, 2.7) { \tiny  {\bf $\bT$-surrogate} };
\node at  (10.7, 3.4) {\color{black} \tiny $y$ };
\begin{scope}[very thick,decoration={
    markings,
    mark={at position 0.6 with {\arrow{latex}}}} ] 
\draw [line width = 0.7, purple, postaction={decorate} ](3.2,0.8)  to  [out=-5, in=200] (4, 1.4)   ; 
\draw [line width = 0.7, red, postaction={decorate}](4,1.4)  to  [out=65, in=260]  (4.4, 2.35)   ;
\draw [line width = 0.7, purple, postaction={decorate} ](4.4,2.35)  to   [out=10, in=190] (6, 2.45)   ; 
\draw [line width = 0.7, red, postaction={decorate}](6,2.45)  to   [out=95, in=260] (6, 3.45)   ;
\draw [line width = 0.7, purple, postaction={decorate}](6,3.45)  to   [out=10, in=190] (10.5, 3.4)   ; 
\end{scope}
\begin{scope}[very thick,decoration={
    markings,
    mark={at position 0.3 with {\arrow{latex}}}} ] 
\draw [line width = 0.7, purple, postaction={decorate}](6,3.45)  to   [out=10, in=190] (10.5, 3.4)   ; 
\end{scope}
\draw[fill, purple]  (4, 1.4) circle [radius = 0.05];
\draw[fill, purple]  (4.4, 2.35) circle [radius = 0.05];
\draw[fill, purple]  (6, 2.45) circle [radius = 0.05];
\draw[fill, purple]  (6, 3.45) circle [radius = 0.05];
\node at  (7.2, 2.7) { \tiny \color{purple} {\bf Stepped $\bT$-path}};
\draw[fill, purple]  (3.2, 0.8) circle [radius = 0.05];
\node at  (3.35, 0.66) {\color{purple} \tiny $x_o'$};
\node at (5, 0) {\color{purple} \tiny \it Integral curves of $\bT$};
\draw[->, line width = 0.3, purple] (5, 0.1) to [out = 90, in= 270] (9, 1.35);
\draw[->, line width = 0.3, purple] (5, 0.1) to [out = 90, in= 270] (5.5, 2.3);
\draw[->, line width = 0.3, purple] (5, 0.1) to [out = 90, in= 270] (8, 3.44);
\node at (1, 3.5) {\color{red} \tiny \it Integral curves of vector  fields $\bW^{(I)}_\a $};
\draw[->, line width = 0.3, red] (1.1, 3) to [out = 0, in= 180] (4.1, 2);
\draw[->, line width = 0.3, red] (1.1, 3) to [out = 10, in= 170] (5.8, 3);
\node at (11, 0) {\color{black} \tiny \it Integral curves of the  fields $\Phi^{\bT}_{s_i}{}_*(\bW^{(I)}_\a)$};
\draw[->, line width = 0.3, black] (11, 0.1) to [out = 90, in= 0] (10, 2);
\draw[->, line width = 0.3, black] (11, 0.1) to [out = 90, in= 0] (10.5, 3);
 \end{tikzpicture}
 }
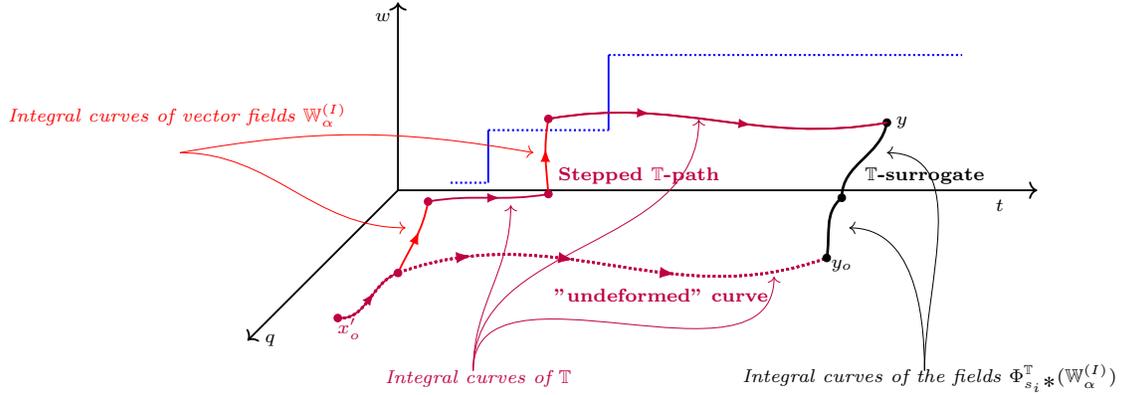
  \centerline{\tiny \bf \hskip 1 cm Fig.\ 3\  The smooth curves constituting  a  stepped $\bT$-path and its  $\bT$-surrogate}
  \ \\[15pt]
  Using these new notions,  the claim of Lemma \ref{lemma4.10} can be stated as follows:  {\it If $\g$ is the  above considered  integral curve  of $\bT$, that  joins $x_o'$ to $y_o$ and whose parameter $t$ runs in $[0, T]$, there is a neighbourhood $\cU$ of $y_o$ such that any of its points,   which can be reached from $x_o'$ via  stepped $\bT$-paths with sum of  $\bT$-lengths equal to  $T$,  can be also reached  starting  from  $y_o$ and moving along   $\bT$-surrogates with $\bT$-depths $(\t_1, \ldots, \t_r)$ with $\t_1 < T$.}\par\smallskip
 Following backwards  the proof of  Lemma \ref{lemma4.10}, one can directly check that  also the inverse of this  claim holds, that is  {\it for any $\bT$-surrogate  $ \f =  \f_1 \ast \f_2 \ast \ldots \ast  \f_{r} $  originating from $y_o$ and with $\bT$-depths  $(\t_1, \ldots, \t_r)$ with $\t_1 < T$,  there exists  a stepped $\bT$-path originating from $x_o'$ and with sum of   $\bT$-lengths  equal to $T$,  for  which $\f$ is the uniquely associated $\bT$-surrogate. }
  Thus, if we denote by  $\Surr(y_o)$ the set  of the points in $\cU$ that can be reached from $y_o$  via $\bT$-surrogates with $\bT$-depths  $(\t_1, \ldots, \t_r)$ with $\t_1 < T$, we conclude  that any point of $\Surr(y_o)$   is  also an endpoint of  a non-negative  $\cD$-path starting from $x_o'$. Thus 
$$\Surr(y_o) \subset \Attain_{x'_o}^{\cC^\o} \subset \Attain_{x_o}^{\cC^\o}\ .$$
   \par
\medskip
  \subsubsection{The $\bT$-surrogate vector  fields}\label{step0}
Let $\cM$, $x_o$, $x_o'$, $y_o$,  $\cU$, $\bT$, $(W_A) = (W^{(I)}_\a, W^{(II \setminus I)}_B)$,  etc. as in the previous subsection and  denote  $N \= \dim \cM$. By  considering a possibly smaller $\cU$,  
 we  may   assume that  
there is a set of coordinates
$ \xi = (x^0,x^1, \ldots, x^{N-1}): \cU \subset \cM \longrightarrow \cV \subset \bR^N$
in which   the  expressions  of $x_o'$, $y_o$  and $\bT$ are 
$$ x'_o \equiv (-T, 0,\ldots, 0)\ ,\qquad y_o \equiv (0, 0, \ldots, 0)\ , \qquad \bT = \frac{\p}{\p x^0}\ .$$
Note  that, 
given a  sufficiently small $\t \in [0, T]$  and a  generator $W^{(I)}_\a$ of the considered set of  $\bT$-adapted  generators, the   vector field 
$W^{(I)\t}_{\a} \= \Phi^\bT_{\t*}(W^{(I)}_\a)$ is equal to   
\beq\label{flusso1}  W^{(I)\t}_\a = W^{(I)}_\a + \sum_{k = 1}^\infty \frac{(-1)^k}{k!} \underset{\text{$k$-times}} {\underbrace{[\bT, [\bT, \ldots, [\bT, W^{(I)}_\a]\ldots]]}}\t^k\ .\eeq
%
This formula can be easily checked as follows.  As a direct consequence of the  definition of Lie derivative, 
the $\t$-parameterised family of vector fields $W^{(I)\t}_{\a}$ is the unique solution to the differential problem $\frac{d W^{(I)\t}_{\a}}{d\t} =  -  [\bT, W^{(I)\t}_\a]$ with the initial condition $ W^{(I)\t= 0}_{\a} =  W^{(I)}_\a$, 
Hence, considering the  coordinate expressions $W^{(I)}_{\a} = W^i_\a(x)  \frac{\p}{\p x^j}$.    
 $W^{(I)\t}_{\a} = \wt W^j_{\a}(x,  \t) \frac{\p}{\p x^j}$,  $\bT = \frac{\p}{\p x^0}$,  the differential problem characterising the family $W^{(I)\t}_{\a}$  corresponds  to   the   differential problem   on  their coordinate components given by 
$$ \frac{\p \wt W^j_{\a}(x,  \t)}{\p \t} = -  \frac{\p \wt W^j_{\a}(x,   \t)}{\p x^0}\ ,\qquad  \wt W^j_\a(x,  0)  = W^j_\a(x)\ .$$
 The  unique solution to this problem is  $   \wt W^j_{\a}(x,  \t) = W^j_\a(x^0 -\t, x^1, \ldots, x^{N-1})$ or, equivalently,  the sum of the power series 
 \beq \label{cris} \wt W^j_{\a}(x,  \t)  =  W^j_\a(x) + \sum_{k = 1}^\infty \frac{(-1)^k}{k!}  \frac{\p^k W^j_\a}{(\p x^0)^k}\bigg|_{(x) } \t^k \ .\eeq
Since these are the  coordinate components of  the  vector field  \eqref{flusso1},   the claim follows.\par
\smallskip
We now observe that, for any  $k \geq 0$,  the vector field $ \underset{\text{$k$-times}} {\underbrace{[\bT, [\bT, \ldots, [\bT, W^{(I)}_\a]\ldots]]}}$  is in the secondary   distribution  $(V^{II},  \cD^{II})$ and has therefore the  form 
\begin{equation} 
 \underset{\text{$k$-times}} {\underbrace{[\bT, [\bT, \ldots, [\bT, W^{(I)}_\a]\ldots]]}} =  \sum_{B= 1} ^{{\bf m}}A_{\a; k}^B  W_B
\end{equation}
 for some appropriate real analytic  functions $A_{\a; k}^B$.  More precisely, since each vector field $W^{(I)}_\a$ has the form  $W^{(I)}_\a= W_{0(a)j}$ for  some   $0 \leq  a \leq \nu$,  
 $1 \leq j \leq R_{a}$,   by  the properties of the $\bT$-adapted generators given in  Lemma \ref{lemmone}, we have that 
\begin{multline} \label{flusso2}
 \underset{\text{$k$-times}} {\underbrace{[\bT, [\bT, \ldots, [\bT, W^{(I)}_\a]\ldots]]}} =
   \underset{\text{$k$-times}} {\underbrace{[\bT, [\bT, \ldots, [\bT, W_{0(a)j}]\ldots]]}} = \\
 = \left\{\begin{array}{ll}W_{k(a)j}  &   \text{if}\ 1\leq k \leq a\ ,
 \\[15 pt]
  \sum_{B = 1}^{{\bf m} } A^B_{\a; k}   W_B  & \text{if} \ k >  a\ .
  \end{array} \right.
\end{multline}
Note that, if the $\bT$-adapted generators $W_B$ are not pointwise linearly independent on a dense open  subset of $\cU$,  there might be several choices for the real analytic functions $A_{\a; k}^B$, $k > a$, defined over the whole $\cU$. In order to eliminate  this ambiguity,  we fix a choice for the  $A_{\a; a+1}^B$ and we   define  
\beq \label{iterateddef} A_{\a; a + r }^B \= \bT(A_{\a; a + r -1}^B) + \sum_{\smallmatrix C = ``k(b)j'' \text{with}\  k \leq b\\
\b \equiv ``0(b)j'' \endsmallmatrix}  A_{\a; a + r -1}^C A^B_{\b; k+ 1}\ , \qquad r \geq 2\ .\eeq
In fact, these functions  might be used as  components in the expansion  \eqref{cris} for any $k =  a + r$, $r \geq 2$, since  the  following   relations hold:
\begin{multline} 
 \underset{\text{$(a+r)$-times}} {\underbrace{[\bT, [\bT, \ldots, [\bT, W^{(I)}_\a]\ldots]]}}  =  [\bT, A_{\a; a + r -1}^B W_B] = \\[-15pt]
= 
\bigg( \bT(A_{\a; a + r -1}^B) + \sum_{\smallmatrix C = ``k(b)j'' \text{with}\  k \leq b\\
\b \equiv ``0(b)j'' \endsmallmatrix} A_{\a; a + r -1}^C A^B_{\b; k + 1}\bigg) W_B\ .
\end{multline}
From \eqref{flusso1}  and \eqref{flusso2},  each $\bT$-pushed forward vector field 
$W^{(I)\t}_\a$   
has  the form  
\begin{multline}  \label{flusso3}  W^{(I)\t}_\a = W^{(I)\t}_{0(a) j} 
= \sum_{B = 1}^{{\bf m}} A_\a^{B}(\t)  W_B \ ,\\
 \text{where} \quad A_\a^{B}(\t) \= \left\{\begin{array}{ll}\frac{{(-1)}^k}{k!} \t^k +  \sum_{k' = a+1}^\infty  \frac{(-1)^{k'}}{k'!} A_{\a; k'}^B   \t^{k'} & \text{if} \   W_B = W_{k (a) j}\ ,\\[10pt]
 \sum_{k' = a+1}^\infty  \frac{(-1)^{k'}}{k'!} A_{\a; k'}^B  \t^{k'} & \text{otherwise}\ ,
 \end{array}\right.\end{multline}
and    is therefore a vector field   in  $(V^{II}, \cD^{II})$. The convergence of $  \sum_{k' = a+1}^\infty  \frac{(-1)^{k'}}{k'!} A_{\a; k'}^B  \t^{k'}$ (and hence the property that  $A_\a^{B}(\t)$ is a well defined real analytic function of $x \in \cM$ and $\t \in [0, \d]$ with $\delta$ small) can be checked  as follows. For any sufficiently small open subset $\cV \subset \cU$, such that all real analytic functions $A_{\a; k}^B$,  $k \leq a +1$ (here,  $a$ is  the integer occurring as index of   $W_\a = W_{0(a)j}$) are identifiable with restrictions to $\bR^N$ of complex analytic functions on some open subset  $\wh \cV \subset \bC^N$, one can determine two constants $C_{\wh \cV}\geq  1$, $r_{\wh \cV}> 0$ such that for any $\a$,  $B$ and  $k \leq a + 1$ 
\beq \label{est} \sup_{x \in \cV}  \left|A_{\a; k}^B|_x \right|< C_{\wh \cV}\ ,\qquad  \sup_{x \in \cV}   \left|\bT(A_{\a; k}^B)\bigg|_x \right| = \sup_{x \in \cV}   \left|\frac{\p A_{\a; k}^B}{\p x^0}\bigg|_x \right|< \frac{C_{\wh \cV}}{r_{\wh \cV}}
\eeq 
(the second estimate  in \eqref{est} is a consequence of the Cauchy integral representation formula for derivatives of holomorphic functions). Combining these estimates with the iterative definition \eqref{iterateddef}, we get that on any sufficiently small open subset $\cV \subset \cU$, the functions $A_{\a; a + 1 + r}^B $, $r \geq 0$,  satisfy 
$$\sup_{x \in \cV}  \left|A_{\a; a + 1+ r}^B|_x \right|< C_{\wh \cV}\left( {\bf m} C_{\wh \cV} + \frac{1}{r_{\wh \cV}}\right)^r\ .$$
This   implies that the series $  \sum_{k' = a+1}^\infty  \frac{(-1)^{k'}}{k'!} A_{\a; k'}^B  \t^{k'}$  converges uniformly on compacta of $\cU \times [0, \d]$ and its sum is a real analytic function of $x \in \cM$ and  $\t \in [0, T]$, as claimed.\par
\smallskip
We conclude stressing   the fact  that all results of previous discussion  are true {\it  for any $\t \in [0,T]$} with no  condition on  the smallness of $\t$. Indeed, in the previous argument the assumption that $\t$ is  sufficiently small was used only to have the possibility of expanding   in power series of $\t$. Decomposing   $[0, T]$ into a finite union of sufficiently small  intervals,  the same arguments allow  to  check that the claim holds for any $\t \in [0, T]$.
 \par \smallskip
These  remarks  motivates  the following
 \begin{definition}\label{surrogatefields}  Given a set of $\bT$-adapted generators  $(W_A ) = (W_\a^{(I)}, W^{(II \setminus I)}_A)$   for $(V^{II}|_{\cU}, \cD^{II}|_{\cU})$, for a given  $W_\a^{(I)}$ and  $\t >0$  the corresponding  vector field  
 \beq 
 \begin{split} & W_\a^{(I)\t} \= \Phi^\bT_{\t*}( W_\a^{(I)})\ ,\\
 \end{split}
 \eeq
 is called {\it elementary $\bT$-surrogate  field of   $\bT$-depth $\t$}.
A {\it $\bT$-surrogate field of $\bT$-depth $\t$} is any linear combinations with constant coefficients of elementary  $\bT$-surrogate fields of  $\bT$-depth $\t$ (i.e. a vector field of the form $\Phi^\bT_{\t*}( \lambda^\a W_\a^{(I)})=  \lambda^\a \Phi^\bT_{\t*}(  W_\a^{(I)})$). 
 \end{definition}
We may therefore say that  the $\bT$-surrogates that  start from $y_o$ and with   $\bT$-depths $(\t_1, \ldots, \t_r)$,  $\t_1 < T$, are 
 piecewise regular curves, whose regular arcs are  integral curves of $\bT$-surrogate fields of $\bT$-depths $\t_1$, $\t_2$, etc..  As we pointed out,  any $\bT$-surrogate field is in $(V^{II}, \cD^{II})$ and, consequently,  {\it any $\bT$-surrogate  that starts from $y_o$ is a $\cD^{II}$-path},  as  we  announced in Remark \ref{remark56}.
\par
\bigskip
\addcontentsline{toc}{section}{Part III}
\centerline{\large \it PART III}
\ \\[-45pt]
\section{The $\bT$-surrogate  fields are generators of  the secondary  distribution}
\label{sect8}
As in Part II,  in  this  section    we assume that  $( \cD, \cD^I,  \cT\= \bT{\!\!\mod \!\cD^I})$   is a real analytic rigged distribution on  a real analytic Riemannian manifold $(\cM, g)$.  We    denote by  $x'_o, y_o$   two fixed  points of $\cM$, joined one to the other by an integral curve $\g(t) = \Phi^\bT_t$, $ t \in [0, T]$, of a  $\bT$ such that  $\cT_x =  \bT_x{\!\!\mod \!\cD^I_x}$ for all $x$ in  a sufficiently small neighbourhood $\cU$ of $y_o$  and  $x_o'$.  We  also assume that on $\cU$ there is   a set of $\bT$-adapted generators $(W_A) = (W^{(I)}_\a,W^{(II \setminus I)}_B)$ for $(V^{II}, \cD^{II})$,  as defined in \S\ref{sect63}.
The next theorem paves the way to the proof of  the main results of Part III.\par
\smallskip
 \begin{theo}  \label{lemma38} Let $\r$ be  a fixed real number,  $\r \in (0,1)$.  There exists a  sufficiently small  $\ve_o \in  (0, T)$  such that, for any $\o \in (0, \ve_o)$,    there is a set of   real numbers $\t_{\ell (a) j} \in [\r \o, \o]$ in bijection with the  adapted generators $(W_A) = (W_{\ell(a)j})$ such that 
%
 \begin{itemize}[leftmargin = 20pt]
 \item[(1)] if   $``\ell(a)j" < ``\wt \ell (\wt a) \wt j"$ in the lexicographic order, then  $ \t_{\ell(a)j} > \t_{\wt \ell(\wt a) \wt j}$; 
 \item[(2)] the ordered set of $\bT$-surrogate fields  $ \bigg( \bW_{\ell(a)j} \=   \Phi^\bT_{\t_{\ell(a) j}*}(W_{0(a)j})\bigg)$
 is a set of  generators for   $(V^{II}|_{\cU}, \cD^{II}|_{\cU})$.
 \end{itemize}
 \end{theo}
 \begin{pf} 
 Due to  \eqref{flusso1} and \eqref{flusso2}, if  $\ve_o>0$ is sufficiently small,  for any $s \in (0, \ve_o)$  the   $\bT$-surrogate field 
  $ W_{0(a)j}^{s} \= \Phi^\bT_{s*}(W_{0(a)j})$  has the form
   \begin{multline} \label{the61} W_{0(a)j}^{s} = W_{0(a )j } - s W_{1(a )j } + s^2 \frac{1}{2} W_{2(a )j } - \ldots +(-1)^{a } s^{a } \frac{1}{a !} W_{a (a )j } + \\
 + s^{a +1} Y_{((a ) j |s)}\end{multline}
 for an appropriate vector field $ Y_{((a)j |s)}$  in  $V^{II}|_\cU$.  We recall that, from the results of \S \ref{step0}, if $W_{0(a)j} $ is the vector field denoted also as $W_{0(a)j} = W^{(I)}_\a$,   the vector field $Y_{((a)j|s)}$ can be expanded in terms of the vector fields $W_A$ as
\beq \label{the61*} Y_{((a)j |s)} = \wh A^B_{(a)j |s} W_B\qquad\text{with}\qquad    \wh A^B_{(a)j |s} (x)=  \sum_{r = 0}^\infty  \frac{(-1)^{a +1+r}}{(a+1+r)!} A_{\a; a + 1 + r }^B(x)  s^{r} \eeq
where the $A_{\a; a +1 +  r}^B$ are the real analytic functions defined inductively in \eqref{iterateddef} and the $ \wh A^B_{(a)j |s}$ are real analytic functions. \par
 \smallskip
  Consider a fixed pair of indices $a = a_o$,  $j = j_o$ and denote by   $v: \cU \to \cD^{II}|_{\cU}$    the vector field defined by 
$$v(x) \= W_{0(a_o)j_o}^{s}|_{x}  \in \cD^{II}|_{x}\ ,\qquad x \in \cU$$
 and    by $\wt v$ the associated field of   projections onto   $\mathbb V_x \= \cD^{II}|_{x}\bigg/\left\langle  W_{\ell(a)j}|_{x}\ , a \neq a_o, j  \neq j_o\right \rangle $:
 $$\wt v(x) \= v(x)  \hskip -10 pt\mod \left\langle  W_{\ell(a)j}|_{x}\ , a \neq a_o, j \neq j_o\right \rangle\ .$$
Being  the $W_B|_x$  generators of the generalised distribution $(V^{II}, \cD^{II})$,   at each  $x \in \cU$ the projections in the quotient $\mathbb V_x$ of  the vectors 
 \beq \label{cris*} W_{0(a_o)j_o}\big|_{x}\ ,\qquad W_{1(a_o)j_o}\big|_{x} \ ,\qquad \frac{1}{2!} W_{2(a_o)j_o}\big|_{x}\ ,\qquad \ldots\qquad   \frac{1}{a_o!} W_{a_o(a_o)j_o}\big|_{x}\eeq
 constitute a  set of (possibly linearly dependent) generators   for the space $\mathbb V_x$.  We may   therefore   expand the field  $\wt v$ (i.e. the field of projections of $v$ onto the quotient spaces $\bV_x$, $x \in \cU$)  in terms of the set  $\cC$, made of the fields of projections   
\beq \label{gen11*}    \frac{1}{k!}\wt{W_{k (a_o)j_o} }\= \frac{1}{k!}W_{k (a_o)j_o} \hskip -7 pt\mod \left\langle  W_{\ell(a)j}|_{x}\ , a \neq a_o, j \neq j_o\right \rangle\ ,\qquad 0 \leq k \leq a_o\ ,\eeq
 with  real analytic coefficients.    
  From  \eqref{the61} and \eqref{the61*},   one gets that  a $(a_o+1)$-tuple of real analytic  coefficients  which can be used to expand $\wt v$ in terms of the set $\cC$ of the  fields of projections  \eqref{gen11*} is   
 \begin{multline} \label{the63} \bigg(1, -s, s^2, \ldots, (-1)^{a_o} s^{a_o}\bigg) + s^{a_o + 1} \bigg(\l_0(s, x), \l_1(s, x), \ldots, \l_{a_o}(s, x)\bigg) =\\
 =  \bigg(1 + s^{a_o+1} \l_0(s, x), s \big(-1 + s^{a_o} \l_1(s, x)\big) , s^2\big(1 + s^{a_o-1} \l_2(s, x)\big) , \ldots\\
 \ldots , s^{a_o}\big((-1)^{a_o } + s \l_{a_o}(s, x)\big)\bigg) \end{multline}
where   $ \Lambda(s, x) =   (\l_0(s, x), \l_1(s, x), \ldots, \l_{a_o}(s, x)) $  denotes    components  of the field of the projections of the vectors $Y_{((a_o)j_o|s)}|_{x}$ in  the spaces $\mathbb V_x$, $x \in \cU$,  in terms of the generators \eqref{gen11*} and it is uniquely determined by the analytic functions  $\wh A^B_{(a)j |s}$.  
\par
Consider now a real number $\o\in (0, \ve_o)$ and,  for any   $0 \leq\ell \leq a_o$,  let $\t_\ell$  be 
  \beq\label{thetau}  \t_\ell \=  \o \s_\ell\qquad \text{with }\qquad \s_\ell \= 1 -  \ell  \frac{1  - \r }{a_o}\ .\eeq
By construction, we have that 
  $\ell < \ell'$ implies  $\t_\ell > \t_{\ell'}$ and that  $\t_0 = \o$, $\t_{a_o} = \r \o$.
 Then, for any $0 \leq \ell \leq a_o$,  let   $v_\ell \= W_{0(a_o)j_o}^{\t_\ell}: \cU \to \cD^{II}$ and denote by 
 $\wt v_\ell$  the corresponding field of projections onto the quotients $\bV_x$, $x \in \cU$. By \eqref{the63}, for any $x \in \cU$,  the equivalence class  
  $\wt v_\ell(x)$ is a linear combination of the above defined generators for $\mathbb V_x$.  
  The coefficients of such linear combinations for $0 \leq \ell \leq a_o$  are given by the entries of the matrix
\begin{multline} \label{semi-vandermond} A(\o, \r, x) = \\
=  {\tiny \left( \begin{array}{ccccc} 1 + \t_0^{a_o+1} c_{0;0} & \t_0 (-1+  \t_0^{{a_o}} c_{0;1}) & \t_1^2 (1   + \t_0^{{a_o}-1}  c_{0;2}) & \ldots & \t_0^{{a_o}}((-1)^{a_o}   + \t_0 c_{0;{a_o}}) \\
 1 + \t_1^{a_o+1} c_{1;0} & \t_1 (-1+  \t_1^{{a_o}} c_{1;1}) & \t_1^2 (1   + \t_1^{{a_o}-1}  c_{1;2}) & \ldots & \t_1^{{a_o}}((-1)^{a_o}   + \t_1 c_{1;{a_o}}) \\
\vdots & \vdots & \vdots& \ddots & \vdots\\
 1 + \t_{a_o}^{a_o+1} c_{a_o;0} & \t_{a_o} (-1+  \t_{a_o}^{{a_o}-1} c_{a_o;1}) & \t_{a_o}^2 (1   + \t_{a_o}^{{a_o}-1}  c_{a_o;2}) & \ldots & \t_{a_o}^{{a_o}}((-1)^{a_o}   + \t_{a_o} c_{a_o;{a_o}}) \\
\end{array}\right) }\ .\end{multline}
where  the $c_{\ell; \ell'} = c_{\ell; \ell'} (\o, x)$ are the components of the vectors   
$$(c_{\ell;0}, \ldots, c_{\ell;a_o})  \=\big( \l_0(\t_\ell,x), \ldots,  \l_{a_o}(\t_\ell,x)\big) = \Lambda(\t_\ell,x)\  , \qquad 0 \leq \ell \leq a_o\ .$$
Since each $\t_\ell$ has the form $\t_\ell = \o \s_\ell$  where   $\s_\ell$   is a real number independent on $\o$   and  in the interval $[\r, 1]$, and since all  terms $c_{\ell ;\ell'} = \l_{\ell'}(\t_\ell,x)$  have a common finite upper bound for small $\omega$, we infer  that    the terms $ \t_\ell^{a_o + 1 -\ell'}  c_{\ell ;\ell'}$, $0 \leq \ell' \leq a_o$,    appearing in the  entries of  $A(\o,\r,  x)$,  tend uniformly in  $x \in \cU$  to  $0$ when $\o \to 0$. Therefore there is a function   $g(\o, x)$  which tends uniformly to $0$ on $\cU$ for $\o \to 0$ such that 
\begin{multline} \label{semi-vandermond-1*}\det A(\o, \r, x) = \\
=  \det{\tiny \left( \begin{array}{cccc} 1 + \t_0^{a_o+1} c_{0;0} & \t_0 (-1+  \t_0^{{a_o}} c_{0;1}) & \ldots & \t_0^{{a_o}}((-1)^{a_o}   + \t_0 c_{0;{a_o}}) \\
 1 + \t_1^{a_o+1} c_{1;0} & \t_1 (-1+  \t_1^{{a_o}} c_{1;1}) & \ldots & \t_1^{{a_o}}((-1)^{a_o}   + \t_1 c_{1;{a_o}}) \\
\vdots & \vdots& \ddots & \vdots\\
 1 + \t_{a_o}^{a_o+1} c_{a_o;0} & \t_{a_o} (-1+  \t_{a_o}^{{a_o}-1} c_{a_o;1}) & \ldots & \t_{a_o}^{{a_o}}((-1)^{a_o}   + \t_{a_o} c_{a_o;{a_o}}) \\
\end{array}\right) } = \\
=   \o^{\frac{(a_o+1)a_o}{2}}\left(\det\left( \begin{array}{ccccc} 1  &- \s_0 & (-\s_0)^2 & \ldots & (-\s_0)^{a_o}\\
 1 & -\s_1 & (-\s_1)^2 & \ldots & (-\s_1)^{a_o} \\
\vdots & \vdots & \vdots& \ddots & \vdots\\
 1 & -\s_{a_o} & (-\s_{a_o})^2& \ldots & (-\s_{a_o})^{a_o}
\end{array}\right)  +  g\big(\o, x\big) \right) = \\
=  \o^{\frac{(a_o+1)a_o}{2}} \left( \prod_{0 \leq j' < j \leq a_o}  (\s_{j'} - \s_{j})  +  g\big( \o, x\big) \right) 
\ .\end{multline}
This implies that,   by considering a sufficiently small $ \ve_o > 0$, for  any choice of $0 <  \o  < \ve_o$,  $0 < \r < 1$ and $x \in \cU$,  the  corresponding determinant of the matrix $A(\o, \r, x)$ is  non-zero. 
Hence, since  the
   $\wt v_0$, $\ldots$, $\wt v_{a_o} \in \mathbb V_x$ are determined  from the generators \eqref{gen11*}  via a matrix which is invertible at all points,  it follows that   also  the elements $\wt v_0(x)$, $\ldots$, $\wt v_{a_o}(x)$ are   generators of  $\bV_x$  for any $x \in \cU$. \par
   \smallskip
   The same argument  can be done on a ``larger scale'' ,   considering  $ a_o$ and $j_o$  not as  fixed, but running freely  in $0 \leq a_o \leq \nu$ and $1 \leq j_o \leq R_{a_o}$.  More precisely,  given $\o\in (0, \ve_o)$ and  $\r \in (0, 1)$,   one can at first  select    
  real numbers $\ \t_{\ell (a) j}$ in the interval $[\r \o, \o]$ by a formula similar to  \eqref{thetau} (where the  denominator $a_o$ is replaced by  the cardinality  ${\bf m} $ of  the set of $\bT$-adapted generators) and in such a way that 
   $$ `` \ell (a)j" < `` \wt \ell(\wt a )\wt j" \qquad \Longrightarrow \qquad  \t_{\ell (a) j} > \t_{ \ell' (a') j'}\ .$$ 
  Second, one may  consider the vector fields 
$ \bW_{\ell(a)j} = W_{0(a)j}^{\t_{\ell(a)j}}$
and observe that, by an argument similar to the  previous,  for any sufficiently small $\ve_o$, $\o \in (0, \ve_o)$ and $\r \in (0,1)$, for any  pair $(a, j)$, $0 \leq a \leq \nu$, $1 \leq j \leq R_a$,  and any $x \in \cU$, the projections onto the quotient space  $ \cD^{II}|_{x}\bigg/\left\langle  W_{\ell(b)k}|_{x}\ , b \neq a\ ,k\neq \ell\right \rangle$ of the vectors   
$$    \bW_{0(a)j}|_x\ ,\ \     \bW_{1(a)j}|_x\ ,\ \ \ldots \ \ , \ \     \bW_{a(a)j}|_x$$
are generators for such  quotient, exactly  as the  projections of the vectors  $ W_{0(a)j}|_x$, $W_{1(a)j}|_x$, \ldots,  $W_{a(a)j}|_x$ are. 
Third, one can observe that the last  property implies that the ordered set
$\left(\bW_{\ell(a)j}\right)$
 is  a set of   generators for  $V^{II}|_{\cU}$.  In order  to check this, consider a vector  field $X \in V^{II}|_{\cU}$  and  expand   it into a linear combination of the generators $(W_{\ell(a)j})$ for $V^{II}|_\cU$
 \\
 \resizebox{0.80\hsize}{!}{\vbox{
\begin{align} 
\nonumber X_x = &\l^{0(0) 1}  W_{0(0)1}|_{x} +\\
\nonumber + & \l^{0(0) 2}  W_{0(0)2}|_{x} +\\
\nonumber +  &\ldots +\\
\nonumber +& \l^{0(0) R_{0}}  W_{0(0)R_{0}}|_{x}  + \\[10pt]
\nonumber +& \l^{0(1) 1}   W_{0(1)1}|_{x} +  \l^{1(1) 1}   W_{1(1)1}|_{x} +\\
\nonumber +& \l^{0(1) 2}   W_{0(1)2}|_{x} +  \l^{1(1) 2}   W_{1(1)2}|_{x} +\\
\nonumber + & \ldots +  \\
\nonumber + & \l^{0(1) R_{1}}  W_{0(1)R_{1}}|_{x} +  \l^{1(1) R_{1}}  W_{1(1){R_{1}}}|_{x} +\\[10pt]
\nonumber +& \ldots + \\[10pt]
\nonumber +& \l^{0(\nu) R_{\nu}-1} W_{0(\nu)R_{\nu}-1}|_{x} + \l^{1(\nu) R_{\nu}-1}  W_{1(\nu)R_{\nu}-1}|_{x} + \ldots
 +  \l^{\nu (\nu) R_{\nu}-1}  W_{\nu(\nu)R_{\nu}-1}|_{x}  \\
\nonumber +& \l^{0(\nu) R_{\nu}}  W_{0(\nu)R_{\nu}}|_{x} + \l^{1(\nu) R_{\nu}}  W_{1(\nu)R_{\nu}}|_{x} + \ldots
+  \l^{\nu (\nu) R_{\nu}}  W_{\nu(\nu)R_{\nu}}|_{x}\ . 
\end{align}
}
}
\ \\[-40pt]
\beq
\label{longona*} 
\eeq
\ \\
Projecting both sides of \eqref{longona*}  onto the quotient space  
\beq \label{space}  \cD^{II}|_{x}\bigg/\left\langle  W_{\ell(b)j}|_{x}\ , b \neq \nu\ ,\ j \neq R_{\nu}\right \rangle\ ,\eeq
all terms of the left hand side, with the only exception of those in the last line,  are mapped into the zero equivalence class.  
On the other hand, we know that   the vectors $  \bW_{\ell(\nu) R_\nu}|_x$ project onto a set of generators for  the  quotient space  \eqref{space}. This implies that in  the last line of \eqref{longona*} we may replace the  linear combination of the vectors   $W_{\ell(\nu) R_\nu}|_x$  
by a corresponding linear combination of the vectors $   \bW_{\ell(\nu) R_\nu}|_x$. This new linear combination is uniquely determined by the vector  $X|_x$. 
We have therefore a new  expansion for $X|_x$ of the form\\
 \resizebox{0.75\hsize}{!}{\vbox{
\begin{align} 
\nonumber X_x = &\l'{}^{0(0) 1}  W_{0(0)1}|_{x} +\\
\nonumber + & \l'{}^{0(0) 2}  W_{0(0)2}|_{x} +\\
\nonumber +  &\ldots +\\
\nonumber +& \l'{}^{0(0) R_{0}}  W_{0(10)R_{0}}|_{x}  + \\[10pt]
\nonumber +& \l'{}^{0(1) 1}   W_{0(1)1}|_{x} +  \l^{1(1) 1}   W_{1(1)1}|_{x} +\\
\nonumber +& \l'{}^{0(1) 2}   W_{0(1)2}|_{x} +  \l'{}^{1(1) 2}   W_{1(1)2}|_{x} +\\
\nonumber + & \ldots +  \\
\nonumber + & \l'{}^{0(1) R_{1}}  W_{0(1)R_{1}}|_{x} +  \l'{}^{1(1) R_{1}}  W_{1(1){R_{1}}}|_{x} +\\[10pt]
\nonumber +& \ldots + \\[10pt]
\nonumber +& \l'{}^{0(\nu) R_{\nu}-1} W_{0(\nu)R_{\nu}-1}|_{x} + \l'{}^{1(\nu) R_{\nu}-1}  W_{1(\nu)R_{\nu}-1}|_{x} + \ldots
 +  \l'{}^{\nu (\nu) R_{\nu}-1}  W_{\nu(\nu)R_{\nu}-1}|_{x}  \\
\nonumber +& \mu^{0(\nu) R_{\nu}}  \bW_{0(\nu)R_{\nu}}|_{x} + \mu^{1(\nu) R_{\nu}}  \bW_{1(\nu)R_{\nu}}|_{x} \ldots
+  \mu^{\nu (\nu) R_{\nu}}  \bW_{\nu(\nu)R_{\nu}}|_{x}\ . \end{align}
}
}
\ \\[-40pt]
\beq \label{longona**}\eeq
 \ \\
 for an appropriate choice of the coefficients $\l'{}^{\ell(a)j}$ and $\mu^{\ell(\nu) R_\nu}$. 
We now project both sides of this equality onto the quotient space 
$$ \cD^{II}|_{x}\bigg/\left\langle  W_{\ell(b)j}|_{x}\ , b \neq \nu\ ,\ j \neq R_{\nu}-1\right \rangle$$
and use a similar argument to  infer that  the linear combination   
$$\l'{}^{0(\nu) R_{\nu}-1} W_{0(\nu)R_{\nu}-1}|_{x} + \l'{}^{1(\nu) R_{\nu}-1}  W_{1(\nu)R_{\nu}-1}|_{x} + \ldots
 +  \l'{}^{\nu (\nu) R_{\nu}-1}  W_{\nu(\nu)R_{\nu}-1}|_{x} $$
 can be  replaced -- in a unique way -- by a linear combination of the vectors  $  \bW_{\ell(\nu)R_{\nu}-1}|_{x} $. Iterating this argument (based on  an appropriate  sequence of projections of   the various version of \eqref{longona*}
onto  quotient spaces) we end up with an expansion of any $X|_x$ as a linear combination of the vectors  $    \bW_{\ell(a)j}|_{x}$, as we needed to show.
 \end{pf}
\begin{definition} \label{def82} Let  $\ve_o, \o, \r$  as in Theorem \ref{lemma38} and denote by  $\t_{ \ell(a)j}$ the corresponding real numbers, associated with $\o$ and $\r$ as in the statement of the theorem. The corresponding  set of $\bT$-surrogate fields    $\bigg(\bW_{\ell(a)j} \= \Phi^{\bT}_{\t_{\ell(a)j}*}(W_{0(a), j})\bigg)$ is called  {\it set of adapted $\bT$-surrogate generators for   $(V^{II}|_{\cU}, \cD^{II}|_{\cU})$  with  $\bT$-depths in the interval $ [\r \o,  \o]$}. 
\end{definition}
In what follows,  a set of adapted $\bT$-surrogate generators   will be denoted with a notation of the kind  $(\bW_A)_{1 \leq A  \leq {\bf m}} = ( \bW_{\ell(a)j})$.   Whenever   we need to   specify the interval $[\r \o, \o]$, to which the  $\bT$-depths  belong, we   enrich the   notation  writing   
$\big(\underset{[\r \o, \o]}{\bW_A}\big)_{1 \leq A  \leq {\bf m}}$. 
\par
\bigskip
\section{$\bT$-surrogate leaflets and   ``good points''}
\label{theeigthsection}
\subsection{Leaflets of a distribution  and the  Chow-Rashevski\u\i-Sussmann Theorem} \label{leaflets} In this subsection we  have to make a little  pause in our discussion of rigged distributions and  introduce a few  convenient  definitions and recall some   properties of certain integral submanifolds of a real analytic generalised  distribution. \par 
\begin{definition}  \label{leaflet} Let   $(\wc V, \wc \cD = \cD^{\wc V})$  be a generalised distribution on a manifold $\cM$ of dimension $N$ and $1 \leq M \leq N$.  A {\it  $\wc \cD$-map  of rank $M$ centred at $x_o \in \cM$} is a  regular smooth  map 
$F: \cV \subset \bR^M \to  \cM$  from a neighbourhood $\cV$ of $0_{\bR^M}$  with the following properties: 
\begin{itemize}[leftmargin = 20pt]
\item[(1)] $F(0) = x_o$ and  $F(\cV)$ is an embedded $M$-dimensional submanifold of $\cM$; 
\item[(2)] there  is a set of  local vector fields $X_1, \ldots, X_{\bf m}$ in $\wc V$ of cardinality ${\bf m} \geq M$,  a set of smooth  functions  $\s_\ell: (- \ve, \ve) \subset \bR \to \bR$,  $1 \leq \ell\leq {\bf m}$,   and an ${\bf m}$-tuple $(i_1, \ldots, i_{\bf m})$ of integers $1 \leq i_\ell  \leq  M$,  such that,  for any $(s^i) \in \cV$, the corresponding point  $ F(s^1, \ldots s^M) \in \cM$ is  given by   a  composition of flows, applied   to $x_o$,  of the form
\beq \label{cond1} F(s^1, \ldots s^M) = \Phi^{X_{\bf m}}_{\s_{\bf m}(s^{i_{\bf m}})} \circ \ldots\circ  \Phi^{X_1}_{\s_1(s^{i_1})}(x_o) \ .\eeq
\end{itemize}
The vector fields $X_\ell$  and the point $x_o$ are called {\it generators} and {\it center of the  $\wc \cD$-map}, respectively. An  embedded submanifold, which is the   image of a  $\wc \cD$-map $F$ is called {\it integral leaflet of the distribution} or  just {\it leaflet},  for short.  If a  leaflet is determined by a $\wc \cD$-map $F$, the generators and the center of the map are called  {\it generators and center} of  the leaflet.
\end{definition}
Due to the condition (2),  all points of  a leaflet  are  joined to its center  by a $\wc \cD$-path and are therefore in the same  $\wc \cD$-path connected component. In general,  the number of generators  of  a leaflet is   larger than   the dimension of the leaflet. However,  if $(\wc V, \wc \cD)$ is  regular  and involutive and  $\cS$ is a maximal  integral   leaf  of $\wc \cD$  passing through a point  $x_o$ (and hence  $\dim \cS = \rank \wc \cD$ and $\cS$ coincides  with   the  $\wc \cD$-path connected component of $x_o$  by the proof of  Frobenius Theorem -- see e.g. \cite{Wa}), any leaflet centred at $x_o$ has dimension less than or equal to  $\dim \cS = \rank \wc \cD$. Actually,  it is not hard to construct  a  leaflet centred at $x_o$ of maximal dimension and hence  equal to  an open subset of  $\cS$.   \par
\medskip
The following theorem is a direct consequence of  a few  facts concerning   moduli of rings of  real analytic functions and families of vector  fields (more precisely, it is a consequence of  \cite[Thm. 3.8]{Ma}) and \cite[Thm.5.16 and Cor. 5.17]{AS}), combined with the celebrated Chow-Rashevski\u\i-Sussmann Orbit Theorem 
(for related work, see also  \cite{GSZ1}). In the next definition, given a finite set   $\{Y_1, \ldots, Y_m\}$ of local vector fields defined  on a common  open subset of $\cM$, for any $r \geq 2$ we denote by $Y_{(i_1, \ldots, i_r)}$ the iterated Lie bracket
$$Y_{(i_1, \ldots, i_r)} \= [Y_{i_1}, [Y_{i_2}, [\ldots [Y_{i_{r-1}}, Y_{i_r}]\ldots]]]\ ,\qquad 1 \leq i_\ell \leq m\ .$$
For $r = 1$, we  set $Y_{(i_1)} \= Y_{i_1}$.  The integer $r$ is called  the {\it depth} of the iterated Lie bracket.

\begin{theo}[Chow-Rashevski\u\i-Sussman Theorem for real analytic distributions] \label{Rash} Let   $(\wc V, \wc \cD = \cD^{\wc V})$ be a real analytic generalised distribution and denote by  $(\wc V^{(\text{\rm Lie})}, \wc \cD^{(\text{\rm Lie})})$ the pair given by the family $\wc V^{(\text{\rm Lie})}$  of  all real analytic vector  fields that are finite combinations  of  iterated Lie brackets $Y_{(i_1, \ldots, i_r)}$ of vector fields $Y_j$ in $\wc V$,  and the associated family of spaces  $\wc \cD^{(\text{\rm Lie})}_y \subset T_y \cM$,  spanned by the vector fields in  $\wc V^{(\text{\rm Lie})}$. Then:
\begin{itemize}[leftmargin = 20pt]
\item[(i)] The pair  $(\wc V^{(\text{\rm Lie})}, \wc \cD^{(\text{\rm Lie})})$ is an involutive  real analytic  generalised distribution; 
\item[(ii)] The $\wc \cD$-path connected component of a  point $x_o$ is  the  maximal integral leaf $\cS$  through $x_o$ of  the  bracket generated distribution  $(\wc V^{(\text{\rm Lie})}, \wc \cD^{(\text{\rm Lie})})$ and 
 any such maximal  integral leaf  $\cS$ is an immersed submanifold of $\cM$;  
 \item[(iii)] Any $\wc \cD$-leaflet  is included in  a  unique maximal integral leaf $\cS$ of  $(\wc V^{(\text{\rm Lie})}, \wc \cD^{(\text{\rm Lie})})$ and for  any point $x_o$ of  such integral leaf there exists a neighbourhood $\cU$  such that   $\cS \cap \cU$   is a $\wc \cD$-leaflet of maximal dimension.
\end{itemize}
\end{theo}
\par
We conclude this subsection introducing the following convenient  notion.  
Let   $(\wc V, \wc \cD = \cD^{\wc V})$ be a  real analytic generalised distribution on a manifold $\cM$ and denote by  $(\wc V^{(\text{\rm Lie})}, \wc \cD^{(\text{\rm Lie})})$ the  corresponding  generalised distribution, as defined in Theorem \ref{Rash}.   A {\it  decomposition of  an open subset $\cU \subset \cM$ into  $\wc \cD$-strata of  maximal $\wc \cD$-depth $\mu$} is a finite family of disjoint subset $\cU_0$, $\cU_1$, \ldots, $\cU_p$   such that $\cU = \cU_0 \cup \cU_1 \cup \ldots \cup \cU_p$  and  with 
 the  following    properties:  
\begin{itemize}[leftmargin = 20pt]
\item[(i)] for each $\wc \cD$-stratum $\cU_j$,    all spaces   $\wc \cD^{(\text{Lie})}|_{y} \subset T_y\cM$, $y \in \cU_j$,   have  the same dimension  and all  maximal integral leaves in $\cU$ of   $(\wc V^{(\text{\rm Lie})}, \wc \cD^{(\text{\rm Lie})})$ passing through the points of $\cU_j$ are  entirely  included  in $\cU_j$; 
\item[(ii)] there is an integer $\mu \geq 1$ (called {\it maximal $\wc \cD$-depth}) and a set of integers $1 \leq \mu_j \leq \mu$,  one per each stratum $\cU_j$, such that  for each $\cU_j$ the spaces  $\wc \cD^{(\text{Lie})}|_{y}$, $y \in \cU_j$, are  generated by the values of a finite number of  iterated Lie brackets $Y_{(i_1, \ldots, i_r)}$,  with  $Y_{j_\ell}\in \wc \cD|_{\cU_j}$ and   $r \leq \mu_j$.   
\end{itemize}
We  recall  that, by  \cite[Thm. 3.8]{Ma}, for any   $x_o \in \cM$ there is  a neighbourhood $\cU$ such that  $\wc \cD^{(\text{Lie})}|_{\cU}$ is spanned by  a finite number of  iterated Lie brackets $Y_{(i_1, \ldots, i_p)}$. {\it We claim that for any such open set $\cU$,  a decomposition   into  $\wc \cD$-strata can be  constructed   as follows. }  First let  $\cU_0 \subset \cU$ be   the maximal set of points  $y \in \cU$   for which $\dim \wc \cD^{(\text{Lie})}|_{y}$ is maximal.   Then,  let  $\cU_1$  be   the maximal  set of  points  $y \in \cU \setminus \cU_0$  for which     $\dim \wc \cD^{(\text{Lie})}|_{y}$ has the second  maximal value,   and so on.  By construction and Theorem \ref{Rash}, each maximal integral leaf of $\wc \cD^{(\text{Lie})}|_{y}$ through a point  $y \in \cU_j$ is necessarily entirely  included in  $\cU_j$ and,  being  $\wc \cD|_{\cU}$  spanned by just a finite number of iterated  Lie brackets of  vector fields in $\wc \cD$,  integers  $1\leq \mu_j \leq \mu = \max \mu_j$  for which (ii) holds can directly be determined. \par
\bigskip
\subsection{$\bT$-surrogate leaflets  of  a secondary distribution} \label{surr-leaflets}
  Let  us go back to   the secondary generalised  distribution   $(V^{II},  \cD^{II})$  on the real analytic Riemannian  manifold $(\cM, g)$, determined by   the real analytic  rigged distribution  $(\cD, \cD^I, \cT = \bT{\!\!\mod \!\cD^I})$ considered in the previous section. Let also $x_o$, $x_o'$, $y_o$,  $\cU$, $\bT$, $(W_A) = (W^{(I)}_\a, W^{(II \setminus I)}_B)$,  etc. as  in \S \ref{stepandsurr}  and denote by 
 \beq \label{cond2}
 \begin{split} & W_\a^{(I)\t} \= \Phi^\bT_{\t*}( W_\a^{(I)})\\
 \end{split}
 \eeq
 the elementary $\bT$-surrogate  fields of   $\bT$-depth $\t$ associated with the generators $W_\a^{(I)}$ of $\cD^I$. 
 \begin{definition} 
 A {\it  $\bT$-surrogate map of rank $M$ centred at $y_o$ and with $\bT$-depths in $(0, T)$}  is a  $\cD^{II}$-map $F: \cV \subset \bR^M \to  \cM$ of rank $M$ and center  $ y_o$, whose  generators  are  $\bT$-surrogate vector fields  $X_{\ell} = \lambda^\a_\ell W_\a^{(I)\t_\ell} $, $1 \leq \ell \leq {\bf m} $, with   $\bT$-depths $\t_\ell$  satisfying 
\beq \label{constrained} T >  \t_1 >  \t_2 >  \ldots > \t_{\bf m} > 0\ .\eeq
The images of the  $\bT$-surrogate  $\cD^{II}$-maps are called  {\it $\bT$-surrogate leaflets}. 
\end{definition}
 Since a  $\bT$-surrogate map  has the form \eqref{cond1} with vector fields that are $\bT$-surrogate fields   constrained  by the inequalities \eqref{constrained}, all  points  of a $\bT$-surrogate  leaflet  are joined to the center by a $\bT$-surrogate curve. This implies that any $\bT$-surrogate leaflet is     in $\Surr(y_o)$ and   is therefore also a subset of   $\Attain_{x'_o}^{\cC^\o} \subset \Attain_{x_o}^{\cC^\o}$.  \par
 \smallskip
 In case  $(\cD, \cD^I, \cT)$ is  the rigged distribution  of  a  control system \eqref{controlsystem} and   $\pi^\cQ: \cM :\to \cQ$  is   the standard projection of $\cM = \bR \times \cQ \times \cK$ onto  $\cQ$, 
 it is  very important  to know  under which conditions  the  $\bT$-surrogate leaflets of maximal dimension are mapped  onto open sets of $ \cQ$ by the projection   $\pi^\cQ$. Indeed,   assume that  $x_o$ and $y_o \in \Attain^{\cC^\o}_{x_o}$ project  onto the points $q_o$,  $\bar q_o$ in $ \cQ$, respectively,  and that the  piecewise regular curve  which joins $x_o$  to $y_o$  has a  final regular arc which is an integral curve of the vector field $\bT$. 
 Then,  not only $\bar q_o$ is in $\Reach_{T} ^{\cC^\o}(q_o)$ for some $T> 0$, but   also all  projections of points of a  $\bT$-surrogate leaflet  centred at  $y_o$   are in  $\Reach_{T} ^{\cC^\o}(q_o)$. This means that, if one can establish that any  $\cQ$-projection of a 
 $\bT$-surrogate leaflet of maximal dimension  is  an  open  subset of  $\cQ$,  one has immediately that  the system has the hyper-accessibility property and  the  small time local controllability property around any  point with  the homing property.\par
\smallskip 
Being the $\bT$-surrogate leaflets a special kind of $\cD^{II}$-leaflets,   their    dimensions  is always  less than or equal to the  maximal dimension of the  $\cD^{II}$-leaflets. By   the discussion of the previous section  and Theorem \ref{Rash}, there exists a neighbourhood $\cU \subset \cM$ of $y_o$ which admits a decomposition in $\cD^{II}$-strata and   the  maximal dimension  for  the  $\cD^{II}$-leaflets through $y_o$ is equal to the dimension of the spaces  $\cE|_y \= \cD^{II(\text{Lie})}|_y$  at the points  $y$ of  the $\cD^{II}$-stratum   containing the leaflet.  Since such a dimension    is computable   by a purely algebraic algorithm,  we   have an easy-to-compute upper bound for the   dimensions of the $\bT$-surrogate leaflets.  But, unfortunately,  at this moment we do not have a convenient algorithm to determine the exact value of the  maximal dimensions for the $\bT$-surrogate leaflets in full generality.  \par
 \smallskip
We therefore focus just on the cases in which these maximal dimensions  coincide with  the  dimensions of the  spaces  $\cE|_y = \cD^{II(\text{Lie})}|_y$ of the  $\cD^{II}$-strata   that contain the leaflets. As it is illustrated  in \cite{GSZ2},  in these cases the checking of the openness  of the  $\cQ$-projections  of   the  maximal dimensional $\bT$-surrogate leaflets  reduces to the (simple) computation of   the ranks of the restrictions $\pi^\cQ_*|_{\cE_x}$ of the differential of the projection $\pi^\cQ$  to the spaces   $\cE_x$. If such  ranks  are equal to $\dim \cQ$ for all $x \in \cM$,  the $\cQ$-projection of  $\bT$-surrogate leaflets  are  open subsets of $\cQ$ by the Inverse Function Theorem. \par
\smallskip
  These observations motivate    the following definition and the  subsequent  corollary.  \par
\begin{definition}\label{definitiongood}  Let  $y_o \in \cM$
and $\bT$ a vector field defined on a neighbourhood $\cU$ of $y_o$  such that $\cT_x = \bT_x  {\!\!\mod \!\cD^I_x}$  for any $x \in \cU$. A point $y_o$   is called {\it $\bT$-good}  if there is a  $\bT$-surrogate leaflet centred at $y_o$,  whose dimension is equal to  the  dimension of the  integral leaves   of the bracket generated distribution $\cE$ of the $\cD^{II}$-stratum of $y_o$ (= the maximal possible dimension for  $\bT$-surrogate leaflets centred at $y_o$).\par
 The point $y_o \in \cM$  is called {\it good} if (a) there is  vector field $\bT$, defined on a neighbourhood $\cU$ of $y_o$    with  $\cT_x = \bT_x  {\!\!\mod \!\cD^I_x}$  for any $x \in \cU$, such that $y_o$ is $\bT$-good, and (b) $y_o$ is $\bT'$-good  for any   vector field $\bT' = \f_*(\bT)$ determined by  the vector field $\bT$ in (a) and  a local diffeomorphism $\f: \cU \to \cM$ that maps the maximal integral leaves of   $\cD^I|_{\cU}$ into themselves. 
\end{definition}
\begin{rem} \label{remrem} A geometric interpretation of the goodness property is the following. A point $y_o \in \cM$ is   good if and only if there is a neighbourhood $\cU \subset \cM$ of $y_o$ with the property that  for each positive $\cC^\o$ curve $\h(t)$, $t \in [0, \ve]$, in $\cU$  with final point $y_o = \h(\ve)$, there is at least one vector field $\bT'$  on a neighbourhood $\cV \subset \cU$ of  the  curve  such that the following three conditions hold:  (i) $\cT_x = \bT'_x  {\!\!\mod \!\cD^I_x}$  for any $x \in \cV$,  (ii) up to a re-parameterisation of the curve, $\dot \h(t) = \bT'_{\h(t)}$ for all $t$ and (iii) $y_o$ is $\bT'$-good.  Such equivalence can be checked as follows.  Assume that  $y_o$ is good and consider a vector field $\bT$ on a neighbourhood $\cU$ as in (a) of Definition \ref{definitiongood}.  For a sufficiently small $\ve$,  the curve $\h_t \= \Phi^{\bT}_t(\Phi^\bT_{-\ve}(y_o))$, $t \in [0, \ve]$,  is a positive real analytic curve in  $\cU$ that  ends at $y_o$ and for which (i) -- (iii) hold. If $\h'(t)$ is any  curve in $\cU$ 
that ends at $y_o$, by Remark \ref{remark36}, there exists a local diffeomorphism $\f$ which preserves the integral leaves of $\cD^I|_{\cU}$ and such that the vector field $\bT' = \f_*(\bT)$ satisfies (i) and (ii). Since $y_o$ is a good point, also condition (iii) holds.  The converse can be checked similarly.
\end{rem}
The existence of vector fields  satisfying  (i) -- (iii)  of Remark \ref{remrem} for all positive $\cC^\o$ curves  ending at a good point   implies that whenever $y_o$ is a good point and it is in  $\Attain^{\cC^\o}_{x_o}$,  then  there is at least one  $\bT'$-surrogate leaflet  with center in $y_o$, which  is not only entirely included in $ \Attain^{\cC^\o}_{x_o}$, but also it has the maximal possible dimension. This and the previous discussion    leads immediately to the next corollary,  the  main reason of interest for the  notion of good points. 
\begin{cor} \label{teoremone1} Let $x_o \in \cM$ and assume that $y_o \in  \Attain_{x_o}^{\cC^\o}$ is a good point. Let also $\cU$ be a neighbourhood   of $y_o$ admitting a decomposition into $\cD^{II}$-strata.  Then  $ \Attain_{x_o}^{\cC^\o}$ contains an open subset  (of the intrinsic topology) of the   maximal integral leaf  $\cS^{(y_o)}$  of  $(V^{II(\text{\rm Lie})},  \cD^{II(\text{\rm Lie})})$ throught $y_o$  (\footnote{ We recall that a maximal integral leaf $\cS$  of an involutive distribution   is the image 
$\cS = \imath(\wt \cS)$ of  an injective immersion $\imath: \wt \cS \to \cM$ of a   manifold $\wt \cS$. The images  $\cV = \imath(\wt \cV)$ of  open  subsets $\wt \cV \subset \wt \cS$ are called  {\it  open sets of    the intrinsic topology of $\cS$}. 
}). 
\end{cor}
%
Apparently,  the  notion of ``good point''   involves an infinite number of conditions, one per each local diffeomorphism $\f: \cU \to \cM$ preserving the integral leaves of $\cD^I|_{\cU}$. The following proposition clarifies that this is not the case. 
\begin{prop} \label{prop76}  A point $y_o$ is good if and only if it is $\bT$-good for at least   one    vector field $\bT$  on a neighbourhood $\cU$ of $y_o$   with   $\cT_x= \bT_x  {\!\!\mod \!\cD^I_x}$ for all $x \in \cU$.
\end{prop}
\begin{pf}  In one direction the implication is trivial. Conversely, assume that $y_o$ is $\bT$-good for at least one  vector field as in the statement and  let   $\bT' = \f_*(\bT)$ where $\f$ is a local diffeomorphism preserving the integral leaves of $\cD^I|_{\cU}$.   
Let also   $(W^{(I)}_\a, W^{(II \setminus I)}_B)$ be the $\bT$-adapted generators, that are used  in the construction of the  $\bT$-surrogate generators and of  the $\bT$-leaflets centred at $y_o$.  The set of vector fields  $(W'_A) = (\f_*(W^{(I)}_\a), \f_*(W^{(II \setminus I)}_B))$ are  $\bT'$-adapted generators and  determine $\bT'$-surrogate generators and  $\bT'$-leaflets.  This  implies that if  there is a $\bT$-surrogate leaflet centred at $y_o$ having     the maximal possible  dimension, then the image of such leaflet  under $\f$ is  a $\bT'$-surrogate leaflet with the same property. \end{pf}

From the previous discussion, we see that  for  applications in control theory it is important to have  manageable  criterions which imply  the ``goodness'' of points. In the next subsections, we provide two such criterions.
\par \bigskip
\subsection{The first criterion  for goodness}
\begin{theo}\label{prop93}   Assume that  $\cU$ is a neighbourhood of $y_o$ admitting a decomposition in $\cD^{II}$-strata$(V^{II}|_{ \cU}, \cD^{II}|_{ \cU})$   and denote by  $\cU_{j_o} \subset \cU$   the $\cD^{II}$-stratum containing $y_o$. If the distribution  $(V^{II}|_{ \cU_{j_o}}, \cD^{II}|_{ \cU_{j_o}})$  is regular and involutive near $y_o$, then   $y_o$  is a good point.
\end{theo}
\begin{pf}  By Proposition \ref{prop76}, the proof reduces to show the existence of  a vector field  $\bT$, which is defined on a neighbourhood $\cU$ of $y_o$,   with  $\cT_x= \bT_x  {\!\!\mod \!\cD^I}_x$, $x \in \cU$,  and such that  there is a  $\bT$-surrogate leaflet centred at $y_o$ of $\bT$-depths in  some interval  $(0, \ve) \subset \bR$,  whose dimension is equal to  the dimension of the  integral leaves   of the bracket generated distribution $\cE^{( \cU_{j_o})} = \cD^{II (\text{(Lie)}}|_{\cU_{j_o}}$ of the $\cD^{II}$-stratum $\cU_{j_o}$ of $y_o$.\par
Pick a vector field $\bT$ on a neighbourhood $\cU$  with   $\cT_x= \bT_x  {\!\!\mod \!\cD^I}$, $x \in \cU$,  and assume that $\cU$ is small enough so that  there is a set of $\bT$-adapted generators $(W_A) = (W^{(I)}_\a,W^{(II \setminus I)}_B)$ for $(V^{II}|_{\cU}, \cD^{II}|_{\cU})$.    Denote by $\wt \ve > 0$ a positive number such  that $\Phi^\bT_s(y)$ is well defined for  all $s \in [-\wt \ve, \wt \ve]$ and  all $y$ in a relatively compact neighbourhood $\cV \subset \cU$ of $y_o$. Then, if  $x_o' \= \Phi^\bT_{-\wt \ve}(y_o)$, we have that $\h_t \= \Phi^\bT_{t}(x'_o)$, $t \in [0, \wt \ve]$ is a real analytic $\cD$-path,  which is positive  and ends at $y_o$,  and $\cU$, $x_o'$, $y_o$, $\bT$,  $(W_A) = (W^{(I)}_\a,W^{(II \setminus I)}_B)$ are as in \S \ref{sect8}.\par
\smallskip
 Let  now $\ve_o, \o, \r$  as in Theorem \ref{lemma38} and denote by  $\t_{ \ell (a) j}$ the corresponding real numbers, associated with $\o$ and $\r$ as in the statement of that theorem. We  denote by  $(\bW_A)_{A = 1 , \ldots, {\bf m}} = \bigg(\bW_{\ell(a)j} \= \Phi^{\bT}_{\t_{\ell(a)j}*}(W_{0(a)j})\bigg)$  the associated  set of   adapted $\bT$-surrogate generators  for   $(V^{II}|_{\cU}, \cD^{II}|_{\cU})$  with  $\bT$-depths in the interval $ [\r \o,  \o]$. We recall that  the  $\bT$-depths $\t_A$    are such that 
 $T > \o = \t_1 > \t_2 > \ldots > \t_{\bf m} =  \r \o > 0$.  Let  $F:  \cV \subset \bR^{\bf m} \to \cM$ be the map, which is defined  on a sufficiently small neighbourhood $ \cV$ of $0_{\bR^{\bf m}}$  by 
 \beq \label{9.6} F(s^1, \ldots, s^{\bf m}) = (\Phi^{\bW_{\bf m}}_{s^{\bf m}} \circ \ldots \circ \Phi^{\bW_1}_{s^1})(y_o)\ .\eeq
 In any  coordinate system around $y_o$,  the columns of the Jacobian matrix  $JF|_{0_{\bR^{\bf m}}}$ coincide with  the coordinate components of the vectors $\bW_A|_{y_o}$, $1 \leq A \leq {\bf m}$,  and are therefore linearly independent. By the Inverse Function Theorem and Frobenius Theorem, 
if $ \cV$ is sufficiently small,   the map $F$ has constant rank and  its image is an ${\bf m}$-dimensional embedded submanifold  contained in  the immersed ${\bf m}$-dimensional submanifold   $\cS^{(y_o)}$, which is the maximal integral leaf of $\cD^{II}$  through $y_o$.  By construction, the submanifold $F(\cV)$ is a $\bT$-surrogate leaflet  centred at $y_o$ and   has the same dimension of  $\cS^{(y_o)}$, as  desired.
\end{pf}
\begin{definition}  The points $y_o \in \cM$  satisfying the conditions of Theorem \ref{prop93}  are called {\it good points of the first kind}.
\end{definition}
\begin{example}
 For  the rigged distribution in Example \ref{example72}  {\it any  point $y_o \in \cM$ is  
 a good point of the first kind}.   Indeed,  in that example   $(V^{II}, \cD^{II})$ is a regular 
 and involutive distribution, hence  with only one stratum  $\cM_0 = \cM$  (in fact,  it is generated by the  vector fields   $W_i^{(\ell)}$, $ 1 \leq i \leq m$,   given   in \eqref{KAL1}, which     are  commuting, globally defined and complete).  So,  Theorem \ref{prop93}  holds for any $y_o \in \cM$. 
 \end{example}
\par
\medskip
\subsection{The second criterion for goodness} \label{sect8.3}
Given a point $y_o \in \cM$, let  $\cU$ be a neighbourhood of $y_o$ with  a vector field  $\bT$, with  $\cT_x  = \bT_x {\!\!\mod \!\cD^I_x}$, $x \in \cU$,  and admitting a decomposition into $\cD^{II}$-strata $\cU_j$. Consider a set of $\bT$-adapted generators  $(W_A) = (W^{(I)}_\a, W^{(II \setminus I)}_B)$ for $(V^{II}|_{\cU}, \cD^{II}|_{\cU})$.   For  each generator $W_\a \= W_\a^{(I)}$ of  $\cD^I|_{\cU}$,  consider also  the rigged distribution $(\cD^{(W_\a)}, \cD^{I(W_\a)}, \cT = \bT \mod \cD^{I(W_\a)})$
 that are determined  by the families of $2$- and $1$- dimensional spaces  $\cD^{(W_\a)}_x$,   $\cD_x^{I(W_\a)}$,  $x \in \cU$,  respectively,  given   by 
$$ \cD^{(W_\a)}_x \=  \langle  \bT_x, W_{\a}|_x \rangle\ ,\qquad \cD^{I(W_\a)}_x \=  \langle  W_{\a}|_x \rangle .$$
The following notion is a fundamental ingredient for  our second criterion for goodness.\par
\begin{definition}  The  {\it  secondary sub-distribution  generated by  $W_\a$} is the  secondary generalised distribution $(V^{II(W_\a)}, \cD^{II(W_\a)})$   determined by the rigged distribution  $(\cD^{(W_\a)}, \cD^{I(W_\a)}, \cT = \bT \mod \cD^{I(W_\a)})$.%
\end{definition}
\begin{theo}  \label{theorem714} Let    $y_o \in \cM$ and  $\cU$ a neighbourhood of $y_o$,  admitting a decomposition into $\cD^{II}$-strata, for which there are  a vector field  $\bT$ such that  $\cT_x = \bT_x  {\!\!\mod \!\cD^I_x} $, $x \in \cU$,  and  a set of $\bT$-adapted generators  $(W_A) = (W^{(I)}_\a, W^{(II \setminus I)}_B)$ for $(V^{II}|_{\cU}, \cD^{II}|_{\cU})$.   Denote   by  $\cE^{(\cU_{j_o})}$ the bracket generated distribution of the $\cD^{II}$-stratum $\cU_{j_o}$
containing $y_o$ with ${\bf m} \= \rank \cE^{(\cU_{j_o})}$. Assume also that:
\begin{itemize}[leftmargin = 15pt]
\item[(1)] the  $\cE^{(\cU_{j_o})}$ has stratum depth $\mu_{j_o} = 2$;  
\item[(2)] there exist  an integer $0 \leq p \leq {\bf m} -  1$ and  vector fields $W_{A_i}$, $1 \leq i \leq p$,   $W_{B_j}, W_{B'_j}$,  $1 \leq j \leq{\bf m}-p$, such that (a)   all vectors fields  $W_{A_i}, W_{B'_j}$ are  among  the $\bT$-adapted generators  $(W_A)$, (b)  each   $W_{B_j}$ is in a set of $\bT$-adapted generators of a   sub-distribution  $(V^{II(W_{\b_j})}, \cD^{II(W_{\b_j})})$  for which  $y_o$ is a $\bT$-good point of the first kind,  and  (c)
the ${\bf m}$-tuple
 \beq\label{genn1} \left(W_{A_1}, \ldots,W_{A_p}, Y_{1} \=  \big[W_{B_1}, W_{B'_1}\big], \ldots, Y_{{\bf m} - p}\= \big[W_{B_{{\bf m}- p}}, W_{B'_{{\bf m}- p}}\big]\right)\eeq
 is a set of generators for $\cE^{(\cU_{j_o})}$ around $y_o$. 
 \end{itemize}
 Then $y_o$ is a good point.
\end{theo}
\begin{definition}  The points $y_o \in \cM$  satisfying the conditions of Theorem \ref{theorem714}  are called {\it good points of the second kind}.
\end{definition}
As we will shortly see, Theorem \ref{theorem714}  is  a consequence of the following  stronger general result on  collections  of $\bT$-surrogate fields and their Lie brackets. 
\par
\smallskip
\begin{theo} \label{criterione}    Let   $y_o$,  $\cU$, $\bT$ and  $(W_A) = (W^{(I)}_\a, W^{(II \setminus I)}_B)$ as in Theorem \ref{theorem714}. Denote by  $\cU_{j_o}  \subset \cU$   the $\cD^{II}|_{\cU}$-stratum containing $y_o$ and by  $\cE^{(\cU_{j_o})}$ the bracket generated distribution of $\cU_{j_o}$ with ${\bf m} \= \rank \cE^{(U_{j_o})}$. The point $y_o$ is  good if the following conditions are satisfied:  
\begin{itemize}[leftmargin = 15pt]
\item[(1)] the  stratum depth of the distribution  $\cE^{(\cU_{j_o})}$ is  $\mu_{j_o} = 2$;  
\item[(2)] there is  an integer $0 \leq p\leq {\bf m} -1$ and 
$1 + {\bf m} - p $ sets   of adapted $\bT$-surrogate generators
$$\left(\bW_A^{(\ell)}\right)_{1 \leq A  \leq {\bf m}} \= \left(\underset{[\r_\ell \o_\ell, \o_\ell]}{\bW_A}\right)_{1 \leq A  \leq {\bf m}}\ ,\qquad 0 \leq \ell \leq {\bf m} - p \ ,$$ 
with  $\bT$-depths in disjoint  intervals $ [\r_\ell \o_\ell,  \o_\ell]$ with $\o_{\ell + 1} < \r_\ell \o_\ell$, from which we can extract  $\bT$-surrogate vector fields 
$$\bW^{(0)}_{A_1}, \ldots,\bW^{(0)}_{A_p}\ ,\   \bW^{(1)}_{B_1}, \bW^{(1)} _{B'_1}\ ,\  \bW^{(2)}_{B_2}, \bW^{(2)} _{B'_2}\\ ,\qquad \ldots \qquad \bW^{({\bf m} - p)}_{B_{{\bf m}- p}}, \bW^{({\bf m} - p)}_{B'_{{\bf m}- p}}$$ 
such that the ${\bf m}$-tuple
 \beq\label{genn1*} \left(\bW^{(0)}_{A_1}, \ldots,\bW^{(0)}_{A_p}, \bY_{1} \=  \big[\bW^{(1)}_{B_1}, \bW^{(1)}_{B'_1}\big], \ldots, \bY_{{\bf m} - p}\= \big[\bW^{({\bf m} - p)}_{B_{{\bf m}- p}}, \bW^{({\bf m} - p)}_{B'_{{\bf m}- p}}\big]\right)\ ,\eeq
 is a set of generators for $\cE^{(\cU_{j_o})}$ and with the property  that  each  vector field $\bW^{(j)}_{B_j}$  belongs  to a    sub-distribution  $(V^{II(W_{\b_j})}, \cD^{II(W_{\b_j})})$  for which  $y_o$  is  a $\bT$-good point of the first kind. 
 \end{itemize}
\end{theo}
The proof of this theorem is  delicate and postponed to the next subsection. In this section we  prove that it  implies Theorem \ref{criterione}. For this we first  need  the following   lemma. 
\par
\smallskip
\begin{lem} \label{lemmone-12} In Theorem \ref{criterione} the conditions $\o_{\ell + 1} < \r_\ell \o_\ell$ can be removed, in the sense that  the claim is  true  on a possibly smaller neighbourhood $\cU$ of $y_o$, whenever the condition  (1) is true and there is  an arbitrary  collection of sets   of adapted $\bT$-surrogate generators
$\left(\bW_A^{(\ell)}\right)_{1 \leq A  \leq {\bf m}}$, which are not necessarily satisfying   the  inequalities $\o_{\ell + 1} < \r_\ell \o_\ell$,  provided that all  other conditions in (2) hold. 
\end{lem}
\begin{pf}
We recall that for any choice of the intervals  $ [\r_\ell \o_\ell,  \o_\ell]$, the corresponding  $\bT$-surrogate fields  $\bW^{(0)}_{A_i}, \bW^{(\ell)}_{B_\ell}$, $\bW^{(\ell)}_{B'_\ell} $ are  vector fields of the form 
\beq \label{76} \bW^{(0)}_{A_i} = W_{\a_i}^{\t^{(0)}_{A_i}} \ ,\qquad  \bW^{(\ell)}_{B_\ell} = W_{\b_\ell}^{\t^{(\ell)}_{B_\ell}}\ ,\qquad \bW^{(\ell)}_{B'_\ell} = W_{\b'_\ell}^{\t^{(\ell)}_{B'_\ell}}\ ,\eeq
for  some $W_{\a_i}$, $W_{\b_\ell}$, $W_{\b'_\ell}$ among the generators $(W_\a) =  (W_\a^{(I)})$ of the distribution $\cD^I$ and $\bT$-depths $\t_{A_i}^{(0)}$ and $\t_{B_\ell}^{(\ell)}$,  $\t_{B'_\ell}^{(\ell)}$ in the corresponding intervals $[\r_0 \o_0, \o_0]$ and $[\r_\ell \o_\ell, \o_\ell]$, respectively.  Consider  the real analytic function
\begin{multline}\bF: (0, T]^{2 {\bf m} - p} \subset \bR^{2 {\bf m} - p} \longrightarrow \bR\ ,\\
\bF(\s^1, \s^2, \ldots, \s^p, \s^{p + 1}, \wt \s^{p+ 1},\ldots, \s^{{\bf m} }, \wt \s^{{\bf m} }) \=\\
=  \det \left( \begin{array}{cccccc} W_{\a_1}^{\s^1|1}(y_o) &\ldots &  W_{\a_p}^{\s^p|1}(y_o) &  Y^{\s^{p+1}, \wt \s^{p+1}|1} (y_o)& \ldots &  Y^{\s^{ {\bf m } }, \wt \s^{{\bf m } }|1} (y_o)\\
W_{\a_1}^{\s^1|2}(y_o) &\ldots &  W_{\a_p}^{\s^p|2}(y_o)  & Y^{\s^{p+1}, \wt \s^{p+1}|2}(y_o)& \ldots &  Y^{\s^{ {\bf m }|2 }, \wt \s^{{\bf m } }|2} (y_o)\\
\vdots & \ddots& \vdots & \vdots  & \ddots & \vdots\\
\end{array}
\right) 
\end{multline}
where we use the notation $W_{\a_i}^{\s^i|h}$ and  $Y^{\s^{r}, \wt \s^{r}|h}$ to denote      the  first $m$ coordinate components (in a set of coordinates, which is adapted to the integral submanifold $\cS^{(y_o)}$ of $\cE^{(\cM_{j_o})}$ through $y_o$) of 
the vector fields $W^{\s^i}_{\a_i}$ and $Y^{\s^r, \wt \s^r}_j \= [W^{\s^r}_{\b_j}$, $W^{\wt \s^{r}}_{\b'_j}]$, respectively.  If we assume that  the vector fields \eqref{76} satisfy  all conditions in (2) of Theorem \ref{criterione} with  the possible exception of the  inequalities $\o_{\ell + 1} < \r_\ell \o_\ell$,   the value of $\bF$ at the point 
$$(\s^1_o, , \ldots, \s^p_o, \s^{p + 1}_o, \wt \s^{p + 1}_o, \ldots) \= (\t^{(0)}_{A_1}, \ldots, \t^{(0)}_{A_p}, ,\t^{(1)}_{B_1}, \t^{(1)}_{B'_1}, \ldots)\ .$$
is non-zero.
By real analyticity,  this implies that the  set of all points on which $\bF$  does not vanish   is  open and dense in $ (0, T)^{2{\bf m} - p}$.  It is  therefore possible   to select  a new $({\bf m} - p + 1)$-tuple  of intervals  $[\r'_k\o'_k, \o'_k]$ and   associated  $\bT$-depths $\wt \t^{(0)}_{A_i}$, $\wt \t^{(\ell)}_{B_\ell}$, $\wt \t^{(\ell)}_{B'_\ell}$, such that not only  $\bF$ is non-zero when the $\s^i$ are set equal to  the new values $ (\wt \t^{(0)}_{A_1}, \ldots,  \wt \t^{(0)}_{A_p}, \wt \t^{(1)}_{B_1}, \wt  \t^{(1)}_{B'_1}, \ldots)$, but also with all   inequalities  $\o'_{\ell + 1} < \r'_\ell \o'_\ell$   satisfied. For such a new set of intervals all conditions in   (2) of  Theorem \ref{criterione} are satisfied. 
\end{pf}
We   are  now  ready to prove  Theorem \ref{theorem714}. \par
\smallskip
\noindent{\it Proof of Theorem \ref{theorem714}.}  Since $\cE^{(\cU_{j_o})}$ is regular and of rank ${\bf m}$, the ${\bf m}$-tuple of vectors of $\cE^{\cU_{j_o}}|_{y_o} \subset T_{y_o} \cM_{j_o}$ 
\beq \label{numerino}\left(w_{1} \= W_{A_1}|_{y_o}, \ldots,w_{p} \= W_{A_p}|_{y_o} , \gy_1 \= Y_{1}|_{y_o} , \ldots, \gy_{{\bf m} - p} \= Y_{{\bf m} - p}|_{y_o} \right)\eeq
is linearly independent.  Since  there might be several possible choices  for the tuple of $\bT$-adapted generators satisfying the hypotheses, we  further  assume that  the integer ${\bf m} - p$ is the smallest among those that    occur in such possibilities. 
\par
Consider    $\wt \ve_o > 0$ such that  $\Phi^{\bT}_s(y)$ is well defined for any $s \in [-\wt \ve_o. \wt \ve_o]$ and  any $y$ in a relatively compact neighbourhood $\wt \cU \subset \cU$ of $y_o$.  Pick $T \in (0, \wt \ve_o)$, $\ve_o \in  (0, T)$, $\r\in (0,1)$, $\o \in (0, \ve_o)$ so that all conditions of Theorem  \ref{lemma38} are satisfied with respect to  the sub-distribution   $(V^{II(W_{\b_1})}|_{\wt \cU}, \cD^{II(W_{\b_1})}|_{\wt \cU})$ which contains the vector field $W_{B_1}$. 
Consider also  a set of $\bT$-adapted generators $\bigg(W_{C}^{(W_{\b_1})}\bigg)$  for  $(V^{II(W_{\b_1})}|_{\wt \cU}, \cD^{II(W_{\b_1})}|_{\wt \cU})$, which contains $W_{B_1}$,  and a corresponding set of adapted $\bT$-surrogate generators 
 $\bigg(\bW^{(W_{\b_1})}_A = \bW^{(W_{\b_1})}_{\ell(a)j}\bigg)$    with  $\bT$-depths in the interval $ [\r \o,  \o]$.    Since both $(W^{(W_{\b_1})}_A)$ and $(\bW^{(W_{\b_1})}_A)$ are   generators for the same generalised distribution, there exists a pointwise invertible  matrix $\cA \= \left(\cA^C_A(y)\right)$ with real analytic  entries, such that 
$$\bW^{(W_{\b_1})}_A|_y = \cA_A^C(y)  W^{(W_{\b_1})}_C|_y$$ 
for any $y \in \wt \cU$ (see the proof of Theorem \ref{lemma38} for a construction of such invertible matrix). 
Consider the $1$-dimensional quotient vector space 
$$V = \cE^{(\cU_{j_o})}_{y_o} / \langle w_1, \ldots, w_p, \gy_2, \ldots, \gy_{{\bf m} - p} \rangle\ ,$$
and,  for any surrogate field $\bW^{(W_{\b_1})}_A  $, let us  denote by $\wt \gy^{(A)}_1$ the projection onto $V$ of the vector $\gy^{(A)}_1$  in $\cE^{(\cU_{j_o})}_{y_o}$ given by 
$$\gy^{(A)}_1  =   \big[\bW^{(W_{\b_1})}_{A}, W_{B'_1}\big]\bigg|_{y_o} = \cA_A^C(y_o)\left[ W^{(W_{\b_1})}_C, W_{B'_1}\right]\bigg|_{y_o} - W_{B'_1}\left( \cA_A^C\right) \bigg|_{y_o} W^{(W_{\b_1})}_C\bigg|_{y_o}\ .$$
We claim that,  for at least one index $A$, the equivalence class $\wt \gy^{(A)}_1$ is non zero. This can be checked as follows. If all equivalence classes  $\wt \gy^{(A)}_1$ were trivial,  we should conclude that,  for any index $A$, the vector $\cA_A^C(y_o)\left[ W^{(W_{\b_1})}_C, W_{B'_1}\right]\bigg|_{y_o} $ is a linear combination 
of the vectors $W^{(W_{\b_1})}_C\big|_{y_o}$ and of the vectors $w_1, \ldots, w_p, \gy_2, \ldots, \gy_{{\bf m} - p} $. Since the matrix  $\cA_A^C(y_o)$ is invertible, this would imply  that each  vector $\left[ W^{(W_{\b_1})}_C, W_{B'_1}\right]\bigg|_{y_o} $ is a linear combination of 
 the vectors $W^{(W_{\b_1})}_C\big|_{y_o}$ and of the vectors $w_1, \ldots, w_p, \gy_2, \ldots, \gy_{{\bf m} - p} $. In particular, also the vector $\gy_	1 = \left[ W_{B_1}, W_{B'_1}\right]\bigg|_{y_o} $ would be such.  Since the tuple \eqref{numerino} is linearly independent and all vector fields $W^{(W_{\b_1})}_C$ are in turn linear combinations of the vector fields $(W_A)$, the above remarks would imply  that we might consider an ${\bf m}$-tuple of vector fields of the form 
  \beq\label{genn2} \left(W_{A_1}, \ldots,W_{A_p}, W_{A_{p+1}},  Y_{2} \=  \big[W_{B_2}, W_{B'_2}\big], \ldots, Y_{{\bf m} - p}\= \big[W_{B_{{\bf m}- p}}, W_{B'_{{\bf m}- p}}\big]\right)\ ,\eeq
 which is (a) still made of vector fields in $\cE^{(\cU_{j_o})}$ (b) is  linearly independent at $y_o$  (hence on a neighbourhood of that point) and 
(c)  is  still a set of generators for $\cE^{(\cU_{j_o})}$. But for such a choice of generators, the integer $p$  would be replaced by the integer $p' = p + 1$. This cannot be because we assumed that the number ${\bf m} - p$ was the smallest possible integer among all choices 
 of generators satisfying the conditions.\par
 This contradiction implies that there exists a $\bT$-surrogate vector field $\bW^{(W_{\b_1})}_{B_1}$ such that the ${\bf m}$-tuple
  \begin{multline}\label{genn3} \bigg(W_{A_1}, \ldots,W_{A_p}, Y'_{1} \=  \big[\bW^{(W_{\b_1})}_{B_1}, W_{B'_1}\big],  \\
  Y_{2} \=  \big[W_{B_2}, W_{B'_2}\big], \ldots, Y_{{\bf m} - p}\= \big[W_{B_{{\bf m}- p}}, W_{B'_{{\bf m}- p}}\big]\bigg)\end{multline}
 is linearly independent at $y_o$  (and hence on a neighbourhood of that point) and it is a set of generators for $\cE^{(\cU_{j_o})}$ near $y_o$. Due to this, with no loss of generality, in  the  hypotheses  of the theorem we may replace  the condition that the set of generators for $\cE^{(\cU_{j_o})}$ has    the form \eqref{genn1} with the new condition that such a set has   the form \eqref{genn3}, where $\bW_{B_1} \= \bW_{B_1}^{(W_{\b_1})}$ is an adapted $\bT$-surrogate generator for the secondary sub-distribution $(V^{II(W_{\b_1})}|_{\wt \cU}, \cD^{II(W_{\b_1})}|_{\wt \cU})$. As above, we may also assume that ${\bf m} - p$ is the smallest integer which might occur in an ${\bf m}$-tuple having the properties of \eqref{genn3}. \par
 \smallskip
 Following essentially the same line of arguments, we may now prove that   there exists a 
  $\bT$-surrogate vector field $\bW_{B'_1}$  in the secondary distribution $(V^{II}|_{\wt \cU}, \cD^{II}|_{\wt \cU})$,  such that the new  ${\bf m}$-tuple
  \begin{multline}\label{gen2} \bigg(W_{A_1}, \ldots,W_{A_p}, Y''_{1} \=  \big[\bW_{B_1}, \bW_{B'_1}\big],  \\
  Y_{2} \=  \big[W_{B_2}, W_{B'_2}\big], \ldots, Y_{{\bf m} - p}\= \big[W_{B_{{\bf m}- p}}, W_{B'_{{\bf m}- p}}\big]\bigg)\end{multline}
   is linearly independent at $y_o$  (and hence on a neighbourhood of that point) and  is a set of generators for $\cE^{(\cU_{j_o})}$ near $y_o$.  Iterating  these  arguments, we  conclude that, 
instead of requiring that  the set of generators for $\cE^{(\cM_{j_o})}$ has    the form \eqref{genn1}  we may assume the existence of  a set of generators   of   the form
  \begin{multline}\label{genN} \bigg(W_{A_1}, \ldots,W_{A_p}, \bY_{1} \=  \big[\bW_{B_1}, \bW_{B'_1}\big],  \\
  \bY_{2} \=  \big[\bW_{B_2}, \bW_{B'_2}\big], \ldots, \bY_{{\bf m} - p}\= \big[\bW_{B_{{\bf m}- p}}, \bW_{B'_{{\bf m}- p}}\big]\bigg)\end{multline}
  where each  $\bW_{B_i}$ is an adapted  $\bT$-surrogate field for a  secondary sub-distribution,  for which   $y_o$ is a good point of the first kind.\par
  \smallskip
  Consider now the linearly independent  ${\bf m}$-tuple of vectors in $\cE^{(\cM_{j_o})}|_{y_o}$ defined by 
$$\left(w_{1} \= W_{A_1}|_{y_o}, \ldots,w_{p} \= W_{A_p}|_{y_o} , \gy''_1 = \bY_{1}|_{y_o} , \ldots, \gy''_{{\bf m} - p} \= \bY_{{\bf m} - p}|_{y_o} \right)$$
and let $\wc T \in (0, \wt \ve_o)$, $\wc \ve_o \in  (0, T)$, $\wc \r\in (0,1)$, $\wc \o \in (0, \wc \ve_o)$ so that all conditions of Theorem  \ref{lemma38} are   satisfied  for  $(V^{II}|_{\wt \cU}, \cD^{II}|_{\wt \cU})$.  This yields the existence of  a set of adapted $\bT$-surrogate generators 
 $\bigg(\bW_A = \bW_{\ell(a)j}\bigg)$ for  $(V^{II}|_{\wt \cU}, \cD^{II}|_{\wt \cU})$ which are related with the $\bT$-adapted generators  $(W_A)$ by means of  a pointwise invertible  matrix $\wc \cA(y) \= \left(\wc \cA^C_A(y)\right)$ with real analytic  entries  and such that 
$$\bW_A|_y = \wc \cA_A^C(y)  W_C|_y$$ 
for any $y \in \wt \cU$. Since $\left( \wc \cA_C^A(y_o) \right)$ is invertible, we may determine $p$ vector fields $\bW_{A_i}$ among the surrogate vector fields of the tuple $(\bW_A)$, such that the vectors 
$$\left(\wc w_1 \= \bW_{A_1}|_{y_o}, \ldots,\wc w_{p} \=  \bW_{A_p}|_{y_o} , \gy''_1 = \bY_{1}|_{y_o} , \ldots, \gy''_{{\bf m} - p} \= \bY_{{\bf m} - p}|_{y_o} \right)$$ 
are linearly independent. It follows that the ${\bf m}$-tuple of vector fields 
 \begin{multline}\label{genN1} \bigg(\bW_{A_1}, \ldots,\bW_{A_p}, \bY_{1} \=  \big[\bW_{B_1}, \bW_{B'_1}\big],  \\
  \bY_{2} \=  \big[\bW_{B_2}, \bW_{B'_2}\big], \ldots, \bY_{{\bf m} - p}\= \big[\bW_{B_{{\bf m}- p}}, \bW_{B'_{{\bf m}- p}}\big]\bigg)\end{multline}
  satisfies all conditions of Lemma \ref{lemmone-12}. By the claim of such a lemma and by  Theorem \ref{criterione}, the conclusion follows. 
  \hfill \qed
  \par
  \medskip
\subsection{Proof of Theorem  \ref{criterione}}
Assume that for any index $1 \leq j \leq {\bf m}  -p$,  any  quadruple $\gq \= (\wt \r, \wt \o,\wt{ \wt \r}, \wt{\wt \o})$ of real numbers with $\wt \r, \wt{\wt \r} \in (0,1)$ and $\wt \o, \wt{\wt \o} \in (0, T)$ and any prescribed arbitrarily small constant $\z > 0$, there is 
a   map $G^{(j, \gq)}: \cU \times (-\ve, \ve) \subset \cU \times \bR \to \cM$  from the cartesian product of an appropriate real interval $(-\ve, \ve)$ and the open set  $\cU$ (possibly smaller than the open set of the statement)  into $\cM$ 
satisfying the following three conditions: 
\begin{itemize}[leftmargin = 15pt]
\item  for any $y \in \cU$,   the  map 
$G^{(y)(j, \gq)} \= G^{(j, \gq)}(y, \cdot) :(-\ve, \ve)  \to \cM$
is a $\bT$-surrogate map of rank $1$ centred at $y$;  in particular, $G^{(y_o)(j, \gq)}(0) = y_o$; 
\item the $\bT$-surrogate  generators $X_\ell$,  the functions $\s_\ell$ and the tuple $(i_\ell)$ which occur in the definitions of the $\bT$-surrogate maps $G^{(y)(j, \gq)} $ do not depend on the point $y$; 
\item   the derivative $ \frac{d G^{(y_o)(j, \gq)}}{d s}\bigg|_{s = 0}$ 
is equal to 
 \beq \label{ult}  \frac{d G^{(y_o)(j, \gq)}}{d s}\bigg|_{s = 0}  =\bY_j|_{y_o} + v^{(\gq)}_j\eeq
 where $v^{(\gq)}_j$ is a vector in $T_{y_o} \cM$ which  depends real analytically   on  $\gq \in \bR^4$  and  such that,  in  case     
 \beq  \label{7.8} \gq =  (\r_j, \o_j, \r_j, \o_j)\ ,\eeq 
the norm $\|v^{(\gq)}_j\|$  is less than or equal to $\z$. 
 \end{itemize}
For the moment,  take the existence of  maps with these properties as granted. Then, 
 for any $({\bf m} - p)$-tuple of quadruples 
$$\gQ \= (\gq_1, \ldots, \gq_{{\bf m} - p}) \in \bigg ((0, 1) \times (0, T) \times  (0, 1) \times (0, T) \bigg)^{{\bf m} - p} \subset \bR^{4({\bf m} - p)}$$
we may consider 
  the  sequence of $\cD^{II}$-maps $F^{(\ell, \gQ)}: \cV \to \cM$, $1 \leq \ell \leq {\bf m} - p$,  defined  iteratively  on $\cV = (- \ve, \ve)^{{\bf m} }$ by 
\begin{align}
\nonumber &  F^{(1, \gQ)}(s) =  G^{(y_o)(1, \gq_1)}(s^{p +1})\ ,\\
\nonumber  & F^{(2, \gQ)}(s) =  G^{(F^{(1, \gQ)}(s))(2, \gq_2)}(s^{p + 2})\ ,\\
\nonumber & F^{(3, \gQ)}(s) =  G^{(F^{(2, \gQ)}(s))(3, \gq_3)}(s^{p + 3})\ ,\\
\nonumber & \vdots\\
\nonumber & F^{({\bf m }- p, \gQ)}(s) =  G^{(F^{({\bf m}- p-1, \gQ))}(s))({\bf m} - p, \gq_{{\bf m} - p})}(s^{\bf m})
\end{align}
and  set 
\beq 
F^{(\gQ)}: \cV \subset \bR^{\bf m} \to \cM\ ,\qquad  F^{(\gQ)}(s) =  \Phi^{\bW^{(0)}_{A_p}}_{s^p} \circ \ldots \Phi^{\bW^{(0)}_{A_1}}_{s^1} \circ  F^{({\bf m }- p, \gQ)}(s)\ . \eeq
Using the  properties of the maps $G^{(j, \gq)}$ and  the classical properties of the flows of vector fields, one can directly check that for any choice of  $\gQ$, the corresponding map $F^{(\gQ)}$ satisfies  the following conditions: 
\begin{itemize}
\item $F^{(\gQ)}(0, \ldots, 0) = y_o$;  
\item $F^{(\gQ)}$ has the form \eqref{cond1} in which  all vector fields $X_\ell$  are $\bT$-surrogate fields; 
\item  The partial derivatives with respect to the variables $s^i$ at the origin of $\bR^{{\bf m} }$ are
\begin{align} \nonumber &\frac{\p F^{(\gQ)}}{\p s^1}\bigg|_{0_{\bR^{\bf m}}} =  \frac{d \Phi^{\bW^{(0)}_{A_1}}_s}{d s}\bigg|_{0} = \bW^{(0)}_{A_1}|_{y_o}\ ,\qquad  \frac{\p F^{(\gQ)} }{\p s^2}\bigg|_{0_{\bR^{\bf m}}} = \frac{d \Phi^{\bW^{(0)}_{A_2}}_s}{d s}\bigg|_{0} =  \bW^{(0)}_{A_2}|_{y_o}\ ,\qquad\ldots \\
\nonumber &\hskip 6 cm  \ldots \qquad   \frac{\p F^{(\gQ)} }{\p s^p}\bigg|_{0_{\bR^{\bf m}}} =  \frac{d \Phi^{\bW^{(0)}_{A_p}}_s}{d s}\bigg|_{0} = \bW^{(0)}_{A_p}|_{y_o}\ ,\\
\nonumber &  \frac{\p F^{(\gQ)} }{\p s^{p+1}}\bigg|_{0_{\bR^{\bf m}}} = \frac{d G^{(y_o)(1, \gq_1)}}{d s}\bigg|_{0} =   \bY_1|_{y_o} + v^{(\gq_1)}_1\ ,\qquad\ldots 
\ldots  \\
\label{ka} & \hskip 4 cm  \ldots \qquad  \frac{\p F^{(\gQ)}}{\p s^{\bf m}}\bigg|_{0_{\bR^{\bf m}}}  =  \frac{d G^{(y_o)({\bf m} - p, \gq_{{\bf m} - p})}}{d s}\bigg|_{0}  =   \bY_{{\bf m} - p}|_{y_o} + v^{(\gq_{{\bf m} - p})}_{{\bf m}-p}\ .
\end{align}
\end{itemize}
We now recall that  when all  quadruples $\gq_j$, $ 1 \leq j \leq {\bf m} - p$,   are  equal to the quadruples  \eqref{7.8},  the  vectors $v^{(\gq_j)}_j$  have a norm which is less than or equal to  the constant $\z > 0$.  Combining this  with    \eqref{ka} and the fact that   ${\bf m}$-tuple \eqref{genn1*} is pointwise linearly independent,  it follows that,    if $\z> 0$ is taken sufficiently small,  there is  at least one choice of  $\gQ$ (namely the one with the quadruples \eqref{7.8})  which makes  the Jacobian  $JF^{(\gQ)}|_0$   of  maximal rank,  i.e. which   makes at least one  minor   of order ${\bf m} $  of  $JF^{(\gQ)}|_0$  to be non zero.  For fixing the ideas,  we may  assume that such a minor is the determinant of the submatrix  consisting of the first ${\bf m} $ entries for each column. Then the map which  sends each tuple of quadruples  $\gQ \in  \bigg( (0, 1) \times (0, T) \times  (0, 1) \times (0, T) \bigg)^{{\bf m} - p} $ into  the corresponding determinant of the submatrix of    $JF^{(\gQ)}|_0$, made of the first ${\bf m}$ entries in each column, cannot be identically zero.  Since  all coordinate components of  $F^{(\gQ)}$ depend  real analytically  on  $\gQ$, we conclude that there exists an open and dense subset of  tuples $\gQ$ in  $\bigg( (0, 1) \times (0, T) \times  (0, 1) \times (0, T) \bigg)^{{\bf m} - p} \subset  \bR^{4 ({\bf m} - p)}$,  for which the corresponding map $F^{(\gQ)}$ has maximal rank at the origin.   Any  such  map $F^{(\gQ)}$ is therefore   a $\cD^{II}$-map. Furthermore,  using the fact that  the corresponding  set of tuples $\gQ$  is an open and dense subset of   $\bigg( (0, 1) \times (0, T) \times  (0, 1) \times (0, T) \bigg)^{{\bf m} - p} $,  we claim that it is always possible to  select at least one  tuple $\gQ_o$, for which  the  corresponding generators $X_\ell$ occurring in the  expression of $F^{(\gQ)}$  have  $\bT$-depths  satisfying \eqref{constrained}. This  means  that the map $F^{(\gQ_o)}$ is not just a $\cD^{II}$-map but also a $\bT$-surrogate map. The corresponding leaflet $F^{(\gQ_o)}(\cV)$ is therefore  a $\bT$-surrogate leaflet, centred at $y_o$ and  of  dimension ${\bf m} = \rank \cE^{(\cM_{j_o})}$. This   would   conclude the proof, provided 
that  the  following two claims (which we assumed to be true) holds: 
\begin{itemize}[leftmargin = 20 pt]
\item[($\a$)] The maps $G^{(j, \gq)}: \cU \times (-\ve, \ve) \subset \cU \times \bR \to \cM$, with the above described properties, do  exist; 
\item[($\b$)] There exists a $({\bf m} - p)$-tuple   $\gQ_o \in \bigg( (0, 1) \times (0, T) \times  (0, 1) \times (0, T) \bigg)^{{\bf m} - p} $ of quadruples $\gq_j$, whose associated map $F^{(\gQ_o)}$ is a $\bT$-surrogate map.
\end{itemize}
\par
\smallskip 
Let us  prove ($\a$).   From now on, we consider a fixed index  $1 \leq j \leq {\bf m}  -p$ and a fixed quadruple $\gq \= (\wt \r, \wt \o,\wt{ \wt \r}, \wt{\wt \o})$  with $\wt \r, \wt{\wt \r} \in (0,1)$ and $\wt \o, \wt{\wt \o} \in (0, T)$. Thus,   for simplicity of notation,  from now on we  use the  shorter  notation $\bY \= \bY_j$, $\bW \= \bW^{(j)}_{B_j}$,  $\bW'\= \bW^{(j)}_{B'_j}$, $\r = \r_j$ and $\o = \o_j$ and  we denote by $W = W_{\beta_j}$ the   vector field in $\cD^I$,  which  determines the sub-distribution    $(V^{II(W)}, \cD^{II(W)})$ of  stratified uniform type and with    $y_o$ as good point of the first kind, in which  $\bW$ takes values.  We also denote by 
  $\cN$  the  $\cD^{II(W)}$-stratum of $\cM$  that   contains $y_o$. The rank of the (locally regular) distribution $(V^{II(W)}|_{\cN}, \cD^{II(W)}|_{\cN})$  is denoted by ${\bf n}$.
\par
\smallskip
We  now consider  two  sets of  $\bT$-surrogate fields with $\bT$-lengths in  $[\wt \r \wt \o,\wt  \o] $  and $[\wt{\wt \r} \wt{\wt \o}, \wt{\wt \o}] $
\beq \label{718} \left(\wt \bW_A \= \underset{[\wt \r \wt \o , \wt \o]}{\bW_A}\right)_{1 \leq A  \leq {\bf n}}\ ,\qquad \left(\wt{\wt \bW}_C \= \underset{[\wt{\wt \r} \wt{\wt \o} , \wt{\wt \o}]}{\bW_C}\right)_{1 \leq C  \leq {\bf m}}\ ,\eeq
which are    generators for the sub-distribution $(V^{II(W)}|_{\cN}, \cD^{II(W)}|_{\cN})$ and for the distribution  $(V^{II}|_{\cU_{j_o}}, \cD^{II}|_{\cU_{j_o}})$,  respectively (we recall that $\cU_{j_o}$ is the $\cD^{II}$-stratum containing $y_o$  and it is in general different from $\cN$, which is  a $\cD^{II(W)}$-stratum).   We assume that the tuples in \eqref{718} are constructed using the  method of the proof of Theorem \ref{lemma38},   starting by the same   $\bT$-adapted bases which lead to the $\bT$-surrogate generators  that contain  $\bW$ and $\bW'$, respectively. \par
  \smallskip
  Then,   we denote by    $\l^A: \cU \to \bR$, $1 \leq A \leq {\bf n}$  and $\mu^C: \cU \to \bR$, $1 \leq C \leq {\bf m}$,  the real analytic functions  which allow to expand $\bW$ and $\bW'$ as 
\beq 
 \bW|_y   = \l^A(y) \wt \bW_A|_y \ \\ ,\qquad  \bW'|_{y}   = \mu^C(y) \wt{\wt \bW}_C|_y \ , \qquad y \in \cU\ , 
\eeq
and we set  
\beq \l^A_o \= \l^A(y_o)\ ,\qquad \mu^C_o \= \mu^C(y_o)\ .\eeq
Since  the set of $\bT$-surrogate generators  $(\wt \bW_A)$ and  the one which contains $\bW$ are both determined  by the same $\bT$-adapted basis and by the same algorithm (uniquely determined by  $\wt \r $ and $\wt \o$),   the functions  $\l^A$ are determined by the real analytic family of matrices which pointwise transform one set of generators into the other. In particular, the $\l^A$ depend real analytically  on the real numbers $\wt \r$ and  $\wt \o$.  Similarly  the $\mu^C$ depend real analytically  on $\wt{\wt \r}$ and  $\wt{\wt \o}$. We also remark that  when  the quadruple $\gq = (\wt \r,\wt \o, \wt{\wt \r} , \wt{\wt \o}) $  tends to the quadruple  $(\r, \o,\r , \o)$,  the  family of matrices which express  one set of generators into the other tend uniformly to the identity matrices at all points. This implies that  the  functions $\l^A(y)$ and $\mu^C(y)$ tend  uniformly to  the constant functions
$\l^A(y) = \d^A_B$ and $\mu^C(y) =  \d^C_{B'}$ for some appropriate indices $B$, $B'$ for  $\gq \to (\r, \o,\r , \o)$ (more precisely, $B = B_j$ and $B' = B'_j$).\par
\smallskip
Now, for     $t \geq  0$ and $\s \geq 0$ running in a   sufficiently small neighbourhood of $0$ (such smallness  requirement is  imposed to make  all of the subsequent formulas  meaningful), let us consider   the two-parameter family $g_{t,\s}: \cU \to \cM$   of local diffeomorphisms
\beq g_{t,\s}(y) \=  \Phi^{\s \mu^1_o  \wt{\wt \bW}_1}_t \circ \ldots \circ   \Phi_{t}^{\s \mu^{\bf q}_o \wt{\wt \bW}_{\bf q}} \circ \Phi_t^{ \l^1_o \wt \bW_1}\circ \ldots \circ   \Phi_t^{ \l^{\bf n}_o \wt \bW_{\bf n}}  \circ \Phi^{-\s\bW'}_{t}  \circ   \Phi^{- \bW}_{t} (y) \ .\eeq
We want to determine the first and the second  derivatives of this family with respect  to $t$ at $t = 0$ and arbitrary $\s$. This can be done with the help of the next technical lemma, whose proof is postponed to Appendix \ref{proofofsecondlemma}. In its statement, we consider  a fixed  system of coordinates $\xi = (x^1, \ldots, x^n)$ on a neighbourhood of  the point $y_o$ 
and for any  two vector fields $A = A^i \frac{\p}{\p x^i}, B = B^j\frac{\p}{\p x^j}$ on such a neighbourhood,  we use the notation  $A(B)$ to denote the vector field 
$$A(B) \= A(B^j) \frac{\p}{\p x^j} \ .$$
\par
\begin{lem}  \label{secondlemma} Given $m+2$  local vector fields   $X_1, \ldots, X_m$,  $Y, Z$    defined on a coordinatizable  neighbourhood of a point $y \in \cM$,  their local flows   satisfy  the following relations:
 \begin{itemize}[leftmargin = 20pt]
 \item[(i)]
\begin{multline} \label{9.13}
  \frac{d}{d s} \left(\Phi^{X_1}_s \circ \ldots \circ \Phi^{X_m}_s(y) \right)\bigg|_{s = 0} = \sum_{i = 1}^m X_i|_{y}\qquad\text{and}\\
\frac{d^2}{d s^2} \left(\Phi^{X_1}_s \circ \ldots \circ \Phi^{X_m}_s(y) \right)\bigg|_{s = 0} 
= \sum_{i = 1}^m \left((X_i + X_{i + 1} + \ldots + X_m)( X_i\right)|_y  + \\
+ \left(X_{i + 1} + \ldots + X_m)( X_i)|_y\right)
= 
 \sum_{i = 1}^m X_i(X_i)|_{y} + 2 \sum_{j = 2}^m \sum_{i = 1}^{j-1} X_j(X_i)|_{y} \ ;
\end{multline}
\item[(ii)]
\begin{multline}
\frac{d}{d s} \left(\Phi^{Z + Y}_s \circ \Phi^{-Z}_s \circ \Phi^{- Y}_s(y) \right)\bigg|_{s = 0} =0\qquad \text{and}\\
 \frac{d^2}{d s^2} \left(\Phi^{Z + Y}_s \circ  \Phi^{-Z}_s \circ \Phi^{- Y}_s(y) \right)\bigg|_{s = 0} = [Y, Z]|_{y}\ ;
 \end{multline}
 \item[(iii)]  Setting $\cA \=  \sum_{i = 1}^m X_i  -  Z  - Y$, then $\frac{d}{d s} \left(\Phi^{X_1}_s \circ \ldots \circ \Phi^{X_m}_s \circ \Phi^{-Z}_s \circ \Phi^{-Y}_s(y)\right)\bigg|_{s = 0} = \cA_y$  
  and 
 \begin{multline}
 \frac{d^2}{d s^2} \left(\Phi^{X_1}_s \circ \ldots \circ \Phi^{X_m}_s \circ \Phi^{-Z}_s \circ \Phi^{-Y}_s(y)\right)\bigg|_{s = 0}  = \\
 \hskip 1 cm = [Y, Z]_y  +  [\cA, Y+ Z]_y + \cA(\cA)|_y -  \sum_{\ell = 1}^m  \sum_{j =  \ell+1}^m X_\ell(X_j) |_y + \sum_{\ell = 2}^m \sum_{j = 1}^{\ell-1} X_\ell(X_j)|_{y}\ .\end{multline}
 \end{itemize}
\end{lem}
With the help of this,  we can prove the following
\begin{lem} \label{firstlemma} On any  coordinatizable  neighbourhood of a  point $y \in  \cN \cap \cU_{j_o}$, 
the family of local diffeomorphisms $g_{t, \s}$   has the form 
\begin{multline} g_{t, \s}(y) = y + t \left(\cA_{\gq}|_y+ \s \cB_{\gq}|_y\right) 
+  t^2  \bigg(\s [\bW, \bW']|_y + \\
+  [\cA_{\gq}+ \s \cB_{\gq}, \bW+ \s \bW']_y 
+ (\cA_{\gq}+ \s \cB_{\gq})(\cA_{\gq} + \s \cB_{\gq})_y+ \\
+ \s^2 Y_{\gq} (y) + \s Y'_{\gq} (y) +  Y''_{\gq} (y)  \bigg) + \gr_{t, \s}(y)\ ,\end{multline}
where 
\begin{itemize}[leftmargin = 20pt]
\item[(a)] $\cA_{\gq}$,  $\cB_{\gq}$, $Y_{\gq}$, $Y'_{\gq}$ and $Y''_{\gq}$ are vector fields,  which are independent of $t$ and $\s$ and such that,  when $\gq$ tends to  $  (\r, \o, \r, \o)$, the vector fields $\cA_{\gq}$,  $\cB_{\gq}$, $Y_{\gq}$,  $Y''_{\gq} $ tend together with their first derivatives uniformly  to   $0$,   while  $Y'_{\gq} $  tends together with its first derivatives uniformly  to   $[\bW, \bW']$; 
\item[(b)] The vector fields  $\cA_{\gq}$ and  $\cB_{\gq}$  are such  that $\cA_{\gq}|_{y_o} = \cB_{\gq}|_{y_o} = 0$; 
\item[(c)]  for any fixed   $\s$  and $y$ we have that $\gr_{t,\s}(y) = o(t^2)$. 
\end{itemize}
\end{lem}
\begin{pf}  Let   $\cA = \cA_{\gq} \= \l^A_o \wt \bW_A-  \bW  $ and $\cB = \cB_{\gq}\=  \mu^C_o\wt{\wt \bW}_C - \bW' $. 
Since $g_{t ,\s}$ is a composition of flows parameterised by $t$,   we have $g_{t = 0,\s}(y) = y$ for any $y$ and $\s$. Moreover, from  Lemma \ref{secondlemma} (iii) 
 \begin{align}
\label{due}&\frac{d g_{t,\s}(y)}{dt}\bigg|_{(t = 0, \s, y)} = \cA|_y+ \s \cB|_y\\
\nonumber &\frac{d^2 g_{t,\s}(y)}{dt^2}\bigg|_{(t = 0, \s, y)}  =\left(\s [\bW, \bW'] +  [\cA+ \s \cB, \bW+ \s \bW']\big|_y + (\cA+ \s \cB)(\cA + \s \cB)\big|_y+ \right.\\
\label{tre} & \hskip 2 cm + \left.\s^2 Y_{\gq} (y) + \s Y'_{\gq} (y) +  Y''_{\gq} (y)  \right)
\end{align}
for some appropriate vector fields $Y_{\gq}, Y'_{\gq}$ and  $Y''_{\gq} $,  given by  appropriate summations between  vector fields  of the form $\l^A_o\mu_o^C \wt W_A( \wt{\wt W}_C)$ and $\l^A_o \mu^C_o \wt{\wt W}_B( \wt W_A)$. Due to  the definition of the constants $\l^A_o$ and $\mu^C_o$, the vector fields $\cA$ and $\cB$  vanish at $y_o$ and   tend uniformly to $0$ together with all their first derivatives when  $\gq= (\wt \r, \wt \o, \wt{\wt \r} , \wt{\wt \o})$ tends to  $(\r, \o, \r, \o)$ because of their definitions and the  fact that the coefficients $\l_o^A$ and $\mu_o^C$ tend to  appropriate coefficients  $1$ or $0$,  when the quadruple $\gq$ becomes  $(\r, \o, \r, \o)$.   When $\gq$ tends to $(\r, \o, \r, \o)$, the vector fields  $Y_{\gq}$ and  $Y''_{\gq} $
tend uniformly to $0$,  while   $Y'_{\gq}$ tends to $[\bW, \bW']$ because for any choice of (a sufficiently small) $\s$ 
the second derivative   $\frac{d^2 g_{t,\s}(y)}{dt^2}\bigg|_{(t = 0, \s, y)}  $ tends to 
$$ \frac{d^2}{dt^2}\bigg|_{(t = 0, \s, y)} \Phi^{\s \bW'}_t \circ  \Phi_t^{ \bW}\circ  \Phi^{-\s\bW'}_{t} \circ   \Phi^{- \bW}_{t}(y) = 2 \s [\bW, \bW'] _y\ .$$
From   these remarks, the claim follows directly.
\end{pf}
\noindent
Now, for any sufficiently small $\d> 0$,  let us denote by  $y^{(\d)}_o $ the point 
$$y^{(\d)}_o \=  \Phi_\d^{ \l^1_o\wt \bW_1}\circ \ldots \circ   \Phi_\d^{ \l^{\bf n}_o\wt \bW_{\bf n}}  \circ \Phi^{-\bW}_{\d}(y_o) \ .$$ 
 We stress the fact   that  $y^{(\d)}_o$ is  in the integral leaf  $\cS^{(W|y_o)}$ passing through $y_o$. 
We also recall that, by  assumptions, $y_o$ is a good point of the first kind for the sub-distribution which contains $W$. By  the proof of Theorem \ref{prop93}, this implies that  there is  a  neighbourhood  of $y_o$  of the intrinsic topology of the  immersed submanifold $\cS^{(W|y_o)}$  which is a $\bT$-surrogate leaflet centred at $y_o$, i.e. the image of a  map of the form
\beq \label{correction} (s^1, \ldots, s^{\bf n}) \in \cV' \subset \bR^{\bf n} \longmapsto \Phi^{-\wh \bW_1}_{s^1} \circ \ldots \circ \Phi^{-\wh \bW_{\bf n}}_{s^{\bf n}}(y_o) \in \cS^{(W|y_o)}\eeq
determined by an appropriate  set of $\bT$-surrogate generators $(- \wh \bW_i)$. The proof holds also if the $\bT$-depths of the vector fields $\wh W_A$ satisfy the inequalities 
$  \t_1 >  \t_2 >  \ldots > \t_{\bf m} $  (instead of $\t_{\bf m} > \ldots > \t_2 > \t_1$, as it is required for the $\bT$-surrogate maps). If we assume that such a sequence of inequalities is satisfied  and that $\d$ is so small that $y_o^{(\d)}$ is in the image of the map \eqref{correction}, by reversing the order of the flows and changing signs to the  vectors, we get that $y_o$ is in the image of the map  
\beq \label{correctionbis} (s^1, \ldots, s^{\bf n}) \in \cV' \subset \bR^{\bf n} \longmapsto \Phi^{\wh \bW_{\bf n
 }}_{s^{\bf n}} \circ \ldots \circ \Phi^{\wh \bW_1}_{s^1}(y_o^{(\d)}) \ ,\eeq 
which is now  a true $\bT$-surrogate map (because now the inequalities satisfied by the $\bT$-depths  of the vector fields $\wh \bW_A$ are the correct ones).\par
Summing up, we conclude that  for any sufficiently small $\d$,  there is a unique 
 ${\bf n}$-tuple   $(\nu^1, \ldots, \nu^{\bf n})$ of real numbers such that 
 the local diffeomorphism 
$\wt G^{(\d, y_o)} =  \Phi_{\nu^{\bf n}}^{ \wh \bW_{\bf n}}\circ \ldots \circ   \Phi_{\nu^1}^{  \wh \bW_1} $   is such that $\wt G^{(\d, y_o)}(y^{(\d)}_o) = y_o$.
 Note    that: 
 \begin{itemize}[leftmargin = 20pt]
 \item for $\d \to 0$, the point $y^{(\d)}_o$ tends to $y_o$ and  the  tuple $(\nu^1, \ldots, \nu^{\bf n})$ tends to $0_{\bR^{\bf n}}$; 
 \item  the Jacobian at $y^{(\d)}_o$ of $\wt G^{(\d, y_o)}$  tends to the identity matrix for $\d \to 0$; 
 \item  the $\bT$-surrogate generators  $(\wh \bW_A)$ that give the leaflet can be chosen  with  $\bT$-depths in any prescribed  interval $[\wh \r \wh \o, \wh \o]$. 
 \end{itemize}
 \smallskip
 Consider now the map $G^{(j, \gq)}: \cU \times (-\ve, \ve) \subset \cU \times \cM \to \cM$ defined by 
\beq G^{(j, \gq)}(y,s)  \= \wt G^{(\d, y_o)} \circ g_{\d, \frac{s}{2\d^2}}(y)\eeq
for some fixed choice of   $\d > 0$. We claim that, for an appropriate choice of $\d$,  this map satisfies the conditions of claim ($\a$). Indeed, by the above remarks, if $\d$ is sufficiently small, 
denoting   $\cA = \cA_{\gq} $ and $\cB = \cB_{\gq}$ and recalling that $\cA|_{y_o} = \cB|_{y_o} = 0$ we have that 
\begin{itemize}[leftmargin = 20pt]
\item[(1)] $G^{(j, \gq)}(y_o,s = 0) = \wt G^{(\d, y_o)}(y^{(\d)}_o) = y_o$; 
\item[(2)]
\begin{multline}  \frac{d}{ds} G^{(j, \gq)}(y_o,s )\bigg|_{s = 0} = \frac{1}{ 2 \d^2} J  \wt G^{(\d, y_o)}|_{y_o} \cdot \frac{d}{d\s} g_{\d, \s}\bigg|_{\s = 0, y = y_o} = \\
=   \frac{1}{2 \d^2} J  \wt G^{(\d, y_o)}|_{y_o} \cdot   \Bigg( \xcancel{ \d   \cB_{y_o}}  + \d^2 [\bW, \bW']_{y_o} +  \d^2[\cB, \bW]_{y_o}   +  \d^2[\cA, \bW']_{y_o}  +\xcancel{ \d^2\cB(\cA)_{y_o}} +  \\
+ \xcancel{\d^2\cA(\cB)_{y_o}} + + \d^2 Y'_{\gq} (y_o) +   \frac{d}{d\s}  \gr_{\d, \s}(y)\bigg|_{\s = 0,  y =  y_o}\Bigg) =\\
= \frac{1}{2} J  \wt G^{(\d, y_o)}|_{y_o} \cdot    \Bigg( [\bW, \bW']_{y_o} -  \bW(\cB)|_{y_o} - \bW'(\cA)|_{y_o} + Y'_{\gq} (y_o) +  \frac{1}{\d^2} \frac{d}{d\s}  \gr_{\d, \s}(y)\bigg|_{\s = 0, y =  y_o}\Bigg) 
  \end{multline}
\end{itemize}
If we now denote   
\beq 
 \wt v \=  \frac{1}{2}\Bigg(  [\bW, \bW']_{y_o} -  \bW(\cB)|_{y_o} - \bW'(\cA)|_{y_o} + Y'_{\gq} (y_o)
 +    \frac{1}{\d^2} \frac{d}{d\s}  \gr_{\d, \s}(y)\bigg|_{\s = 0, y =  y_o} \Bigg) \ ,
\eeq
  and if  we  express the remainder $\gr_{\d, \s}(y)$ in  integral form, one can directly check  that this vector $\wt v$ of $T_{y_o} \cM$ depends in a real analytic way from the quadruple $\gq = (\wt \r, \wt \o, \wt{\wt \r},  \wt{\wt \o})$ and,   for $\gq \to ( \r, \o, \r, \o)$, it tends to
\beq \bY \=  [\bW, \bW']_{y_o} + \wt{\wt v}^{(\d)}\eeq
 where $ \wt{\wt v}^{(\d)}$ is a vector, whose components are infinitesimal  of the same order of $\d$ for $\d \to 0$.  Since  the matrix $ J  \wt G^{(\d, y_o)}|_{y_o}$  tends to the identity map for $\gq$ tending to $( \r, \o, \r, \o)$,  if we choose a sufficiently small $\d$, we conclude that  $\frac{d}{ds} G^{(j, \gq)}(y_o,s)\bigg|_{s = 0, y_o}$ has the form \eqref{ult}, as we needed to prove. \par
\smallskip
It remains to prove ($\b$), i.e. that one can determine each quadruple $\gq_i$ of the full ordered set of quadruples $\gQ = (\gq_j)$ and construct the associated maps $G^{(j, \gq)}$ in such a way that  $F^{(\gQ_o)}$ is a $\bT$-surrogate map. This can be done by choosing the quadruple $\gq_j$ inductively as follows. First, select a quadruple $\gq_1$  so that the intervals $[\wt \r_1 \wt \o_1, \wt \o_1]$ and  $[\wt{\wt \r}_1 \wt {\wt \o}_1, \wt{\wt \o}_1]$ are with $\wt {\wt \o}_1 < {\wt \r}_1  {\wt \o}_1$ and included in the open interval $(\o_2, \r_1 \o_1)$ (in this way  the interval that  contains the $\bT$-depths of the vector fields $\bW_{B_1}^{(1)}$ and $\bW_{B'_1}^{(1)}$ 
are disjoints with  the intervals $[\wt \r_1 \wt \o_1, \wt \o_1]$ and  $[\wt{\wt \r}_1 \wt {\wt \o}_1, \wt{\wt \o}_1]$ and 
all $\bT$-depths of the vector fields $\wt{\wt \bW}_A$ are strictly smaller than those of the vector fields  $\wt \bW_A$ and the latter  strictly smaller than those of $\bW_{B_1}^{(1)}$ and $\bW_{B'_1}^{(1)}$). Moreover, with no loss of generality, we may assume that the vector fields $\bW_{A}^{(1)}$ are ordered so that the $\bT$-depth of $\bW_{B'_1}^{(1)}$  is strictly smaller than the one of $\bW_{B'_1}^{(1)}$.   Finally, we may  choose the generators $\wh \bW_A$, which appear in the construction of the map  $G^{(1, \gq_1)}$ with $\bT$-depths in the open  interval  $(\o_2, \wt{\wt \r}_1 \wt{\wt \o}_1 )$. All these conditions  on the intervals containing the various $\bT$-depths  are imposed in order to  guarantee that  the  $\bT$-depths  of the vector fields used in the construction of the map $G^{(1, \gq_1)}$
satisfy the necessary inequalities to  make   $G^{(1, \gq_1)}$ a $\bT$-surrogate  map. After these conditions are imposed,  we can make  similar choices for what concerns  quadruple $\gq_2$, i.e. such that the intervals $[\wt \r_2 \wt \o_2, \wt \o_2]$ and  $[\wt{\wt \r}_2 \wt {\wt \o}_2, \wt{\wt \o}_2]$ are with $\wt{\wt \o}_2 < \wt \r_2 \wt \o_2$  and included in the open interval $(\o_3, \r_2 \o_2)$ and appropriate choices for the interval containing the  $\bT$-depths of the remaining $\bT$-surrogate vector fields which appear in the definition of   $G^{(2, \gq_2)}$, so that also the latter is a $\bT$-surrogate map generated by vector fields with $\bT$-depths larger than $\o_3$. And so on.  After a finite number of steps, one fixes all quadruples and in the definition of the map $F^{(\cQ)}$, all the $\bT$-depths of the vector fields which generate the map are ordered as desired. Since the tuples  of quadruples $\gQ \in  \bigg( (0, 1) \times (0, T) \times  (0, 1) \times (0, T) \bigg)^{{\bf m} - p} $, whose associated  maps $F^{(\gQ)}$ have  Jacobian of  maximal rank at the origin, belong to  open and dense subset of   $ \bigg( (0, 1) \times (0, T) \times  (0, 1) \times (0, T) \bigg)^{{\bf m} - p}$\!\!\!\!, the above described  finite  collection of  quadruples $\gq_j$ and of set of  generators $(\wh W_A)$, used  in construction of the maps $G^{(j, \gq)}$,  can be done in such a way that  $F^{(\gQ)}$  satisfies the conditions of being a $\bT$-surrogate map. \par
\bigskip
\section{Applications and possible developments}\label{section8}
As we pointed out in \S \ref{surr-leaflets},  the hyper-accessibility  and the small-time local controllability of a non-linear control  real analytic system can be established by means of  some relatively simple computations, provided that all points of the extended space-time $\cM = \bR \times \cQ \times \cK$ are good points. The first and the second criteria for goodness given  in \S \ref{theeigthsection} allow to easily  identify the systems with  such a property.  In   \cite{GSZ2} we give several examples of how this circle of  ideas apply. More precisely:  
\begin{itemize}[leftmargin = 20pt]
\item[(1)] We show how the classical Kalman criterion for the {\it local} controllability of linear control systems can be obtained as an immediate  corollary of our first criterion of goodness and of the fact that the integral leaves of the secondary distribution of a linear system  with an $m$-dimensional control space are  parallel affine subspaces, whose  dimension is  determined by   the rank of the Kalman matrix; 
\item[(2)] We discuss  a few elementary non-linear control systems  and, for each of them,  check  whether  our criterions for goodness do or do not apply. Since for these  very simple examples  the controllability or the non-controllability can be  also established  by other elementary methods,  this discussion is helpful to have a deeper insight  on how  our new  approach  is consistent with (and actually improves) the  classical methods;  
\item[(3)] We give an explicit  example of  a non-linear real analytic control system, for which the classical Kalman linear test is inconclusive, while our second criterion for goodness  and our theory of surrogate $\bT$-leaflets    is able to establish the small-time controllability at the points of stability of that system;  
\item[(4)] We use the results of this paper to prove that the   systems  of the   {\it  controlled Chaplygin sleigh}  and  of its {\it hydrodynamical} variant (\cite{FN, SZ}) have the hyper-accessibility property and hence the small-time local controllability property at each stable  point.  At the best of our knowledge, this is  the first place where  this  property is proved  for such    classical control system. \
\end{itemize} \par
\smallskip 
For what concerns possible future  developments, we  would like to stress   that the  idea of considering  the surrogate $\bT$-leaflets  and their projections onto the state space $\cQ$ can be used to discuss  not only the local controllability but also  the   global controllability. For instance, the analysis in \cite{GSZ2} of the  linear control systems  indicates that   the surrogate $\bT$-leaflets  can be arbitrarily enlarged and fill the orbits  of an appropriate  abelian Lie group of translations.  When the dimensions of such orbits is sufficiently large, the projections onto  $\cQ$  of the leaflets  coincide with the whole space and all points of $\cQ$ are therefore  reachable. This leads to a new proof of  the Kalman criterion for {\it global} controllability.   It is reasonable to  expect that a similar argument can be used to determine  the global controllability of many other  types of  non-linear  control systems, namely  for those  whose   secondary distributions admit generators that constitute   Lie algebras of nilpotent Lie groups with surjective exponential maps. \par
\smallskip
Other tools  for establishing  global controllability  properties might  come   from estimates of the sizes of  the surrogate $\bT$-leaflets and of their projections on $\cQ$.   Uniform lower  bounds for  the radii of the  balls contained in such  $\cQ$-projections can be exploited  in  ``open and closed"  arguments   and hence used to show that, in certain settings,  the  reachable sets   coincide with the whole state space. Estimates  of this kind  for arbitrary compositions of flows have been determined in \cite{GSZ1}.  We expect that  similar estimates  can be  determined   for  $\bT$-surrogate maps and $\bT$-leaflets as well. \par
\smallskip 
We would like to conclude this paper with a brief comment on the real analyticity assumptions that are here considered. Since they have been heavily exploited in  many arguments, it is hard to believe that they might be  removed in a straightforward way. Nonetheless, we conjecture   that the new notions  of rigged distributions, secondary distributions,  $\bT$-leaflets etc., might  be quite  useful   for studying   also  several kinds of non-linear control systems  of class $\cC^\infty$ or lower. Classical approximation techniques might be  appropriate tools  for extending  some of the results  of this paper to  control problems  of such lower regularity (see f.i.  \cite{CGS}).\par
\vskip  2 cm 
\appendix

\addcontentsline{toc}{section}{APPENDIX}
\centerline{\large \it APPENDIX}
\ \\[-45pt]

 \section{The proof of Lemma \ref{lemmone}} 
 \label{appendixb}
 Consider a set of generators $\{X_\b\}$ for the regular distribution $ \cD^{I(0)}|_\cU\= \cD^{I}|_\cU$ for a sufficiently small open subset $\cU \subset \cV$ and, for   any $ 0 \leq \ell \leq \nu$ and  $z \in \cU$,  let  us denote by $\cD^{I(\ell)}|_{z} $   the subspace of $ T_{z} \cM$,   which is spanned by the values at $z$ of the local vector fields 
 \beq 
Y =  \sum_{k = 0}^\ell r^{\b}_{(k)} X_\b^{(k)}\quad \text{with}\  r^{\b}_{(k)}: \cU \to \bR\ \text{real analytic}\  \text{and}\ \  X_\b^{(k)} \= \underset{\text{$k$-times}}{\underbrace{[\bT, [\bT, \ldots, [ \bT, }}X_\b]\ldots]]\ .
\eeq
In general,  for a given $\ell$,  the corresponding vector spaces  $\cD^{I(\ell)}|_{z} $, $z \in \cU$, do not  have all the same dimension. However, by the semicontinuity property of the rank function and real analyticity, there is an open and dense subset $\cU' \subset \cU$,  on which  the vector spaces  $\cD^{I(\ell)}|_{z} $, $z \in \cU'$, $0 \leq \ell \leq \nu$,  have constant maximal dimension for any $\ell$. Let $z_o$ be a fixed point of $\cU'$ and   denote by $\bV = \bV^0 \oplus \ldots \oplus \bV^\nu$  the direct  sum of the 
 vector spaces 
 \begin{multline*} \bV^0 = \cD^{I(0)}|_{z_o}\ ,\ \ \bV^1 = \cD^{I(1)}|_{z_o}/ (\cD^{I(0)}|_{z_o})\ , \qquad\ldots \\
  \ \ldots \qquad \bV^\nu \= \cD^{I(\nu)}|_{z_o}/(\cD^{I(\nu -1)}|_{z_o})\ .\end{multline*}
As a  vector space,  $\bV = \bV^0 \oplus \ldots \oplus \bV^\nu$ is  isomorphic to  $  \cD^{II}|_{z_o} \simeq \bR^{{\bf m}}$ with  ${\bf m} \= \max_{z \in \cU}  \dim  \cD^{II}|_z$. 
For  each  subspace $\bV^{\ell} \subset \bV$, its elements are called  {\it homogeneous  elements of degree  $\ell$}.\par
 The natural action of the Lie derivative operator  $\cL_{\bT} = [\bT, \cdot]$  on the (local) vector fields in $(V^{II}, \cD^{II})$  induces a natural  grade  shifting   map on the homogeneous vectors of $\bV$, which we now describe.  
Given $v^{(\ell)} \in \cD^{I(\ell)}|_{z_o} $,  we  denote by $[v^{(\ell)}]$  the corresponding equivalence class  in   $\bV^\ell \= \cD^{I(\ell)}|_{z_o}/\cD^{I(\ell-1)}|_{z_o}$. 
For each  class $[v^{(\ell)}]$,   let us also consider  
a real analytic   vector field  $X^{[v^{(\ell)}]}$   with values in the spaces  $\cD^{I(\ell)}_z$   and such that 
\beq \label{conditions} \left[X^{[v^{(\ell)}]}|_{z_o}\right] = [v^{(\ell)}] \ .\eeq
 Vector fields  $X^{[v^{(\ell)}]}$  satisfying \eqref{conditions} can be  constructed  as follows.  By   the maximality of  $\dim  \cD^{I(j)}|_{z_o} $ and  considering a finite set of  real analytic generators  $\{Y_m\}$  for $\cD^{II}|_{\cU} = \cD^{I(\nu)}|_{\cU}$ of   the form $[\bT, \ldots, [\bT, [\bT, X_j]]\ldots]$,  up to a reordering, we may assume that  such a  set of generators splits  into a disjoint union of subsets of the form  
 \begin{multline*}\{Y_m\} = \big\{Y^{(0)}_{A_0}, \ 1 \leq A_0 \leq \dim  \cD^{I(0)}|_{z_o}\big\}
  \cup \\
  \cup \big\{Y^{(1)}_{A_1}\ ,\  \dim  \cD^{I(0)}|_{z_o} + 1 \leq A_1 \leq \dim  \cD^{I(1)}|_{z_o}\big\}\cup \ldots\\
  \ldots \cup  \big\{Y^{(\nu)}_{A_\nu}\ ,\  \dim  \cD^{I(\nu-1)}|_{z_o} + 1 \leq A_\nu \leq \dim  \cD^{I(\nu)}|_{z_o}\big\}
 \end{multline*} 
 such that  {\it for any given  $0 \leq \ell \leq \nu$,  all vector fields $Y^{(\ell')}_{A_{\ell'}}$ with  $0 \leq \ell' \leq \ell$   are in the spaces $\cD^{I(\ell)}|_z$,  $z \in \cU$
  and the  tuple  
 \beq \label{firstb} \left(Y^{(\ell')}_{A_{\ell'}}\big|_{z_o}\right)_{\smallmatrix 1 \leq \ell' \leq \ell\ ,\hfill\\
 \dim  \cD^{I(\ell'-1)}|_{z_o} +1 \leq A_{\ell'}  \leq  \dim  \cD^{I(\ell')}|_{z_o} \endsmallmatrix}\eeq 
 is a basis for   $\cD^{I(\ell)}|_{z_o} $}. 
This implies that  any vector $v^\ell \in \cD^{I(\ell)}|_{z_o}$  admits a unique expansion  
 $$v^\ell = \sum_{A, k} v^\ell{}^{A_k}_{k} Y^{(k)}_{A_k}|_{z_o}\ ,\qquad v^\ell{}^{A_k}_{k} \in \bR\ .$$
 Therefore,  any real analytic vector field  of the form 
\beq \label{form} X^{[v^{(\ell)}]} \= \sum_{k = 0}^\ell \sum_{A_k} v^\ell{}^{A_k}_{k} Y^{(k)}_{A_k} + Z\ ,\eeq 
 for some vector field   $Z$   such that $Z\big|_{z_o}$ is in the space 
\beq \label{form1}  \Span_\bR\left \langle  Y^{(k)}_{B_k}\big|_{z_o}\ ,\ \  0\leq k \leq \ell-1\ ,\ 
 1 \leq B_k \leq \dim  \cD^{I(\ell-1)}|_{z_o} \right\rangle\ , \eeq
 satisfies \eqref{conditions}. From now on, {\it we assume that $X^{[v^{(\ell)}]}$ has the  form \eqref{form}}.\par 
Now, given a vector field $X^{[v^{(\ell)}]}$, $0 \leq \ell \leq \nu$,  we may consider  the class $[w^{(\ell+1)}]$ in  $\bV^{\ell+1}$ determined by the vector   
\beq \label{definiz}  w^{(\ell + 1)} \= [\bT, X^{[v^{(\ell)}]}]|_{z_o}\ . \eeq
We claim that the equivalence class  $[w^{(\ell+1)}]$ is well defined. Indeed,  if $X^{[v^{(\ell)}]}$  is replaced by any other vector field $ X'{}^{[v^{(\ell)}]}$ of the form  \eqref{form}  (i.e.  with  $Z =  X'{}^{[v^{(\ell)}]} -  X^{[v^{(\ell)}]}$  such that $Z|_{z_o} $ is in the space \eqref{form1}), 
then   $w^{(\ell + 1)}$  changes into  
$$
w'{}^{(\ell +1)} =  [\bT, X^{v^{(\ell)}} + Z]|_{z_o} 
= w^{(\ell +1)}  +  \underset{\in \cD^{I(\ell)}|_{z_o}}{\underbrace{[\bT,  Z]|_{z_o}} } \simeq
 w^{(\ell +1)} \hskip -7 pt{\mod \cD^{I(\ell)}|_{z_o} }\ , 
$$
showing that   the equivalence class   $[w^{(\ell +1)}]$ does not change and   is   uniquely associated with  $[v^{(\ell)}]$. We denote such  equivalence class  by 
$$[w^{(\ell +1)}] = \ad_{\bT}([v^{(\ell)}]) \ .$$  Notice that, in case   $\ell = \nu$, the equivalence class   $\ad_{\bT}([v^{(\nu)}]) $ is trivial for any choice  of  $[v^{(\nu)}] \in \bV^{\nu}$. 
The  unique linear endomorphism of $\bV$,  which is defined on  homogeneous vectors as above, is called  
the  {\it adjoint  $\bT$-action on $\bV$}. We denote it by $\ad_\bT: \bV \to \bV$. \par
By construction, for any $1 \leq s \leq \nu$,  the iterated linear map $(\ad_\bT)^s: \bV \to \bV$ shifts the degrees of  homogeneous elements by $s$ units and, consequently, 
 $$(\ad_\bT)^{\nu+1}(\a) = 0\qquad \text{for any}\ \a \in \bV\ .$$ 
We now construct a distinguished basis for $\bV$ as follows. Let $R_{\nu} \= \dim \bV^0/\ker (\ad_\bT)^{\nu}$
 and select an ordered $R_{\nu}$-tuple of  vectors  $([w_{0(\nu)1}], \ldots, [w_{0(\nu)R_{\nu}}])$ in $\bV^0$, which project onto a basis for the quotient $\bV^0/\ker (\ad_\bT)^{\nu}$. Then, for any $0 \leq \ell \leq \nu$,  let us denote by $([w_{\ell(\nu)1}], \ldots, [w_{\ell(\nu)R_{\nu}}])$  the  $R_{\nu}$ vectors in $\bV^\ell$  
 defined by 
 \beq [w_{\ell(\nu)j}] \= (\ad_\bT)^{\ell}\left([w_{0(\nu)j}]\right)\ .\eeq
Since  each map, which is   induced on $\bV^0/\ker (\ad_\bT)^{\nu}$  by  an iterated adjoint $\bT$-action $(\ad_\bT)^{\ell}$, $1 \leq \ell\leq \nu$,     has   trivial kernel,  the  vectors $[w_{\ell(\nu)j}]$, $ 1 \leq j \leq R_\nu$, $0 \leq \ell\leq \nu$,   are     linearly independent. 
Let us now  define
$R_{\nu -1} = \dim  \left( \ker (\ad_\bT)^{\nu}|_{\bV_0}\Big/ \ker (\ad_\bT)^{\nu-1}|_{\bV_0} \right)$
  and select an ordered $R_{\nu -1}$-tuple of   vectors  $([w_{0(\nu-1)1}], \ldots, [w_{0(\nu-1)R_{\nu-1}}])$ in $\ker (\ad_\bT)^{\nu}|_{\bV_0}$, which project onto a basis of the quotient $  \ker (\ad_\bT)^{\nu}|_{\bV_0}/ \ker (\ad_\bT)^{\nu-1}|_{\bV_0}$. Then, for any $1 \leq \ell \leq \nu -1$, 
we may consider  the linearly independent  $R_{\nu -1}$-tuple  $([w_{\ell(\nu -1)1}], \ldots, [w_{\ell(\nu-1)R_{\nu -1}}])$   in $\bV^\ell$  
 defined by
 \beq [w_{\ell(\nu-1)j}] \= (\ad_\bT)^{\ell}\left([w_{0(\nu-1)j}]\right)\ .\eeq
In a similar way,  for any  other integer $1 \leq s \leq \nu-2 $ we may select a  $R_{s}$-tuple  $([w_{0(s)1}], \ldots, [w_{0(s)R_{s}}])$  in $\bV^1$, which projects  onto a basis for
$$ \ker (\ad_\bT)^{s+1}|_{\bV_1}\Big/ \ker (\ad_\bT)^{s}|_{\bV_1}$$  and determine   the associated  linearly independent $R_{s}$-tuples in  the  spaces $\bV^{\ell}$, $1 \leq \ell \leq s$, obtained as images under  the  linear maps $(\ad_\bT)^{\ell}$. 
Finally, we may determine a  basis $([w_{0(0)1}], \ldots, [w_{0(0)R_{0}}])$ for $\ker \ad_\bT$ (to which we  do not associate any element in the subspaces of degree higher than $0$). 
\par \smallskip
One can  check that   the full  collection of homogeneous elements  selected as above
$$  \left\{[w_{\ell(a)j}]\right\}_{\smallmatrix 0 \leq a \leq \nu\\
 0 \leq\ell \leq  a ,\ 1 \leq j \leq R_{a}
\endsmallmatrix}$$
 is a  basis  of   $\bV \simeq   \cD^{II}|_{z_o} $ and that any associated  set of vectors  $ \left\{w_{0(a)j}\right\}_{0  \leq  a \leq \nu, 1 \leq j \leq R_{a}}$
 is a basis of  $\cD^{I}|_{z_o}$.  Now, let us denote by $W_{0(a)j}$, $0 \leq  a \leq \nu, 1 \leq j\leq R_{a}$, some real analytic vector fields  in $\cD^{I}|_{\cU}$ that have the form \eqref{form} with $Z = 0$ and   satisfy the condition
\beq \label{lalla} W_{0(a)j}|_{z_o} =  w_{0(a)j}\ .\eeq
Since we are  assuming  that the generators $\{Y_m\} = \{Y^{(\ell)}_{A_\ell}\}$ are all vector fields having  the form $[\bT, \ldots, [\bT, [\bT, X_j]]\ldots]$ for appropriate generators  $X_j$  for  $\cD^{I}|_{\cU}$,  the vector fields  $W_{0(a)j}$ are   linear combinations with constant coefficients of the generators  $X_j$.
Finally, for any $1 \leq \ell \leq a$, let $W_{\ell(a)j}$ be the  vector fields defined by 
\beq \label{lallina}   W_{\ell(a)j} \= \underset{\ell\text{\rm -times}}{\underbrace{[\bT, [\bT, [\bT,  \ldots [\bT, }}W_{0(a)j}] \ldots ]]]\ .\eeq
 By construction and  previous remarks, each vector field $W_{\ell(a)j}$ takes values in the spaces $ \cD^{I(\ell)}|_z$, $z \in \cU$,  it is a linear 
combination with constant coefficients of an appropriate  set of  generators (possibly not linearly independent at all points of $\cU$) $X^{(k)}_{m}  = \underset{k\text{\rm -times}}{\underbrace{[\bT, [\bT, [\bT,  \ldots [\bT, }}X_{m}] \ldots ]]]$ and its value $W_{\ell(a)j}\big|_{z_o}  \in \cD^{II}\big|_{z_o}$ projects onto 
the  homogeneous element $[w_{\ell(a)j}] \in \bV^\ell$.  The last property   implies that the  constant  matrix, given by  the coefficients which express $W_{\ell(a)j}$ in terms of an appropriate subset of  the set $\{X^{(\ell)}_m\}$ (namely of a subset  of the latter set which is  linearly independent at $z_o$) is invertible. Thus, we may replace the elements of such a  subset  by the vector fields  $W_{\ell(a)j}$ and, up to a re-ordering,  get  a new  collection of generators with the following properties: (a) the  first generators are  exactly the vector fields $W_{\ell(a)j}$; (b)
all other generators are vector fields  having the form $Y_r = X^{(k_r)}_{m_r} $ for some appropriate integers $k_r$ and $m_r$. With no loss of generality, we may assume that the cardinality of the  set $\{Y_r\}$ is the smallest possible.  If the set $\{Y_r\}$ is empty, we conclude that the  collection  $\big\{W_{\ell(a)j}\big\}$ is already  a set of  real analytic generators for $\cD^{II}|_{\cU}$ which satisfies (2) and (3) by  construction.  Claim (1)  is a direct consequence of the fact that the $W_{0(a)j}$ are linear combinations with constant coefficients of the $X_j$ and the lemma is proved in this case.\par
 If  the complementary set  $\{ Y_r\}$ is not empty,  $\big\{W_{\ell(a)j}\big\}$  is of course no longer  a  set of generators, but it  is still a set that satisfies (1).  Replacing the vector fields  $X_m$ by  linear combinations with constant coefficients of the generators $W_{0 (a) j}$ of $\cD^I|_{\cU}$,  we may assume that the  $Y_r$ have the form  
 $$ Y_r = \underset{(b_r)\text{\rm -times}}{\underbrace{[\bT, [\bT, [\bT,  \ldots [\bT, }}W_{0(a_r)j_r}] \ldots ]]] \qquad \text{for appropriate  integers}\  a_r, j_r \ \text{and}\  b_r\ .$$
Since the cardinality of  $\{Y_r\}$ is the smallest possible, each  $b_r$ is greater than the associated  integer $a_r$ (otherwise  $Y_r$ would coincide with one of the other generators $W_{\ell(a_r) j_{r}}$, $\ell \leq a_r$, and it could be omitted). Now,  for any  $0\leq a\leq \nu$,  we  define
  $$f[a]\= \max\bigg(\{ a\} \cup\{  b_r > a\  \text{for some of the}\ Y_r\ \text{for which}\ a_r = a\}\bigg)$$
 and  enlarge the  collection of vector fields $\{W_{\ell(a) j}\}$ into a new collection which includes  the vector fields 
 $$W_{\ell (f[a]) j }  \=  \underset{(\ell)\text{\rm -times}}{\underbrace{[\bT, [\bT, [\bT,  \ldots [\bT, }}W_{0(a)j}] \ldots ]]] \ ,\  \ \ a + 1 \leq \ell \leq f[a]\ .$$
 Such   enlarged collection  includes not only the previously constructed vector fields $W_{\ell(a) j}$, $\ell \leq a$,  but also  the 
 vector fields $Y_r$ and is therefore a collection of generators of $(V^{II}|_{\cU}, \cD^{II}|_{\cU})$ as desired. Replacing $\nu$ by $\nu' = \max\big(\{\nu\} \cup \{f[a], \ 0 \leq a \leq \nu\}\big)$ and  appropriately renaming  the triples $``\ell(a)j$'',  one can  check that such  enlarged set of vector fields is a set of generators  satisfying  (1), (2) and (3). 
\par
 \medskip
 \section{Proof of Lemma \ref{secondlemma}}
 \label{proofofsecondlemma}
 \begin{pf}  We prove  claim (i) by induction on $m$. For $m = 1$ and for any $s_o$,   we have that   $ \frac{d}{d s} \Phi^{X_1}_s(y)\bigg|_{s_o} =  X_1|_{\Phi^{X_1}_{s_o}(y)}$ by the  definition of a  flow.  This immediately implies  the first relation in \eqref{9.13}  and that 
$\frac{d^2}{d s^2} \left(\Phi^{X_1}_s(y) \right)\bigg|_{s = 0}  = \frac{\p X_1^j}{\p x^\ell}\bigg|_{y} X^\ell_1|_y = X_1(X_1)|_y$.
Assume now that (i) has been proved for $m -1$. Then, for any $s_o$ we have 
\begin{multline*} \frac{d}{d s} \left(\Phi^{X_1}_s  \circ \ldots \circ \Phi^{X_m}_s(y) \right)\bigg|_{s = s_o} \hskip - 0.5cm =  X_1|_{ \Phi^{X_1}_{s_o}  \circ \ldots \circ \Phi^{X_m}_{s_o}(y)  } 
 + \Phi^{X_1}_{s_o*} \left( \frac{d}{d s} \Phi^{X_2}_s  \circ \ldots \circ \Phi^{X_m}_s(y) \right)\bigg|_{s = s_o}\hskip - 0.3cm  \ . 
  \end{multline*}
  This and  inductive hypothesis yield 
  $ \frac{d}{d s} \left(\Phi^{X_1}_s  \circ \ldots \circ \Phi^{X_m}_s(y) \right)\bigg|_{s = 0}\hskip - 0.5 cm  = X_1|_y + \sum_{i = 2}^n X_i|_y$ and  
  \begin{multline*} \frac{d^2}{d s^2} \left(\Phi^{X_1}_s  \circ \ldots \circ \Phi^{X_m}_s(y) \right)\bigg|_{s = s_o} =  (X_1 + X_2 + \ldots + X_m)(X_1)|_y + \\
  +   (X_2 + \ldots + X_m)(X_1)|_y  + \\
  +  \sum_{i = 2}^{m} \left((X_i + X_{i + 1} + \ldots + X_m)( X_i)|_y  + (X_{i + 1} + \ldots + X_m)( X_i)|_y\right)\ , 
  \end{multline*}
  from which  claim (i)  follows immediately. Claim
  (ii) is an immediate consequence of (i). In particular, the second identity follows directly by applying \eqref{9.13} to the case $X_1 =Z +  Y$, $X_2 = -Z$ and  $X_3 = -Y$.  For  (iii), 
the claim on the first derivative is a direct consequence of \eqref{9.13}. For the second derivative,  by (ii) and  the definition  of $\cA$ 
\begin{multline}
 \frac{d^2}{d s^2} \left(\Phi^{X_1+ \ldots + X_m}_s \circ \Phi^{-Z}_s \circ \Phi^{-Y}_s(y)\right)\bigg|_{s = 0}  = \\
 =  \frac{d^2}{d s^2} \left(\Phi^{Y + Z + \cA}_s  \circ  \Phi^{- Y - Z }_s \right)   \circ \left( \Phi^{Y + Z }_s \circ \Phi^{-Z}_s \circ \Phi^{-Y}_s(y)\right)\bigg|_{s = 0}  =\\
 = \frac{d^2}{d s^2} \left(\Phi^{Y + Z + \cA}_s  \circ  \Phi^{- Y - Z }_s(y) \right) +  \frac{d^2}{d s^2} \left(\Phi^{ Y + Z}_s \circ \Phi^{-Z}_s \circ \Phi^{-Y}_s(y)\right)\bigg|_{s = 0}   =\\
 =\xcancel{ 2(Y + Z)(Y+ Z)|_y} + (Y+Z)(\cA)|_y + \cA(Y+ Z)|_y + \cA(\cA)|_y -\\
 -  \xcancel{2(Y+Z)(Y+Z)|_y} - 2 (Y+Z)(A)|_y +  [Y, Z]_y  = \\
 = [\cA, Y+ Z]_y + \cA(\cA)|_y +  [Y, Z]_y  \ .\end{multline} 
  On the other hand, by (i), 
  \begin{multline*}
 \frac{d^2}{d s^2} \left(\Phi^{X_1+ \ldots + X_m}_s \circ \Phi^{-Z}_s \circ \Phi^{-Y}_s(y)\right)\bigg|_{s = 0}  = \\
 =  \left(X_1+ \ldots + X_m\right)\left(X_1+ \ldots + X_m\right)|_y + Y(Y)  |_y+ Z(Z)  |_y
 + 2 Y(Z) |_y - \\
 - 2 Y( X_1+ \ldots + X_m) |_y  - 2 Z( X_1+ \ldots + X_m) |_y =\\
 \!\!=\!\! \sum_{i = 1}^m X_i(X_i)|_y {+}\!\!\!\!\! \sum_{\smallmatrix 1 \leq  \ell  \leq m\\ j \neq \ell \endsmallmatrix} X_\ell(X_j)|_y {+} Y(Y)|_y {+} Z(Z) |_y
 {+} 2 Y(Z)  |_y
 {-} 2 Y(\sum_{i = 1}^m X_i) |_y \! {-} 2 Z( \sum_{i = 1}^mX_i) |_y.\end{multline*} 
 This implies that 
 \begin{multline*}
 [Y, Z]_y +  [\cA, Y+ Z]_y + \cA(\cA)_y  - \frac{d^2}{d s^2} \left(\Phi^{X_1}_s \circ \ldots \circ \Phi^{X_m}_s \circ \Phi^{-Z}_s \circ \Phi^{-Y}_s(y)\right)\bigg|_{s = 0} = \\
 =  \frac{d^2}{d s^2} \left(\Phi^{X_1+ \ldots + X_m}_s \circ \Phi^{-Z}_s \circ \Phi^{-Y}_s(y)\right)\bigg|_{s = 0} \hskip - 5 pt  - \frac{d^2}{d s^2} \left(\Phi^{X_1}_s \circ \ldots \circ \Phi^{X_m}_s \circ \Phi^{-Z}_s \circ \Phi^{-Y}_s(y)\right)\bigg|_{s = 0}  \hskip - 5 pt = \\
  = 
   \xcancel{\sum_{i = 1}^m X_i(X_i)|_y }+  \sum_{\ell = 2}^m  \sum_{j = 1}^{\ell -1} X_\ell(X_j) |_y + \sum_{\ell = 1}^m  \sum_{j =  \ell+1}^m X_\ell(X_j) |_y +\xcancel{ Y(Y)|_y} + \xcancel{Z(Z)|_y} 
 +\xcancel{ 2 Y(Z)|_y} - \\- \xcancel{2 Y( \sum_{i = 1}^m X_i) |_y}  - \xcancel{2 Z( \sum_{i = 1}^m X_i) |_y}  -
\\
  - \xcancel{\sum_{i = 1}^m X_i(X_i)|_y}  - \xcancel{Y(Y)|_y} - \xcancel{Z(Z)|_y} - \xcancel{2 Y(Z)|_y} +\xcancel{ 2  \sum_{i = 1} Y(X_i)|_y} + \xcancel{ 2  \sum_{i = 1} Z(X_i)|_y} -
\\
   - 2 \sum_{\ell = 2}^m \sum_{j = 1}^{\ell-1} X_\ell(X_j)|_{y} 
  =   \sum_{\ell = 1}^m  \sum_{j =  \ell+1}^m X_\ell(X_j) |_y - \sum_{\ell = 2}^m \sum_{j = 1}^{\ell-1} X_\ell(X_j)|_{y}
\end{multline*}
and  claim (iii) follows.
\end{pf}

\bigskip
\bigskip
\font\smallsmc = cmcsc8
\font\smalltt = cmtt8
\font\smallit = cmti8
\hbox{\parindent=0pt\parskip=0pt
\vbox{\baselineskip 9.5 pt \hsize=3.7truein
\obeylines
{\smallsmc 
Marta Zoppello
Dipartimento di Scienze Matematiche 
``G. L. Lagrange'' (DISMA)
Politecnico di Torino
Corso Duca degli Abruzzi, 24, 
10129 Torino 
ITALY
}\medskip
{\smallit E-mail}\/: {\smalltt  marta.zoppello@polito.it 
\ 
}
}
\hskip 0.0truecm
\vbox{\baselineskip 9.5 pt \hsize=3.7truein
\obeylines
{\smallsmc
Cristina Giannotti \& Andrea Spiro
Scuola di Scienze e Tecnologie
Universit\`a di Camerino
Via Madonna delle Carceri
I-62032 Camerino (Macerata)
ITALY
\ 
}\medskip
{\smallit E-mail}\/: {\smalltt cristina.giannotti@unicam.it
\smallit E-mail}\/: {\smalltt andrea.spiro@unicam.it}
}
}
\end{document}